\documentclass[preprint,11pt]{elsarticle}
\usepackage[T1]{fontenc}
\usepackage{lmodern}
\usepackage{amsmath,amssymb,amsthm,mathtools,mathrsfs}
\usepackage{booktabs,array,longtable,enumitem}
\usepackage{geometry}
\usepackage{microtype}
\usepackage[hidelinks,unicode]{hyperref}
\hypersetup{pdfauthor={Chen Cao, Hebai Chen, Yilei Tang}}
\setlist{nosep,leftmargin=*}
\allowdisplaybreaks[2]
\newtheorem{theorem}{Theorem}[section]
\newtheorem{lemma}[theorem]{Lemma}
\newtheorem{proposition}[theorem]{Proposition}
\newtheorem{corollary}[theorem]{Corollary}
\theoremstyle{definition}

\theoremstyle{remark}
\newtheorem{remark}[theorem]{Remark}
\newcommand{\R}{\mathbb R}
\DeclareMathOperator{\sech}{sech}
\DeclareMathOperator{\artanh}{artanh}

\begin{document}
\begin{frontmatter}
\title{New criteria for upper bounds on the number of limit cycles in Li\'enard systems}
\author[sjtu]{Chen Cao}
\ead{aceaccce@sjtu.edu.cn}
\author[csu]{Hebai Chen}
\ead{chen_hebai@csu.edu.cn}
\author[sjtu]{Yilei Tang\corref{corresponding}}
\ead{mathtyl@sjtu.edu.cn}
\cortext[corresponding]{Corresponding author.}
\affiliation[sjtu]{organization={School of Mathematical Sciences, CMA-Shanghai, Shanghai Jiao Tong University},
 city={Shanghai}, postcode={200240}, country={PR China}}
\affiliation[csu]{organization={School of Mathematics and Statistics, Hunan Research Center of the Basic Discipline for Analytical Mathematics, HNP-LAMA, Central South University},
 city={Changsha}, postcode={410083}, country={PR China}}
\begin{abstract}
We establish criteria for at most one or at most two limit cycles of
generalized Li\'enard systems in a single potential well. Two independent
uniqueness tests use modified quotients of the damping and restoring
forces; both allow a degenerate restoring force, and one extends to
continuous potential profiles with corners. Opposite convexity of the
two branches of the damping quotient in the potential coordinate gives
an upper bound of two. Further criteria apply to weak damping with an
analytic normalized profile on a bounded potential interval under an
explicit smallness condition, to small profile perturbations vanishing
near the equilibrium on a bounded range of half-widths, and to parameter
families satisfying curvature, monotonicity and boundary conditions.
Applications give another proof of the classical quartic bound. For Selkov
systems, we obtain a two-cycle bound in a specified negative-parameter
region and prove existence and uniqueness in a positive-parameter region.
We also prove Conjecture~2.1 in [J. Math. Anal. Appl., 2025]. Examples attain
the two-cycle bound.
\end{abstract}
\end{frontmatter}

\section{Introduction}\label{sec:introduction}
The classical existence and uniqueness theory of Li\'enard and
Levinson--Smith~\cite{Lienard1928,LevinsonSmith1942} relates self-sustained
oscillations to the damping and restoring forces. For the equation
and its planar form
\[
 \ddot x+f(x)\dot x+g(x)=0,\qquad
 \dot x=y-F(x),\quad \dot y=-g(x),\quad F'=f,
\]
the counting problem is to find conditions on $f$ and $g$ that restrict
the number of isolated periodic orbits. We are concerned with upper
bounds of one and two in a single potential well.

Comparison of nested trajectories is a recurring method in uniqueness
proofs. The damping conditions of Sansone~\cite{Sansone1949}, the
homothetic comparison of Massera~\cite{Massera1954}, and its extension to
planar systems by Coppel~\cite{Coppel1964} give classical examples.
For generalized Li\'enard equations, this approach was developed by
Zhang~\cite{Zhang1986}; related uniqueness criteria are surveyed by Zeng, Zhang and
Gao~\cite{ZengZhangGao1994}. Later criteria allow asymmetry
\cite{CarlettiVillari2005}, formulate the problem on strips
\cite{XiaoZhang2008}, or use energy-level comparisons and weaker sign
restrictions on the damping~\cite{VillariZanolin2016,VillariZanolin2018}.

Uniqueness for cycles surrounding several equilibria has been studied
in~\cite{DumortierLi1996,ChenChen2015,ChenTangXiao2025}.
For nonsmooth equations, Li, Chen and Chen~\cite{LiChenChen2020} treated
crossing periodic orbits with one or several equilibria, while Chen,
Llibre and Tang~\cite{ChenLlibreTang2021} allowed specified discontinuities
in a strip criterion applied to Higgins--Selkov systems. Our assumptions
keep the restoring force directed towards one equilibrium. In the
interface version of our uniqueness criterion, $F$ and $g$ remain
continuous, while the potential profiles may have corners.

Upper bounds of two have been obtained under several different types of
structure. Rychkov's sharp bound concerns an odd quintic $F$ with
$g(x)=x$~\cite{Rychkov1975}. Two-cycle criteria and an example attaining
the bound were given in~\cite{GasullGiacominiLlibre2009}; symmetric systems
allowing discontinuities and an asymmetric three-zone piecewise linear
family were treated in~\cite{ChenHanXia} and~\cite{ChenTang2019},
respectively. Dulac functions give another way to control the number of
cycles~\cite{GasullGiacomini2002}. Combined with stability arguments, they
yield a sharp bound of two, counting multiplicities, for a one-parameter
Wilson Li\'enard-type family~\cite{GasullGiacomini2021}. Our two-cycle
criterion instead imposes opposite convexity on the two branches of
$f/g$ as functions of the potential. It applies without a prescribed
polynomial degree or parity, but requires these shape conditions
throughout the well under consideration.

The shape of $F$ also enters criteria for larger cycle counts: Odani
\cite{Odani1996} obtained exactly $N$ hyperbolic periodic orbits under
suitable shape conditions when $g(x)=x$, and Chen, Llibre and Zhang
\cite{ChenLlibreZhang2007} gave conditions for at least $n$ or exactly $n$
limit cycles. Local cycle counts and small-amplitude constructions are
studied in~\cite{BlowsLloyd1984,ChristopherLynch1999,YuHan2006,HanTianYu2013}.
Higher-order averaging~\cite{LlibreMereuTeixeira2010} and the constructions
of Han and Romanovski~\cite{HanRomanovski2013} yield further lower bounds
for polynomial families. An upper bound throughout a potential well must
also exclude additional cycles away from the bifurcation region.

For the classical polynomial system $\dot x=y-F(x)$, $\dot y=-x$, the
degree convention is $n=\deg F$.
Lins, de Melo and Pugh~\cite{LinsDeMeloPugh1977} proposed the bound
$\lfloor(n-1)/2\rfloor$ and established the low-degree cases up to three.
Li and Llibre~\cite{LiLlibre2012} proved the quartic case: there is at most
one limit cycle, and it is hyperbolic if it exists. The application in
Section~\ref{ap:quartic-section} gives another proof of this result.
For degree five,
Li and Lu~\cite{LiLu2014} obtained cyclicity at most two near nondegenerate
slow--fast cycles.
The proposed degree bound fails in higher degrees. Dumortier, Panazzolo
and Roussarie~\cite{DumortierPanazzoloRoussarie2007} constructed examples
with at least $\lfloor(n-1)/2\rfloor+1$ cycles for $n\ge7$.
De Maesschalck and Dumortier~\cite{DeMaesschalckDumortier2011} obtained at
least $\lfloor(n-1)/2\rfloor+2$ cycles for $n\ge6$, and De Maesschalck
and Huzak~\cite{DeMaesschalckHuzak2015} constructed at least $n-2$
hyperbolic cycles for every $n\ge6$.

Write $G(x)=\int_0^x g(s)\,ds$. Our first uniqueness criterion tests the
strict monotonicity of $f/g-KG$ on each side of the equilibrium, with the same
constant $K$ on both sides. The choice $K=0$ tests $f/g$ itself; allowing
$K\ne0$ admits quotients that change their direction of monotonicity.
For smooth coefficients with $F'(0)<0$, the second criterion requires
the derivative of $(f/g-d)/G$ to be positive on both sides for a common
constant $d$; when $F'(0)>0$, the derivative must be negative.
The two conditions are independent, as the examples in
Section~\ref{sec:applications} show. In the proofs, the coordinate $u=G(x)$
turns each half-orbit into a scalar equation. We compare its endpoints
with those of quadratic profiles to control the number of matched orbits.

The potential intervals may be finite and unequal; the two basic
uniqueness criteria also allow a degenerate restoring force at the
equilibrium. Joint comparison of the two profiles gives further
two-cycle criteria when the separate third-order sign conditions fail.
Weakly nonlinear approximations were studied by L\'opez and
L\'opez-Ruiz~\cite{LopezLopezRuiz2009}, and analytic properties of the return
map by Fran\c{c}oise~\cite{Francoise2002}. Our weak-damping argument estimates
the nonlinear remainder in the difference of the fitted curvatures.
It assumes an analytic normalized damping profile on a bounded potential
interval and an explicit smallness condition. A separate criterion treats
$C^2$-small profile perturbations that vanish near the equilibrium; its
two-cycle bound holds on a prescribed bounded range of half-widths.
For the Selkov application, we combine branch
curvature with parameter monotonicity and endpoint conditions to control
the whole family of periodic orbits.

We study the generalized Li\'enard system
\begin{equation}\label{cr:system}
 \dot x=y-F(x),\qquad \dot y=-g(x),\qquad (x,y)\in S:=I\times\mathbb R,
\end{equation}
where $I=(a,b)$, $-\infty\le a<0<b\le\infty$. Unless stated otherwise,
our basic assumptions are
\begin{equation}\label{cr:basic}
 F,g\in C(I),\quad F\in C^1(I\setminus\{0\}),\quad
 F(0)=g(0)=0,\quad xg(x)>0\quad(x\ne0).
\end{equation}
Write $I_-=(a,0)$, $I_+=(0,b)$, $f=F'$ on $I_-\cup I_+$, and
\begin{equation}\label{cr:potential}
 G(x)=\int_0^xg(s)\,ds,\qquad
 U_- =\lim_{x\to a^+}G(x),\qquad U_+=\lim_{x\to b^-}G(x).
\end{equation}
At an infinite endpoint, the corresponding limit is taken as $x\to-\infty$
or $x\to+\infty$, respectively.
Both $U_\pm$ belong to $(0,\infty]$; they need not be equal or infinite.
Let $x_\pm:(0,U_\pm)\to I_\pm$ be the inverse branches of $G$, and put
\begin{equation}\label{cr:profiles}
 \phi_\pm(u)=F(x_\pm(u)),\qquad
 k_\pm(u)=\phi_\pm'(u)=\frac{f(x_\pm(u))}{g(x_\pm(u))}.
\end{equation}
The profiles extend continuously to $u=0$ with value zero. Throughout,
a periodic orbit is understood to be nonconstant. A limit cycle
is an isolated periodic orbit, counted once as a geometric
curve. Nonhyperbolic isolated orbits, including semistable cycles, are included;
nonisolated orbits in a period annulus are not. All assertions concern orbits
entirely contained in $S$. No invariance or completeness assumption on $S$ is
made. Under the shape conditions below, the transverse return map is $C^1$
even when the original vector field is only continuous at the origin.
Hyperbolicity refers to its multiplier being different from one;
a cycle is neutral when this multiplier equals one.

\begin{theorem}\label{cr:unique}
Assume \eqref{cr:basic}. If one constant $K\in\mathbb R$ makes
\[
 \frac{f(x)}{g(x)}-KG(x)
\]
strictly increasing on each of $I_-$ and $I_+$, then
\eqref{cr:system} has at most one nonconstant periodic orbit in $S$.
If it exists, it is a hyperbolic attracting limit cycle.
If the two functions are strictly decreasing, the same upper bound holds,
and the orbit, if present, is hyperbolic and repelling.
The same $K$ is required on the two sides; no comparison across $x=0$ is intended.
\end{theorem}

\begin{corollary}\label{cr:unique-derivative}
Under \eqref{cr:basic}, suppose additionally that $f,g\in
AC_{\mathrm{loc}}(I_-\cup I_+)$. Define almost everywhere
\begin{equation}\label{cr:U}
 \mathcal U_K=f'g-fg'-Kg^3.
\end{equation}
The attracting conclusion of Theorem~\ref{cr:unique} follows if
$\mathcal U_K>0$ almost everywhere. More generally, it suffices that
$\mathcal U_K\ge0$ almost everywhere and
$\int_c^d\mathcal U_K/g^2\,dx>0$ on every nondegenerate compact interval
$[c,d]$ lying on one side of zero. Reversing both inequalities gives the
repelling conclusion. In particular, $F\in C^2$, $g\in C^1$ off zero and
$\mathcal U_K>0$ pointwise suffice.
\end{corollary}

For an upper bound of two, we impose opposite convexity conditions on the
two slope functions $k_+$ and $k_-$. The following formulation uses secant
slopes and therefore also covers profiles without a classical third derivative.

\begin{theorem}\label{cr:two}
Assume \eqref{cr:basic}, and put $q=f/g$. Suppose that, for every
$x_1<x_2<x_3$ in the same one of $I_-,I_+$,
\begin{equation}\label{cr:secants}
 \frac{q(x_2)-q(x_1)}{G(x_2)-G(x_1)}
 \le \frac{q(x_3)-q(x_2)}{G(x_3)-G(x_2)}.
\end{equation}
Then \eqref{cr:system} has at most two limit cycles in $S$.
The same conclusion holds if all the inequalities in \eqref{cr:secants}
are reversed. Equivalently, $k_+$ is convex and $k_-$ is concave, or vice versa.
Convexity and concavity need not be strict, and affine subintervals are allowed.
\end{theorem}

\begin{corollary}\label{cr:two-derivative}
Under \eqref{cr:basic}, suppose that $F\in C^2$, $g\in C^1$ off zero and
$F'',g'\in AC_{\mathrm{loc}}(I_-\cup I_+)$. If the almost-everywhere expression
\begin{equation}\label{cr:N}
 \mathcal N=g^2f''-3gg'f'+\bigl(3(g')^2-gg''\bigr)f
\end{equation}
is nonnegative almost everywhere, or nonpositive almost everywhere, then there are at most two limit cycles in $S$.
In particular, the pointwise test applies when $F\in C^3$, $g\in C^2$
off zero. No derivative at zero is required.
\end{corollary}

\begin{theorem}\label{cr:intercept}
Assume \eqref{cr:basic}, $F,g\in C^\infty(I)$ and $F'(0)<0$.
If one constant $d\in\mathbb R$ satisfies
\begin{equation}\label{cr:intercept-condition}
 \frac{d}{dx}\left(\frac{f(x)/g(x)-d}{G(x)}\right)>0
 \qquad(x\in I_-\cup I_+),
\end{equation}
then there is at most one nonconstant periodic orbit in $S$; any such orbit
is a hyperbolic attracting limit cycle. Equivalently, one may check
\begin{equation}\label{cr:I}
 \mathcal I_d:=G(f'g-fg')-fg^2+dg^3>0\qquad(x\ne0).
\end{equation}
The restoring force may be degenerate: $g'(0)>0$ is not assumed.
If $F'(0)>0$ and the inequality in \eqref{cr:intercept-condition} is
reversed, the corresponding orbit is hyperbolic and repelling.
\end{theorem}

The intercept criterion requires smooth coefficients, a strict inequality,
and $F'(0)\ne0$. For common curvature, the convexity formulation below
allows corners in the potential profiles.

\begin{theorem}\label{cr:interface}
Assume only $F,g\in C(I)$, $F(0)=g(0)=0$, and $xg(x)>0$ for $x\ne0$.
Define $G$, $x_\pm$ and $\phi_\pm$ as above. If one constant $K$ makes
\begin{equation}\label{cr:interface-condition}
 \phi_+(u)-\frac K2u^2\ \text{strictly convex on }(0,U_+),\qquad
 \phi_-(u)-\frac K2u^2\ \text{strictly concave on }(0,U_-),
\end{equation}
then there is at most one nonconstant periodic orbit in $S$.
Any such orbit has a $C^1$ transverse return map with multiplier in $(0,1)$,
and is therefore a hyperbolic attracting limit cycle.
Reversing convexity and concavity gives a hyperbolic repelling orbit.
The strict convexity is meant in the usual chord sense and does not
require the first derivative to exist everywhere.
\end{theorem}

\begin{corollary}\label{cr:interface-piecewise}
Under the continuity and restoring-force assumptions of
Theorem~\ref{cr:interface}, suppose $F$ is locally Lipschitz off zero,
$C^1$ between a locally finite set of interfaces there, and has finite
one-sided first derivatives at every interface. If $q-KG$ is strictly
increasing on each open piece, for the same $K$, and
\[
 q(c+)-q(c-)=\frac{F'(c+)-F'(c-)}{g(c)}\ge0
\]
at every interface $c\ne0$, the attracting conclusion holds.
Strict decrease on the pieces and nonpositive jumps give the repelling
conclusion. Interfaces may accumulate only at zero or at a strip boundary.
\end{corollary}

The interface extension is specific to common-curvature uniqueness:
allowing slope jumps of one sign does not extend Theorem~\ref{cr:two}.
Existence will be treated separately in the applications, where several
examples attain the stated upper bounds.

The proofs use quadratic profiles as comparison models. Their global
fitting property and the required shape inequality are established in
Section~\ref{cr:sec-two}.

\paragraph{Joint criteria}
The preceding third-order test treats the two potential profiles separately.
A joint comparison can give a two-cycle bound even when both third
derivatives change sign. The following theorem gives an explicit
weak-damping range for this comparison.

\begin{theorem}\label{jt:weak}
Let $I=(\alpha,\beta)$, $F,g\in C^\infty(I)$,
$F(0)=g(0)=0$, $xg(x)>0$ for $x\ne0$, and $g'(0)>0$. Set
\[
 G(x)=\int_0^x g(s)\,ds,\qquad
 z(x)=\operatorname{sgn}(x)\sqrt{2G(x)}.
\]
Suppose $z(I)=(-a,b)$, where $0<a,b\le3$, and
$F(x)=\varepsilon p(z(x))$, where $0<\varepsilon\le10^{-210}$.
Assume that $p$ is holomorphic on $|w|<10$, is real on the real axis,
$p(0)=0$, and
\[
 \sup_{|w|<10}|p(w)/w|\le10^{40},\qquad
 wp(w)>0\quad(0<|w|\le1/16, w\in\mathbb R).
\]
The quotient at zero is interpreted by continuity. Define
\begin{align}
 D_p(r)&=\frac12\int_0^\pi\sin\theta
       [p(r\sin\theta)-p(-r\sin\theta)]\,d\theta,\label{jt:D}\\
 C_p(r)&=\frac{15}{2r^3}\left(D_p'(r)-\frac{2D_p(r)}r\right).
       \label{jt:C}
\end{align}
If $C_p'(r)\ge1$ throughout $[1/16,4]$, or $C_p'(r)\le-1$
throughout this interval, then the system has at most two limit cycles in
$I\times\mathbb R$. Nonhyperbolic isolated cycles are included, each counted
once. No equality of the potential endpoints or symmetry of $p$ is assumed.
\end{theorem}

For perturbations supported away from the equilibrium, the two-cycle bound
persists on a bounded interval of half-widths under a sufficiently small
perturbation in the $C^2$ norm. The required smallness depends on the
reference system and the half-width interval.
For a profile $\phi$, let $M_\phi(r)$ denote the midpoint of the two
zero-height endpoints of its transverse half arc of half-width $r$ for
$u_y=\phi(u)-y$. For the quadratic profile $\lambda u+\kappa u^2/2$ write
this midpoint as $H(\lambda,\kappa,r)$, and let
$\mathcal K(M,v,r)$ be the curvature component of the inverse fit
$(H,H_r)=(M,v)$. These quantities are uniquely defined for
$r>0$, $|M|<r$, $|v|<1$ by Lemma~\ref{cr:fit}. Put
\[
 \mathscr T(M,v,w,r)=\mathcal K_Mv+\mathcal K_vw+\mathcal K_r.
\]

\begin{theorem}\label{jt:compact}
Consider a smooth system $\dot x=y-F(x)$, $\dot y=-g(x)$ on an open interval
containing zero, with $F(0)=g(0)=0$ and $xg(x)>0$ off zero.
Write its potential profiles as $\phi_\pm=\phi^0_\pm+\psi_\pm$.
Assume that the reference profiles are continuous at zero, smooth at positive
height, vanish at zero, and have transverse half arcs at all half-widths
$0<r\le R_1$, where $R_1>R>0$. Assume
\[
 (\phi^0_+)'''>0,\quad (\phi^0_-)'''<0,\qquad
 (\phi^0_\pm)'\in L^1(0,\eta)
\]
for some $\eta>0$. Fix $a>0$ and
$0<r_0<\min(R,\sqrt{2a})$. Suppose $\psi_\pm=0$ on $[0,a]$.
Choose heights $T_\pm>a$ above all reference peaks on a neighborhood of
$[r_0,R]$, within common height domains of reference and candidate profiles.
Let $C_\pm,\epsilon_\pm>0$ be second-jet stability constants from
Lemma~\ref{jt:jets} for these data. Around the reference triples
$(M^0_\pm,(M^0_\pm)',(M^0_\pm)'')$ choose compact tubes, at each fixed
$r$, convex and contained in $|M|<r$, $|v|<1$. Let $L_\pm$ bound
$|\mathscr T_M|+|\mathscr T_v|+|\mathscr T_w|$ on the respective tubes.
Choose $m>0$ such that, on $[r_0,R]$,
\[
 \mathscr T(M^0_+,(M^0_+)',(M^0_+)'',r)
 -\mathscr T(M^0_-,(M^0_-)',(M^0_-)'',r)\ge m.
\]
If $\delta_\pm=\|\psi_\pm\|_{C^2([a/2,T_\pm])}<\epsilon_\pm$,
$C_\pm\delta_\pm$ are smaller than the tube clearances, and
\[
 L_+C_+\delta_++L_-C_-\delta_-<m,
\]
then at most two limit cycles have half-width in $(0,R]$. This is a bound
for all cycles in the strip only if an independent argument places all such
cycles in this half-width interval.
\end{theorem}

Theorem~\ref{jt:compact} gives a stability estimate with reference-dependent
constants, whereas Theorem~\ref{jt:weak} uses analyticity to obtain a fixed
explicit bound.

We next combine curvature conditions on the two branches of
$f/g$ with a parameter-family argument. The curvature conditions determine
the direction of a neutral return map. To deduce an upper bound for all
cycles, we also control how periodic orbits are connected as the parameter
varies. For the parameter-family results, write
\begin{equation}\label{sk:system}
 \dot x=y-F(x,\theta),\qquad \dot y=-g(x),\qquad x\in J=(\ell,\rho),
 \quad \ell<0<\rho,
\end{equation}
where either endpoint of $J$ may be infinite, and write
\begin{equation}\label{sk:coefficients}
 f=F_x,\quad G(x)=\int_0^xg(s)\,ds,\quad
 L=\frac fg,\quad A=L_h=\frac{L_x}{g},\quad
 \Xi=LL_{hh}-L_h^2,
 \qquad h=G(x).
\end{equation}
The derivatives in $h$ are taken separately on the two inverse branches
of $G$, away from $x=0$.

\begin{theorem}\label{sk:two}
Let $F,g$ be $C^7$ on $J\times[\theta_-,\theta_H]$, where
$\theta_-<\theta_H$, with $g$ independent
of $\theta$ and with one-sided $C^7$ extensions at the parameter endpoints.
Assume the following conditions.
\begin{enumerate}
\item[(S1)] $F(0,\theta)=g(0)=0$, $g_1:=g'(0)>0$, and $xg(x)>0$ for
 $x\ne0$. Both $G(x)$ and $F(x,\theta)$ tend to $+\infty$ at both
 endpoints of $J$, for every parameter. Moreover $f(0,\theta)^2<4g_1$
 throughout the closed parameter interval.
\item[(S2)] For $\theta<\theta_H$, $f(0,\theta)<0$ and there is
 $x_f(\theta)>0$ such that $f<0$ on $(\ell,x_f)$ and $f>0$ on
 $(x_f,\rho)$, with $f(x_f,\theta)=0$. At $\theta_H$,
 $f(0,\theta_H)=0$, $f_x(0,\theta_H)>0$, and $xf(x,\theta_H)>0$
 for $x\ne0$. Throughout the closed parameter interval,
 \[
 A>0\quad(x>0),\qquad \Xi>0\quad(x<0),\qquad
 \Xi<0\quad(x>0,\ f>0).
 \]
\item[(S3)] $xF_\theta(x,\theta)>0$ for $x\ne0$. If $x_L(h)<0<x_R(h)$
 are the inverse branches of $G$, then
 \[
 F(x_R(h),\theta_-)\le F(x_L(h),\theta_-)\qquad(h>0).
 \]
\item[(S4)] The removable extension $L_H=f(\cdot,\theta_H)/g$ satisfies
 $L_{H,x}(0)>0$ and $(F(\cdot,\theta_H)L_H)_x<0$ on $(\ell,0)$.
\end{enumerate}
Then every system in the closed family has at most two periodic orbits
in $J\times\mathbb R$, and has no period annulus. Every nonhyperbolic
periodic orbit at an interior parameter is a limit cycle of multiplicity
two. At most one interior parameter admits such a cycle. If two cycles
coexist, the inner cycle is hyperbolic and attracting and the outer cycle
is hyperbolic and repelling. The lower endpoint has no periodic orbit;
the upper endpoint has at most one, which is hyperbolic and repelling
if it exists.
\end{theorem}

For a fixed system, the curvature conditions in Theorem~\ref{sk:two}
determine the direction of its neutral return maps.

\begin{proposition}\label{sk:neutral}
For a fixed system $\dot x=y-F(x)$, $\dot y=-g(x)$ on $J$, assume
$F(0)=g(0)=0$, $f=F'$, $f,g\in C^6$, $g'(0)>0$, and $xg(x)>0$
for $x\ne0$. Suppose that some $x_f>0$ satisfies
$f<0$ on $(\ell,x_f)$, $f(x_f)=0$, and $f>0$ on $(x_f,\rho)$,
and that
\[
 A>0\ (x>0),\qquad \Xi>0\ (x<0),\qquad \Xi<0\ (x>x_f).
\]
For the first return map $\Pi$ on $x=0$, $y=q>0$, every
periodic orbit satisfies
\[
 \Pi(q)=q,\quad\Pi'(q)=1\quad\Longrightarrow\quad\Pi''(q)>0.
\]
This proposition alone does not assert an upper bound for the number
of hyperbolic cycles of a fixed system.
\end{proposition}

We obtain connectedness of the periodic-orbit family by imposing a
condition at the parameter boundary, which
need not be a Hopf boundary. For each parameter and each half-plane, consider
the orbit arc whose maximum value of $G$ is $t>0$. Under (S1) this arc
meets $x=0$ at $y=A_i(t,\theta)<0<B_i(t,\theta)$, and its half-width
$(B_i-A_i)/2$ increases bijectively from zero to infinity. Let
$M_i(r,\theta)$ be its midpoint when its half-width is $r$, and put
$\Delta=M_R-M_L$. These constructions are justified in
Section~\ref{sk:family-proof}.

\begin{theorem}\label{sk:family}
Let $F,g$ be $C^7$ on $J\times[\theta_-,\theta_+]$, where
$\theta_-<\theta_+$, with $g$ independent of the parameter and
one-sided $C^7$ extensions at the parameter endpoints.
Assume (S1) and (S3) on $[\theta_-,\theta_+]$, with
$f(0,\theta_-)<0$. Assume also that
\[
 \{r>0:\Delta(r,\theta_+)>0\}
\]
is empty or an interval, and that every neutral periodic orbit at an
interior parameter has $\Pi''>0$ on the positive section. Then every
interior system has at most two periodic orbits. All its neutral cycles
have multiplicity two. All interior periodic orbits form one $C^2$
graph $\theta=\beta(r)$ over that interval; this graph has at most one
stationary point, necessarily a nondegenerate minimum. If two cycles
coexist, the inner one is hyperbolic and attracting and the outer one
is hyperbolic and repelling.
\end{theorem}

The boundary condition and the neutral-direction condition are the two
essential hypotheses of Theorem~\ref{sk:family} in addition to parameter
monotonicity. Theorem~\ref{sk:two} gives coefficient conditions that imply
both. The Hopf-boundary argument also yields the following uniqueness
criterion.

\begin{theorem}\label{sk:hopf}
Let $F,g\in C^7(J)$, $F(0)=g(0)=0$, $g'(0)=g_1>0$,
$xg(x)>0$ for $x\ne0$, and let $F,G\to+\infty$ at both endpoints.
Suppose $f=F'$ satisfies $f(0)=0$, $f'(0)>0$, and $xf(x)>0$ for
$x\ne0$. Extend $L=f/g$ removably to zero and suppose $L_x(0)>0$.
Using the left and right inverse branches of $F$, regard $L_L,L_R$
as functions of $\omega=F(x)>0$. If
\[
 (L_R-L_L)_{\omega\omega}<0,\qquad
 (\omega L_L)_\omega>0\qquad(\omega>0),
\]
then the system has at most one periodic orbit in $J\times\mathbb R$,
and that orbit is hyperbolic and repelling if it exists.
\end{theorem}

The proofs combine the branch curvature conditions with an area comparison
along the orbit arcs. In the Selkov application, the coefficient conditions
can be verified explicitly, and the two-cycle bound is attained on a
nonempty parameter region.

Sections~\ref{cr:sec-unique}--\ref{cr:sec-interface} develop the comparison
of half-orbits and prove the first four criteria.
Sections~\ref{jt:proof-section}--\ref{sk:coefficient-proof} treat the joint
estimates and the parameter-family arguments.
Section~\ref{sec:applications} applies the criteria to quartic systems,
a negative-parameter Selkov region, Conjecture~2.1 of
Chen--Li--Xiao--Yang, and further examples. It also gives a separate
existence and uniqueness proof for a positive-parameter Selkov region.

\section{Proof of the common-curvature criterion}\label{cr:sec-unique}
The proof of Theorem~\ref{cr:unique} compares the two half-orbits that
make up a periodic orbit. Their half-widths and midpoints provide coordinates
for matching them across the axis. We first establish these coordinates
for profiles that are only $C^1$ at positive height, then compare them
with smooth quadratic profiles.

\begin{lemma}\label{cr:flow}
Under \eqref{cr:basic}, the initial-value problem is locally unique away
from $(0,0)$, with continuous dependence on initial data and $C^1$
geometric transverse maps. If $xF(x)\ge-Lx^2$ near zero for some $L\ge0$,
the equilibrium is forward unique; if $xF(x)\le Lx^2$, it is backward
unique. Either condition excludes a nonconstant periodic solution through
$(0,0)$. The shape assumptions of either Theorem~\ref{cr:unique} or
Theorem~\ref{cr:two} imply one of these bounds.
\end{lemma}
\begin{proof}
For the first bound set $H=F+Lx$, $w=y+Lx$ and
\[
 V=G(x)+L\int_0^xH(s)\,ds+\tfrac12w^2.
\]
Here $xH(x)\ge0$, so $V\ge0$, with equality only at the origin. Along any
solution,
\[
 \dot V=-gH-LH^2+Lw^2\le2LV.
\]
Gronwall's inequality proves forward uniqueness at the equilibrium.
The transformation $(t,y)\mapsto(-t,-y)$ replaces $F$ by $-F$ and proves
the second assertion. One-sided uniqueness is enough for a periodic
solution through the equilibrium to be constant.

If $q-KG$ increases on both sides, its values have a finite upper bound
on the right near zero and a finite lower bound on the left. Since $g$
is bounded and has the sign of $x$, it follows that $f\le L$.
Integrating on either side and using $F(0)=0$ gives $xF(x)\le Lx^2$.
Decreasing functions give the opposite bound. Under \eqref{cr:secants},
$k_+$ is convex and $k_-$ concave. A finite convex function has a finite
lower affine bound near its left endpoint: for $0<u<u_0<u_1$,
\[
 k(u)\ge k(u_0)+(u-u_0)\frac{k(u_1)-k(u_0)}{u_1-u_0}.
\]
Thus $k_+$ is bounded below and $k_-$ above near zero. Multiplication by
$g$ gives $f\ge-L$ and hence the other energy bound. All directions reverse
under the reversed shape assumptions.

For $x\ne0$, $u=G(x)$ is a $C^1$ local diffeomorphism and $y$ is a valid
orbit coordinate. The equation becomes
\begin{equation}\label{cr:profile-ode}
 u_y=\phi(u)-y,
\end{equation}
a scalar $C^1$ equation. Its solutions depend $C^1$ on the initial value;
original time is recovered uniquely from
$dt/dy=-1/g(x(u(y)))$, whose denominator is bounded away from zero locally.
At an axis point $(0,y_0)$ with $y_0\ne0$, use instead
\begin{equation}\label{cr:axis-graph}
 y_x=\frac{g(x)}{F(x)-y},\qquad
 t-t_0=\int_{x_0}^x\frac{d\xi}{y(\xi)-F(\xi)}.
\end{equation}
In a small rectangle the denominators stay away from zero. The first
right-hand side is continuous in $x$ and smooth in $y$, with all $y$
derivatives locally bounded. Its integral equation proves uniqueness,
continuous dependence, and smooth dependence on the initial ordinate.
The time integral is finite and has a nonzero continuous derivative.
These local orbit coordinates give $C^1$ geometric flow boxes; replacing
the orbit coordinate by the strictly monotone continuous time integral
gives topological flow boxes. Finitely many of them cover a periodic
orbit away from the equilibrium and produce its $C^1$ return map.
The $C^1$ regularity obtained here concerns the geometric transverse maps.
\end{proof}

\begin{lemma}\label{cr:arcs}
Let $\phi\in C^1((0,U))\cap C([0,U))$, $\phi(0)=0$, $0<U\le\infty$.
For $0<h<t<U$, the solution of \eqref{cr:profile-ode} with
$u(\phi(t);t)=t$ meets height $h$ at two points $A(h,t)<B(h,t)$.
Define
\begin{equation}\label{cr:arc-data}
 r=\frac{B-A}{2},\qquad M=\frac{B+A}{2}-\phi(h),\qquad
 v=M_r,\qquad q_0=r^2-M^2.
\end{equation}
Then $A<\phi(h)<B$, $A_t<0<B_t$, and $r_t>0$.
The map $(h,t)\mapsto(h,r)$ is a $C^1$ diffeomorphism onto its actual open
image $\Omega_\phi=\{0<h<U,\ 0<r<R_\phi(h)\}$, and
\begin{equation}\label{cr:arc-identities}
 \begin{gathered}
 |M|<r,\quad |v|<1,\quad r^2\ge2(t-h),\qquad
 v=\frac{B_t+A_t}{B_t-A_t},\\
 A_h=\frac1{r-M},\quad B_h=-\frac1{r+M},\quad
 r_h\big|_t=-\frac r{q_0},\quad M_h\big|_t=\frac M{q_0}-\phi'(h).
 \end{gathered}
\end{equation}
There are finite limits $A_0(t)\le0\le B_0(t)$ as $h\to0^+$.
The set $\{t:A_0(t)<0<B_0(t)\}$ is empty or an interval $(\tau_\phi,U)$;
its half-width image $I_\phi$ is an open interval. On $I_\phi$, $M(0,r)$ is
$C^1$, $|M|<r$, $|M_r|<1$, and $(M(h,r),M_r(h,r))$ converges locally
uniformly to its axis value. For smooth positive-height profiles all
positive-height data are smooth, and this axis convergence holds for each
fixed derivative order with respect to $r$. The endpoints of $I_\phi$ may be finite.
\end{lemma}
\begin{proof}
Along an orbit $u_{yy}=\phi'(u)u_y-1$; every critical point is a
nondegenerate maximum. Both branches decrease away from the prescribed
peak and cannot encounter a second critical point. Boundedness of $\phi$
on $[h,t]$ forces them to reach $h$ in finite $y$; there can be no blow-up
while the height remains in this compact interval. The variation
$\xi=u_t$ satisfies
\[
 \xi(\phi(t);t)=1,\qquad \xi_y=\phi'(u)\xi,
\]
where the moving initial coordinate contributes nothing because $u_y=0$
at the peak. Hence $\xi>0$ and
$A_t=-\xi(A)/(r-M)<0$, $B_t=\xi(B)/(r+M)>0$.
Implicit differentiation at the transverse endpoints gives
\eqref{cr:arc-identities}, including $|v|<1$.
Integrating $d(r^2)/dh=-2r^2/q_0\le-2$ from the peak proves the half-width bound.
As $t\to h^+$, $r\to0^+$; strict increase in $t$ and the inverse
function theorem give the stated actual open image. For a profile smooth
on all of $\mathbb R$, these assertions also hold for any real $h<t$;
the half-width bound shows that either $t\to\infty$ at fixed $h$, or
$h\to-\infty$ at fixed $t$, gives all positive half-widths.

For fixed $t$, put $C=\max_{[0,t]}|\phi|$. The endpoints move outwards
as $h\to0^+$ and satisfy
$-C-\sqrt{2t}\le A<B\le C+\sqrt{2t}$: outside $[-C,C]$, integrate
$u_y\le C-y$ or its left-hand counterpart. Thus the limits exist.
Monotonicity in $t$ makes the transverse peak set upward closed.
At a nonzero limiting endpoint, choose a rectangle
$[0,\delta]\times(y_0-\eta,y_0+\eta)$ on which $|\phi(u)-Y|\ge m>0$
and solve the inverse equation
\[
 Y_u=\frac1{\phi(u)-Y}.
\]
Choose $\delta/m<\eta/2$ and initial values in the inner half rectangle;
the solution remains there to $u=0$. Its peak variation is
\[
 Y_t(u,t)=Y_t(\delta,t)
 \exp\!\left(\int_\delta^u\frac{dv}{(\phi(v)-Y(v,t))^2}\right),
\]
with exponential between $e^{-\delta/m^2}$ and $1$.
The strict endpoint-derivative signs therefore persist at the axis.
Transversality is open, proving the asserted interval description.
The inequality $r_t>0$ supplies a $C^1$ inverse peak coordinate.
On a compact interval of half-widths at the axis its peak image is compact below $U$;
a finite covering by these rectangles gives a uniform positive lower
bound for $r_t$ and uniform convergence of endpoints and variations.
This proves the stated convergence of $M,v$. In the smooth case, higher
variations of the inverse equation have bounded continuous coefficients
and inductively controlled inhomogeneous terms, so the same argument
works at every fixed derivative order with respect to $r$. Continuity of $\phi$ at
zero suffices for this argument.
\end{proof}

\begin{lemma}\label{cr:angles}
Along a fixed-peak half-orbit, parameterize outwards by its half-width $s$ and put
$k=\phi'$, $\theta=\artanh(M/s)$ and $L=\artanh v$. Then
\begin{equation}\label{cr:angle-ode}
 h'=-s\sech^2\theta,\qquad
 \theta'=k(h)-\frac{\sinh(2\theta)}s,\qquad
 L'=\frac{\sinh(2\theta)}s.
\end{equation}
At a peak of height $t$,
\begin{equation}\label{cr:angle-start}
 \theta(s)=\frac{k(t)}3s+o(s),\quad
 L(s)=\frac{2k(t)}3s+o(s),\qquad
 L(s)=\int_0^s\frac{\sinh(2\theta(\sigma))}{\sigma}\,d\sigma.
\end{equation}
The integral holds up to a transverse axis endpoint. For smooth profiles,
$M(h,r)=k(h)r^2/3+O(r^3)$ with a smooth remainder, including smooth
finite-dimensional profile parameters.
\end{lemma}
\begin{proof}
Set $y=\phi(t)+\eta$ at the peak. Successive substitution in
$u_{yy}=k(u)u_y-1$, using only continuity of $k$, gives
\[
 u=t-\eta^2/2-k(t)\eta^3/6+o(\eta^3).
\]
With $\rho=\sqrt{2(t-h)}$, its inverse branches satisfy
$\eta_\pm=\pm\rho-k(t)\rho^2/6+o(\rho^2)$.
Consequently $r=\rho+o(\rho^2)$ and $M=k(t)\rho^2/3+o(\rho^2)$,
which proves the first angular initial term.
For a smooth profile write $t-u=\eta^2\mathcal A(t,\eta)/2$ with smooth
$\mathcal A(t,0)=1$; invert $\eta\sqrt{\mathcal A(t,\eta)}$.
This constructs smooth remainders and their parameter derivatives.
The first two differential equations follow from \eqref{cr:arc-identities}.
To obtain the third without assuming second parameter derivatives, write
\[
 \mathcal L(h,t)=\tfrac12\log\frac{B_t}{-A_t}.
\]
The explicit first-variation formulas imply
$(A_t)_h=A_t/(r-M)^2$, $(B_t)_h=B_t/(r+M)^2$; the terms in $k(h)$ cancel.
Hence $\mathcal L_h=-2rM/q_0^2$, and division by $r_h=-r/q_0$ gives
$L'=2M/q_0$. At the peak the variation factors tend to one and the two
endpoint speeds have the same leading magnitude, so $L\to0$.
Integrating proves the remaining assertions. The initial integrand is
bounded, and transversality makes the terminal angles finite.
\end{proof}

\begin{lemma}\label{cr:reference}
For $K,z\in\mathbb R$, take the profile $zu+Ku^2/2$ on the entire real
axis, its peak at $u=0$, and its outward half-width $s>0$. Write $h_K$ for its
base height, $j_K=z+Kh_K$, and $\Theta_K$ for its midpoint angle.
All these functions are smooth and
\begin{equation}\label{cr:reference-ode}
 \begin{gathered}
 (j_K)_s=-sK\sech^2\Theta_K,\qquad
 (\Theta_K)_s=j_K-\frac{\sinh(2\Theta_K)}s,\\
 |j_K-z|\le |K|s^2/2,\qquad (j_K)_z>0,\qquad(\Theta_K)_z>0.
 \end{gathered}
\end{equation}
For fixed $s>0$, $z\mapsto\Theta_K(s,z)$ is a diffeomorphism onto
$\mathbb R$. Define its inverse $Z_K(s,\theta)$ and
\[
 J_K(s,\theta)=j_K(s,Z_K(s,\theta));\qquad (J_K)_\theta>0.
\]
\end{lemma}
\begin{proof}
The whole-line conclusion of Lemma~\ref{cr:arcs} gives every half-width.
The angular equations prove the first identities and the bound. If $K=0$,
$(j_K)_z=1$. If $K\ne0$, regard the reference as the profile $Ku^2/2$
with peak $t=z/K$, translated to height zero. At fixed half-width its base
height $\bar h(t,s)$ satisfies
$\bar h_t=-r_t/(r_h|_t)>0$. Thus $j_K=K\bar h(z/K,s)$ and
$(j_K)_z=\bar h_t>0$. Differentiating the angle equation gives
\[
 \partial_s(\Theta_K)_z+\frac{2\cosh(2\Theta_K)}s(\Theta_K)_z=(j_K)_z.
\]
The smooth initial expansion $(\Theta_K)_z=s/3+O(s^2)$ and a positive
integrating factor prove strict positivity for all $s>0$.
For surjectivity fix $r>0$. On $(0,r]$, $j_K\ge z-|K|r^2/2$.
For large positive $z$, the angle is positive. While $0\le\Theta_K\le N$
on $[r/2,r]$, its derivative is at least
$z-|K|r^2/2-2\sinh(2N)/r$. If this exceeds $2N/r$, the angle must pass
$N$ and cannot cross it downwards. Thus $\Theta_K(r,z)\to\infty$.
The reflection $\Theta_K(s,z)=-\Theta_{-K}(s,-z)$ proves the other limit.
The inverse function theorem gives the joint smooth inverse and
$(J_K)_\theta=(j_K)_z/(\Theta_K)_z>0$.
\end{proof}

\begin{lemma}\label{cr:comparison}
If $k(u)-Ku$ is strictly increasing, then along any actual fixed-peak
half-orbit the label $z(s)=Z_K(s,\theta(s))$ strictly decreases.
For a transverse terminal half-width $R$ and its terminal label $z_R$,
\begin{equation}\label{cr:angle-compare}
 \theta(s)>\Theta_K(s,z_R)\qquad(0<s<R).
\end{equation}
Strict decrease of $k-Ku$ reverses both conclusions; weak monotonicity
produces weak conclusions.
\end{lemma}
\begin{proof}
Set $E=k(h)-J_K(s,\theta)$, $a=(J_K)_\theta(s,\theta)>0$ and
$\eta=k(h)-Kh$. Since $h$ decreases, $\eta$ is a continuous locally BV
function, and the reference identities imply
\begin{equation}\label{cr:error}
 z'=\frac{E}{(\Theta_K)_z},\qquad dE+aE\,ds=d\eta.
\end{equation}
Indeed differentiation of $J_K$ along the actual angle gives
$dJ_K/ds=Kh'+aE$. The initial expansions, bracketed by reference labels
$k(t)\pm\delta$, imply $z(s)\to k(t)$ and $E(s)\to0$.
Integrating from $e>0$ and letting $e\to0^+$ gives
\begin{equation}\label{cr:error-kernel}
 E(s)=\int_{(0,s]}\exp\!\left(-\int_\sigma^s a(\tau)\,d\tau\right)d\eta(\sigma).
\end{equation}
The limit exists because the initial term is bounded by $|E(e)|$,
the kernel lies in $(0,1]$, and $\eta$ has finite variation on $[0,s]$.
Thus the initial term vanishes even if $a$ is not integrable at zero.
Strict increase of $k-Ku$ makes $d\eta$ negative and
gives negative mass to every interior interval.
The kernel has a positive minimum on each such interval, so $E(s)<0$.
Thus $z'<0$. At the axis, Lemma~\ref{cr:arcs} and the reference inverse
give a finite terminal label. For $s<\sigma<R$,
$z(s)>z(\sigma)\ge z_R$, proving \eqref{cr:angle-compare}.
All reversed and weak cases follow from the same integral.
For a smooth profile the measure equation becomes
\begin{equation}\label{cr:error-smooth}
 E'+aE=-s\bigl(\phi''(h)-K\bigr)\sech^2\theta.
\end{equation}
\end{proof}

\begin{proposition}\label{cr:matching}
Suppose \eqref{cr:basic} holds and the equilibrium is unique in at least
one time direction. Let $J=I_{\phi_+}\cap I_{\phi_-}$ and
\begin{equation}\label{cr:delta}
 \Delta(r)=M_+(0,r)-M_-(0,r),\qquad r\in J.
\end{equation}
Then $J$ is empty or an open interval, and nonconstant periodic orbits
correspond bijectively to zeros of $\Delta$. Limit cycles correspond to
isolated zeros. At a zero of half-width $R$, put $v_\pm=M_{\pm,r}(0,R)$.
The return multiplier is
\begin{equation}\label{cr:multiplier}
 \Pi'=\frac{(1-v_+)(1+v_-)}{(1+v_+)(1-v_-)}
 =\exp\{2(\artanh v_--\artanh v_+)\}>0.
\end{equation}
\end{proposition}
\begin{proof}
Every periodic orbit avoids the equilibrium and has points of both signs
of $x$, since $\int g(x(t))\,dt=0$. Axis intersections have $y\ne0$ and
are transverse. There are finitely many in a period: an accumulation
point would be a nontransverse point of the compact orbit. Right and left
half-return maps are strictly decreasing by $A_t<0<B_t$; their composition
is strictly increasing. An increasing interval map has no finite periodic
sequence of length greater than one, since a strict increase or decrease
at the first iterate persists. Therefore there are exactly two axis
intersections, one above and one below zero, and the two half-orbits have
the same half-width and midpoint.

Conversely, a zero of $\Delta$ joins two half-orbits with common endpoints
$A=M-r<0<B=M+r$. Away from the axis, original time is finite because
$|g|$ is bounded below on compact subarcs. Near the endpoints the inverse
profile equation and \eqref{cr:axis-graph} give unique crossing in finite
time. Thus the joined curve is an actual periodic orbit.
The ordinate $B=r+M_+(0,r)$ has derivative $1+v_+>0$, so half-widths give
local $C^1$ section coordinates. Differentiating the two half-return maps
in these coordinates gives
$dA_+/dB_+=(v_+-1)/(v_++1)$ and
$dB_-/dA_-=(v_-+1)/(v_--1)$, proving \eqref{cr:multiplier}.
Local flow boxes from Lemma~\ref{cr:flow} and compactness give a return
tube around the orbit. Any periodic orbit in that tube meets its section;
the increasing return map has only fixed points as finite periodic
sequences. Hence isolated fixed points, isolated matching zeros, and
isolated periodic orbits are equivalent. Conversely, nearby matching
zeros produce nearby closed curves by continuous dependence and finite
continuous return times. The correspondence holds for nonhyperbolic
cycles as well.
\end{proof}

\begin{lemma}\label{cr:count}
On an open interval, a $C^1$ function whose derivative is positive at every
zero has at most one zero. The same is true for a negative derivative.
If $w'=bq$ with continuous $b>0$ and continuous monotone $q$, then $w$
has at most two isolated zeros.
\end{lemma}
\begin{proof}
Two positive-derivative zeros would provide intermediate points $a<b$
with $w(a)>0>w(b)$. The first zero between them has a nonpositive left
derivative, a contradiction. For the second assertion, when $q$ increases,
$w$ decreases where $q<0$, is constant where $q=0$, and increases where
$q>0$. The two strict intervals each contain at most one zero. A zero
plateau, including its interior endpoints, contains only nonisolated
zeros and precludes any other zeros on the strict intervals. A single
zero at the minimum is the only zero. Replace $w$ by $-w$ for decreasing $q$.
\end{proof}

\begin{proof}[Proof of Theorem~\ref{cr:unique} and Corollary~\ref{cr:unique-derivative}]
The increasing condition is precisely that $k_+-Ku$ increases and
$k_--Ku$ decreases; the negative inverse potential branch reverses order.
Lemma~\ref{cr:flow} and Proposition~\ref{cr:matching} apply.
At a matching half-width $R$, the two terminal midpoint angles, and hence their
reference labels $z_R$, agree. Lemma~\ref{cr:comparison} gives
\[
 \theta_+(s)>\Theta_K(s,z_R)>\theta_-(s),\qquad 0<s<R.
\]
The angular integrals imply $\artanh v_+>\artanh v_-$, so
$\Delta'(R)>0$ and $0<\Pi'<1$. Lemma~\ref{cr:count} on the connected
actual interval $J$ proves uniqueness. Reversing the conditions gives
$\Delta'(R)<0$ and $\Pi'>1$. Finally,
\[
 (q-KG)'=\frac{\mathcal U_K}{g^2}\quad\hbox{almost everywhere}.
\]
Local absolute continuity and the fundamental theorem of calculus give
exactly the stated monotonicity tests, including the interval-strict version.
\end{proof}

\section{Proof of the convex--concave two-cycle criterion}\label{cr:sec-two}
For Theorem~\ref{cr:two}, we fit a quadratic profile to the midpoint and
its derivative with respect to half-width on each half-orbit. Convexity or concavity of the
profile slope then controls the fitted curvature as the half-width varies.
We first prove that this fit exists and is unique, then derive a
transport inequality for its curvature. In the model calculations,
partial derivatives keep the other model parameters fixed.

\subsection{Global quadratic fitting}
For the whole-line profile $\lambda u+\kappa u^2/2$, let
$H(\lambda,\kappa,r)$ denote its relative midpoint at base height zero
and half-width $r>0$. Lemma~\ref{cr:arcs} gives this smooth function for
all $(\lambda,\kappa,r)\in\mathbb R^2\times(0,\infty)$. Put
\begin{equation}\label{cr:model-notation}
 P=H_\lambda,\quad Q=H_\kappa,\quad R=P_r,\quad S=Q_r,
 \qquad \mathcal J=PS-QR.
\end{equation}

\begin{lemma}\label{cr:fit}
The map $(\lambda,\kappa,r)\mapsto(H,H_r,r)$ is a smooth diffeomorphism
onto $\{(M,v,r):r>0,\ |M|<r,\ |v|<1\}$. Everywhere in the model domain,
$P,Q,R,S,\mathcal J>0$.
For fixed $r,M$, let $\lambda_*(r,M,\kappa)$ solve $H=M$, and define
$\widehat v(r,M,\kappa)=H_r(\lambda_*,\kappa,r)$. Then
\begin{equation}\label{cr:conditional}
 \widehat v_M=R/P>0,\qquad \widehat v_\kappa=\mathcal J/P>0.
\end{equation}
\end{lemma}
\begin{proof}
We first establish a variation formula. For a smooth profile based at zero,
with endpoints $a=M-r<0<b=M+r$, vary $\phi$ by $\varepsilon\psi$,
where $\psi(0)=0$, while fixing $r$. The variation $\eta$ satisfies
$\eta_y=\phi'(u)\eta+\psi(u)$ and
$\eta(a)=a\,\delta M$, $\eta(b)=b\,\delta M$. Thus
\begin{equation}\label{cr:variation}
 \delta M=\frac{\int_a^b K_0(y)\psi(u(y))\,dy}
 {bK_0(b)-aK_0(a)},\qquad
 K_0(y)=\exp\left(-\int_a^y\phi'(u(z))\,dz\right).
\end{equation}
The denominator is positive. Taking $\psi=u,u^2/2$ proves $P,Q>0$.

For fixed curvature, Lemma~\ref{cr:reference} uniquely supplies the
reference label having terminal angle $\artanh(M/r)$. Translate its base
to zero to obtain a model with midpoint $M$. Since $P>0$, its slope is
unique and depends smoothly on $(r,M,\kappa)$; this defines $\lambda_*$.
Let
\[
 L_\kappa(s,z)=\int_0^s\frac{\sinh(2\Theta_\kappa(q,z))}{q}\,dq.
\]
Differentiation under the integral is legitimate because
$(\Theta_\kappa)_z=q/3+O(q^2)$ at zero, and gives
\begin{equation}\label{cr:Lz}
 (L_\kappa)_z=\int_0^s
 \frac{2\cosh(2\Theta_\kappa)}q(\Theta_\kappa)_z\,dq>0.
\end{equation}
At fixed $r,M$, write $\widehat L=\artanh\widehat v$.
Comparing curvatures $\kappa+\delta$ and $\kappa$ by
\eqref{cr:error-smooth}, with the same terminal midpoint, shows that
the larger-curvature angle is strictly larger at every interior half-width.
Consequently $\widehat L$ strictly increases with $\kappa$.
The following estimate strengthens this comparison to $\widehat L_\kappa>0$.
On fixed $0<a<b<c<r$, the comparison error satisfies
\[
 E_\delta'+a_\delta E_\delta
 =-\delta s\sech^2\theta_\delta,\qquad E_\delta\le0.
\]
Smooth dependence on the compact interval $[a,r]$ gives uniform bounds
$a_\delta\le A$ and $s\sech^2\theta_\delta\ge q_*>0$ for small
$\delta\ge0$. Keeping only the forcing on $[a,b]$ in the integrating
factor formula gives $E_\delta(s)\le-c_*\delta$ for $b\le s\le r$,
where $c_*=q_*(b-a)e^{-Ar}>0$.
The reference inverse puts all relevant labels in a compact set, on which
$0<m_*\le(\Theta_\kappa)_z\le B_*$. Since their terminal labels agree,
$z_\delta'=E_\delta/(\Theta_\kappa)_z$ yields
\[
 \theta_\delta(s)-\theta_0(s)
 \ge m_*c_*(r-c)\delta/B_*\qquad (b\le s\le c).
\]
The difference is nonnegative elsewhere. Integrating the angle identity
therefore gives
$\widehat L(\kappa+\delta)-\widehat L(\kappa)
\ge2m_*c_*(r-c)B_*^{-1}\log(c/b)\delta$.
Hence $\widehat L_\kappa>0$. Differentiating $H(\lambda_*,\kappa,r)=M$
now proves $\widehat v_\kappa=S-RQ/P=\mathcal J/P>0$.

To determine the range of $\widehat L$, suppose it stayed bounded above
along an increasing sequence $\kappa_n\to+\infty$. The corresponding interior angles
$\theta_n$ would increase pointwise, and, with
$f_n(s)=\sinh(2\theta_n(s))/s$,
\[
 f_n-f_1\ge0,\qquad \int_0^r(f_n-f_1)\,ds\le C.
\]
Monotone convergence implies a finite angle limit almost everywhere.
Choose $0<a<b<c<d<r$ where these limits are finite.
Since $j_n=\theta_n'+f_n$, the integrals of $j_n$ over $[a,b]$ and
$[c,d]$ remain bounded above and below. More explicitly, the endpoint
terms $\theta_n(b)-\theta_n(a)$ and $\theta_n(d)-\theta_n(c)$ converge,
and the integrals of $f_n$ are bounded by the preceding nonnegative
integral estimate. Discard finitely many terms so that $\kappa_n>0$.
Then $j_n'=-s\kappa_n\sech^2\theta_n<0$, and hence
\[
 j_n(b)\le\frac1{b-a}\int_a^b j_n(s)\,ds,\qquad
 j_n(c)\ge\frac1{d-c}\int_c^d j_n(s)\,ds.
\]
These two estimates bound $j_n(b)-j_n(c)$ above.
But
\[
 j_n(b)-j_n(c)=\kappa_n\int_b^c s\sech^2\theta_n(s)\,ds\longrightarrow\infty:
\]
the integral tends by dominated convergence to a strictly positive number.
This contradiction proves $\widehat L\to\infty$. Reflection gives
$\widehat L(r,M,\kappa)=-\widehat L(r,-M,-\kappa)$, hence the negative
limit is $-\infty$. Thus any $v\in(-1,1)$ is fitted by exactly one curvature.

It remains to prove $R,S>0$. Set $D_0=r^2-H^2$, $v=H_r$.
Translation of the base height in a quadratic profile gives
\begin{equation}\label{cr:model-base}
 \kappa P+\lambda=\frac{rv+H}{D_0},\qquad
 \kappa R=\frac{r(1-v^2)}{D_0}
 \left((\artanh v)_r-\frac{2H}{D_0}\right).
\end{equation}
The second formula follows by differentiating the first and substituting
$D_{0,r}=2(r-Hv)$. If $t=t(\lambda,\kappa,r)$ is the model peak height,
then $t_r>0$ and
$\artanh v=L_\kappa(r,\lambda+\kappa t)$.
The angular equation gives
$(\artanh v)_r=2H/D_0+\kappa t_r(L_\kappa)_z$.
For $\kappa\ne0$, cancellation in \eqref{cr:model-base} gives
$R=r(1-v^2)t_r(L_\kappa)_z/D_0>0$; for $\kappa=0$ directly
$R=(1-v^2)(L_0)_z>0$.
Then $S=(\mathcal J+QR)/P>0$.
The map is therefore bijective with Jacobian $\mathcal J>0$,
so its local smooth inverses define a global smooth inverse.
Finally differentiating with respect to $M$ proves the other identity
in \eqref{cr:conditional}.
\end{proof}

For any smooth actual profile, use this inverse to fit its data:
\begin{equation}\label{cr:actual-fit}
 M(h,r)=H(\lambda(h,r),\kappa(h,r),r),\quad
 M_r(h,r)=H_r(\lambda(h,r),\kappa(h,r),r).
\end{equation}
Write $k=\phi'$, $\alpha=r/(r^2-M^2)$ and
$D=\partial_h-\alpha\partial_r$; $D$ is differentiation in base height
along a fixed-peak half-orbit. The angular equations give
$DM=M/(r^2-M^2)-k$ and
$Dv=-2Mr(1-v^2)/(r^2-M^2)^2$.
Subtracting the corresponding quadratic-model equations from the
chain-rule derivatives of \eqref{cr:actual-fit} proves
\begin{equation}\label{cr:fit-transport}
 P\lambda_r+Q\kappa_r=0,\qquad
 D\kappa=-\frac R{\mathcal J}(\lambda-k),\qquad
 D\lambda-\kappa=\frac S{\mathcal J}(\lambda-k).
\end{equation}
Indeed the last two identities are the inverse of the matrix equation
$\left(\begin{smallmatrix}P&Q\\R&S\end{smallmatrix}\right)
(D\lambda-\kappa,D\kappa)^T=(\lambda-k,0)^T$.

\subsection{The quadratic-model inequality}
\begin{lemma}\label{cr:model-inequality}
For all real $\lambda,\kappa$ and all $r>0$,
\begin{equation}\label{cr:model-target}
 \left.\partial_r\left(\frac R{\mathcal J}\right)\right|_{\lambda,\kappa}\le0.
\end{equation}
\end{lemma}
We prove the inequality by a transport equation for a determinant of the
model derivatives. Set
\[
 W=RS_r-SR_r,\quad \alpha=\frac r{r^2-H^2},\quad
 \gamma=\frac{r^2+H^2+2rHH_r}{(r^2-H^2)^2},\quad
 D=\kappa\partial_\lambda-\alpha\partial_r.
\]
Since $\mathcal J_r=PS_r-QR_r$,
\begin{equation}\label{cr:model-ratio}
 (R/\mathcal J)_r=-PW/\mathcal J^2.
\end{equation}
From $DH=H/(r^2-H^2)-\lambda$, differentiation of both $H$ and the
coefficients of $D$ gives $DP=\gamma P-1$, $DQ=\gamma Q-P$.
Using $D(f_r)=(Df)_r+\alpha_rf_r$ twice gives
\begin{align}
 DR&=(\gamma+\alpha_r)R+\gamma_rP,&
 DS&=(\gamma+\alpha_r)S+\gamma_rQ-R,\notag\\
 D(R_r)&=(\gamma+2\alpha_r)R_r+(2\gamma_r+\alpha_{rr})R+\gamma_{rr}P,\notag\\
 D(S_r)&=(\gamma+2\alpha_r)S_r+(2\gamma_r+\alpha_{rr})S+\gamma_{rr}Q-R_r.
 \label{cr:model-transport}
\end{align}
Consequently
\begin{equation}\label{cr:DW}
 DW=(2\gamma+3\alpha_r)W+\gamma_r\mathcal J_r-\gamma_{rr}\mathcal J.
\end{equation}
For later use the smooth scaling $u=r^2U$, $y=rY$ gives
$H=r\mathscr H(\lambda r,\kappa r^3)$, with $\mathscr H$ smooth near zero.
The variation formula \eqref{cr:variation} at the parabola
$U=(1-Y^2)/2$ gives $\mathscr H_\lambda(0)=1/3$ and
$\mathscr H_\kappa(0)=1/15$. Thus, with differentiable remainders uniform
on compact parameter sets,
\begin{equation}\label{cr:model-small}
 \begin{array}{lll}
 P=r^2/3+O(r^3),&R=2r/3+O(r^2),&R_r=2/3+O(r),\\
 Q=r^4/15+O(r^5),&S=4r^3/15+O(r^4),&S_r=4r^2/5+O(r^3),
 \end{array}
\end{equation}
and $W=16r^3/45+O(r^4)>0$ initially.
At a zero of $W$, set $\sigma_0=R_r/R$. Then
$S_r=\sigma_0S$, $\mathcal J_r=\sigma_0\mathcal J$ and
\begin{equation}\label{cr:DW-zero}
 DW=-\mathcal J(\gamma_{rr}-\sigma_0\gamma_r).
\end{equation}
It remains to show that the parenthesis is positive for $\kappa>0$
at every zero of $W$.

Write $\kappa=2b>0$, $\lambda=d$, and write the endpoints as $-m,n$,
where
\[
 m=r-H>0,\quad n=r+H>0,\quad
 \ell=\frac{dn}{dm}=\frac{1+H_r}{1-H_r}>0,\quad
 w=\frac{\ell_m}{\ell},\quad
 \mathscr S=\frac{\ell_{mm}}\ell-\frac32w^2.
\]
Since $m_r=1-H_r>0$, $m$ is a legitimate coordinate.
Define
\begin{align}
 X&=\frac{\ell/m-1/n}{1+\ell}=\kappa P+d,&
 Z&=\frac{\ell m-n}{1+\ell}=dP+3\kappa Q,\notag\\
 \gamma&=\frac{\ell/m^2+1/n^2}{1+\ell},\label{cr:XZ}\\
 D_*&=(m-n)\mathscr S+
 \left[1-\frac{2n}{m}+\ell\left(\frac{2m}{n}-1\right)\right]w
 +\frac{m+n}{2}w^2,\notag\\
 w_*&=\frac{2(\ell m^4-2\ell m^3n+2mn^3-n^4)}{m^2n^2(m+n)}.
 \label{cr:Dstar}
\end{align}
The formula for $X$ is \eqref{cr:model-base}; that for $Z$ follows from
the scaling identity $rH_r-H=\lambda P+3\kappa Q$.

\begin{lemma}\label{cr:endpoint-algebra}
The following identities hold:
\begin{equation}\label{cr:endpoint-det}
 X_mZ_{mm}-Z_mX_{mm}
 =\frac{\ell^2(m+n)^2}{m^2n^2(1+\ell)^3}D_*
 =\frac{3\kappa^2}{m_r^3}W.
\end{equation}
At $W=0$, $m\ne n$, putting $\sigma=X_{mm}/X_m$ is legitimate and gives
\begin{equation}\label{cr:endpoint-gamma}
 \gamma_{rr}-\sigma_0\gamma_r
 =m_r^2(\gamma_{mm}-\sigma\gamma_m),\qquad
 \gamma_{mm}-\sigma\gamma_m
 =\frac{\ell(m+n)^2}{(1+\ell)m^2n^2}\frac{w_*-w}{m-n}.
\end{equation}
\end{lemma}
\begin{proof}
Let $B=n^{-2}-m^{-2}$. Differentiating $n_m=\ell$, $\ell_m=\ell w$ gives
\[
 Z_m=\frac{w\ell(m+n)}{(1+\ell)^2},\quad
 X_m=\frac{\ell B}{1+\ell}+\frac{Z_m}{mn},\quad
 B_m=2m^{-3}-2\ell n^{-3},\quad w_m=\mathscr S+w^2/2.
\]
Substitution into the determinant, and division by
$\ell^2(m+n)/(1+\ell)^3$, leaves
\[
 B\mathscr S+
 \left(\frac{B(1+\ell)}{m+n}-B_m\right)w+
 \frac{(m+n)^2}{2m^2n^2}w^2.
\]
Here $B=(m-n)(m+n)/(m^2n^2)$ and the coefficient of $w$ equals
$(m+n)[1-2n/m+\ell(2m/n-1)]/(m^2n^2)$, proving the first identity.
The second follows by differentiating \eqref{cr:XZ} with respect to $r$.
Also $X_m=\kappa R/m_r>0$. At $W=0$, $Z_{mm}=\sigma Z_m$ and
$\sigma m_r^2=\sigma_0m_r-m_{rr}$, proving the first identity in
\eqref{cr:endpoint-gamma} by the chain rule.

For the remaining algebra, for a function $f$ near $m,-n$, set
\[
 T_f=\frac{\ell f(m)+f(-n)}{1+\ell},\quad
 a_0=\frac\ell{1+\ell},\quad \beta=\frac w{1+\ell},\quad
 \rho=(\log a_0)_m-\sigma.
\]
Differentiation twice gives
\begin{align*}
 \frac{(T_f)_{mm}-\sigma(T_f)_m}{a_0}
 ={}&f''(m)+\ell f''(-n)+\beta\{f'(m)+\ell f'(-n)\}\\
 &+\rho\{f'(m)-f'(-n)\}
 +(\beta_m+\rho\beta)\{f(m)-f(-n)\}.
\end{align*}
The choices $f(x)=x,1/x$ have left-hand side zero, since $T_x=Z$ and
$T_{1/x}=X$. They give
\[
 \beta_m+\rho\beta=-\frac{\beta(1+\ell)}{m+n},\quad
 \rho B=\frac{2\ell}{n^3}-\frac2{m^3}
 +\beta\left(\frac1{m^2}+\frac\ell{n^2}+\frac{1+\ell}{mn}\right).
\]
For $f(x)=x^{-2}$, $T_f=\gamma$ and the expression is
\[
 \frac6{m^4}+\frac{6\ell}{n^4}
 +\beta\left(-\frac2{m^3}+\frac{2\ell}{n^3}\right)
 -2\rho\left(\frac1{m^3}+\frac1{n^3}\right)
 +\frac{\beta(1+\ell)B}{m+n}.
\]
Eliminating $\rho$ leaves constant term
$2(m+n)N/[m^4n^4(m-n)]$ and coefficient of $\beta$
$-(1+\ell)(m+n)^2/[m^2n^2(m-n)]$, where
$N=\ell m^4-2\ell m^3n+2mn^3-n^4$.
The two identities after clearing denominators are
\begin{align*}
 6(n^4+\ell m^4)(m^2-n^2)-4(m^3+n^3)(\ell m^3-n^3)
 &=2(m+n)^2N,\\
 (m+n)\bigl[2(\ell m^3-n^3)(m-n)-2(m^3+n^3)(\ell m+n)&\\
 {}+(1+\ell)mn(m-n)^2\bigr]
 &=-(1+\ell)mn(m+n)^3.
\end{align*}
Substitution of $\beta=w/(1+\ell)$ proves the formula. The only additional
division is by $B\ne0$, so the calculation also covers $w=0$ and
$\beta=0$.
\end{proof}

\begin{lemma}\label{cr:riccati}
For the quadratic profile $du+bu^2$, $b>0$, its endpoint Schwarzian satisfies
\begin{equation}\label{cr:schwarzian}
 \mathscr S=2b(m+n\ell^2)+\frac{d^2}{2}(\ell^2-1)
 -d\left(\frac1m+\frac{\ell^2}n\right)
 +\frac32\left(\frac{\ell^2}{n^2}-\frac1{m^2}\right).
\end{equation}
There is a positive function $V$ on $[-m,n]$ such that
\begin{equation}\label{cr:w-integral}
 \ell=\frac{mV(-m)^2}{nV(n)^2},\qquad
 w=\frac1m-\frac\ell n+d(1+\ell)
 +2bmV(-m)^2\int_{-m}^nV(y)^{-2}\,dy.
\end{equation}
If $d<0$, $d^2\ge4bm$, and $\nu=\sqrt{d^2/4-bm}$, then
\begin{equation}\label{cr:w-bound}
 w\le\frac1m-\frac\ell n+d\ell-2\nu.
\end{equation}
\end{lemma}
\begin{proof}
Along the actual arc define
$V(y)=\exp[-\int_{-m}^y(bu(s)+d/2)\,ds]>0$ and $U=uV$.
Direct differentiation gives the regular linear system
\begin{equation}\label{cr:linear-riccati}
 U'=dU/2-yV,\qquad V'=-bU-dV/2.
\end{equation}
Near the nonzero endpoints, elimination gives
$U''=U'/y+(by+d^2/4-d/(2y))U$.
Complete $(U_1,V_1)=(U,V)$ to a fundamental pair; its determinant
$\Omega=U_1V_2-U_2V_1$ is a nonzero constant. At both endpoints,
$U_1=0$, $V_1>0$, so $U_2=-\Omega/V_1\ne0$.
Thus the same ratio $f_0=U_1/U_2$ is defined near both endpoints, and
$f_0'=y\Omega/U_2^2\ne0$ there. Nearby arcs correspond to linear
combinations $U_1-qU_2$, whence
\begin{equation}\label{cr:common-label}
 f_0(n(m))=f_0(-m).
\end{equation}
The regular fundamental system is defined across $y=0$. The ratio is
needed only near the endpoints and may have poles in the interior.
For $\{f,y\}=f'''/f'-\tfrac32(f''/f')^2$, the quotient of two solutions
of $a''=Q_0a$ has Schwarzian $-2Q_0$, as follows from
$f''/f'=-2a_2'/a_2$. The chain rule is
$\{f\circ h,m\}=\{f,h(m)\}(h')^2+\{h,m\}$.
Substitute $U=\sqrt{|y|}\widetilde U$ in the scalar endpoint equation.
It gives
$\widetilde U''=(by+d^2/4-d/(2y)+3/(4y^2))\widetilde U$;
hence $\{f_0,y\}=-2by-d^2/2+d/y-3/(2y^2)$.
Apply the chain rule to \eqref{cr:common-label} to get
$\mathscr S=\{f_0,y\}|_{-m}-\ell^2\{f_0,y\}|_n$,
which is \eqref{cr:schwarzian}.

Differentiation of $u(-m;m)=0$ gives $u_m(-m)=m$.
The variational equation yields $u_m(y)=mV(-m)^2/V(y)^2$.
Differentiating $u(n(m);m)=0$ and using $u_y(n)=-n$ proves the formula
for $\ell$. Differentiate
$\log\ell=\log m-\log n+\int_{-m}^n(d+2bu)\,dy$ to obtain
\eqref{cr:w-integral}.
Finally $V''=(by+d^2/4)V$, $V'(-m)=-dV(-m)/2$.
Writing $s=y+m$ and $v_0(s)=V(-m+s)/V(-m)$ gives
$v_0''=(\nu^2+bs)v_0$, $v_0(0)=1$, $v_0'(0)=-d/2$.
The positive comparison function
\[
 f_\nu(s)=
 \begin{cases}\cosh(\nu s)-\dfrac d{2\nu}\sinh(\nu s),&\nu>0,\\
 1-ds/2,&\nu=0
 \end{cases}
\]
satisfies $v_0-f_\nu=b\int_0^sK_\nu(s-t)t\,v_0(t)\,dt\ge0$,
where $K_\nu(t)=\sinh(\nu t)/\nu$ or $t$ when $\nu=0$.
Therefore
\[
 V(-m)^2\int_{-m}^nV^{-2}\,dy
 \le\int_0^\infty f_\nu(s)^{-2}\,ds=\frac1{\nu-d/2}.
\]
For $\nu>0$ the last integral follows by $q=\tanh(\nu s)$; for $\nu=0$
it is $-2/d$. Since $2bm/(\nu-d/2)=-d-2\nu$, substitution proves
\eqref{cr:w-bound}.
\end{proof}

\begin{lemma}\label{cr:zero-sign}
If $\kappa>0$ and $W=0$, then $\gamma_{rr}-\sigma_0\gamma_r>0$.
\end{lemma}
\begin{proof}
By \eqref{cr:endpoint-det}, $D_*=0$ at such a point.
Along its fixed-peak reference arc, put $\theta=\artanh(H/r)$ and
$L=\artanh H_r$. Then $j'=-s\kappa\sech^2\theta<0$,
$\theta'=j-\sinh(2\theta)/s$, $L'=\sinh(2\theta)/s$.
The integrating factor for
$\theta'+a(s)\theta=j(s)$, with
$a=\sinh(2\theta)/(s\theta)>0$ (value $2/s$ at zero angle), is valid
from $e>0$; the initial contribution tends to zero as $e\to0^+$.
We shall also use $(L_\kappa)_z>0$ from \eqref{cr:Lz}.

First suppose $m\le n$. The terminal angle is nonnegative, so the peak
slope $z$ must be positive; otherwise $j(s)<0$ and the preceding positive
integrating kernel would make the angle negative. The angle starts
positive. If it first reached zero before the terminal half-width, then $j\le0$
there and subsequently $j<0$, preventing a return to a nonnegative angle.
Thus $\theta>0$ before the endpoint and $L(r)>0$, so $\ell>1$.
At fixed $(d,\kappa)$, $t_r>0$ and $z=d+\kappa t$ give
\[
 L_r=\frac{\sinh(2\theta)}r+\kappa t_r(L_\kappa)_z>0,
 \qquad w=2L_r/m_r>0.
\]
If $m=n$, $D_*=(\ell-1)w+mw^2>0$, a contradiction.
If $m<n$, put $\zeta=m/n\in(0,1)$. The numerator of $w_*$ has the sign of
\[
 \ell\zeta^3(\zeta-2)+2\zeta-1
 \le\zeta^3(\zeta-2)+2\zeta-1
 =(\zeta-1)^3(\zeta+1)<0.
\]
Thus $(w_*-w)/(m-n)>0$, giving the claim by
\eqref{cr:endpoint-gamma}.

Next suppose $m>n$, $\ell\ge1$. The negative terminal angle forces $d<0$:
if $d\ge0$, the decreasing $j$ is positive earlier, forcing a positive
angle. Every group in \eqref{cr:schwarzian} is now nonnegative and one
is strictly positive; hence $\mathscr S>0$.
We justify $w\ge0$ in this region. The coordinate $m$ ranges over an
interval $(0,M_*)$ and initially
$\ell=1+4dr/3+O(r^2)<1$. On any component $(a,c)$ of
$\{m-n>0,\ell>1\}$, $(m-n)_m=1-\ell<0$; its left endpoint is an
interior positive point with $m-n>0$ and $\ell=1$.
There $\ell_m\ge0$. Since $w_m=\mathscr S+w^2/2>0$, $w>0$ within
the component. At a point $m>n$, $\ell=1$, a value $w<0$ would put
$\ell>1,w<0$ just to its left, a contradiction. A value $w=0$ has
$w_m=\ell_{mm}=\mathscr S>0$ and yields the same contradiction just
to the left. Thus $w\ge0$ wherever required. With $\zeta=m/n>1$,
$1-2/\zeta+\ell(2\zeta-1)\ge2\zeta-2/\zeta>0$.
It follows that $D_*>0$, so there is no zero in this case.

It remains to consider $m>n$, $0<\ell<1$, again with $d<0$.
We prove $w\ge w_*\Longrightarrow D_*>0$. Scale at this point by the
constant $N=n$: $y=N\widetilde y$, $u=N^2\widetilde u$.
Then $\widetilde d=Nd$, $\widetilde b=N^3b$,
$\widetilde w=Nw$, $\widetilde{\mathscr S}=N^2\mathscr S$,
$\widetilde D_*=ND_*$ and $\widetilde w_*=Nw_*$.
The signs are unchanged, so take $n=1$, $m=z>1$. Now
\begin{equation}\label{cr:wstar-scaled}
 w_*=\frac{2[\ell z^3(z-2)+2z-1]}{z^2(z+1)}>0.
\end{equation}
For $z\ge2$ positivity is immediate from the displayed terms. For
$1<z<2$, replacing $\ell<1$ by $1$ lowers the numerator to
$(z-1)^3(z+1)>0$.
As a quadratic in $w$, $D_*$ has at $w_*$ a derivative whose product
with $z^2$ is
$J_0=2\ell z^4-2\ell z^3-\ell z^2+z^2+2z-2$.
It is affine in $\ell$; with $\xi=z-1>0$ its endpoint values are
$\xi^2+4\xi+1$ and $2\xi(\xi+2)(\xi^2+\xi+1)$.
Thus $D_*$ strictly increases on $w\ge w_*$.
Put
\[
 T=2\ell z^3-\ell z^2-z+2,\quad
 A=2\ell z^3-3\ell z^2+3z-2,\quad
 \mathscr S_0=\tfrac32(\ell^2-z^{-2}).
\]
Expansion of \eqref{cr:wstar-scaled} gives
\begin{align}
 (z-1)\mathscr S_0+[1-2/z+\ell(2z-1)]w_*+\tfrac{z+1}{2}w_*^2
 &=\frac{z-1}{2z^4}AT,\label{cr:AT}\\
 \Delta_0:=1/z-\ell-w_*&=-\frac{z-1}{(z+1)z^2}T.\label{cr:Delta0}
\end{align}
Multiplying the left side of \eqref{cr:AT} by $2z^4(z+1)$ gives
\[
 (z^2-1)\bigl[\ell^2z^4(2z-3)(2z-1)
 +4\ell z^2(z^2-1)+(3z-2)(2-z)\bigr],
\]
and the bracket equals $AT$. Moreover $A>0$: its endpoint values in
$\ell$ are $3\xi+1$ and $\xi(2\xi^2+3\xi+3)$.

Assume $w\ge w_*$. If $T\ge0$ and $d^2\ge4bz$,
\eqref{cr:w-bound} implies
$\Delta_0\ge1/z-\ell-w\ge-d\ell+2\nu>0$,
contrary to \eqref{cr:Delta0}. Thus $d^2<4bz$, and
\[
 \mathscr S-\mathscr S_0
 =(2bz-d^2/2)+\ell^2(2b+d^2/2-d)-d/z>0.
\]
Equation \eqref{cr:AT} is nonnegative, so the strict increase in $w$
and $z-1>0$ imply $D_*>0$.
If $T<0$, then $z-2>\ell z^2(2z-1)>0$, so $z>2$ and $z-2>3\ell$.
By $X=d+\kappa P>d$,
\[
 d<\frac{\ell/z-1}{1+\ell},\qquad -dz>\frac{z-\ell}{1+\ell}>2.
\]
Rewrite the Schwarzian as
\[
 \mathscr S=(2bz-d^2/2)+\ell^2(2b+d^2/2-d+3/2)
 -d/z-3/(2z^2).
\]
If the first parenthesis is nonnegative this is positive.
Otherwise $\nu=\sqrt{d^2/4-bz}>0$, and \eqref{cr:w-bound} yields
$2\nu<\Delta_0<1/z$. Hence its first parenthesis is
$-2\nu^2>-1/(2z^2)$; with $-d/z>2/z^2$ this again gives
$\mathscr S>0$. Since $z>2$ makes the coefficient of $w$ positive,
$D_*>0$. We have proved the implication in all cases.
Therefore a zero in this final region must have $w<w_*$.
Now $m-n>0$ and \eqref{cr:endpoint-gamma} prove the required strict sign.
\end{proof}

\begin{proof}[Proof of Lemma~\ref{cr:model-inequality}]
First let $\kappa>0$. Every finite model point lies on a fixed-peak
characteristic whose slope tends to a finite peak slope as $s\to0^+$.
The uniform expansion \eqref{cr:model-small} gives $W>0$ initially.
Along this characteristic the outward derivative is $-D/\alpha$.
At a first zero of $W$, \eqref{cr:DW-zero} and Lemma~\ref{cr:zero-sign}
would give
\[
 \frac{dW}{ds}=\frac{\mathcal J}{\alpha}
 (\gamma_{rr}-\sigma_0\gamma_r)>0,
\]
whereas a first zero reached from positive values has derivative at most
zero. Thus $W>0$ throughout. Reflection
$H(-\lambda,-\kappa,r)=-H(\lambda,\kappa,r)$ preserves $P,Q,R,S$
and $W$, proving the case $\kappa<0$. Smooth dependence as
$\kappa\to0$ gives $W\ge0$ at zero curvature. Equation
\eqref{cr:model-ratio} converts this sign of $W$ into the
fixed-parameter inequality \eqref{cr:model-target}.
\end{proof}

\subsection{Shape comparison and the axis limit}
\begin{lemma}\label{cr:shape}
For $\phi\in C^\infty((0,U))$, the fitted parameters in
\eqref{cr:actual-fit} satisfy
\begin{equation}\label{cr:fit-small}
 \begin{split}
 \lambda(h,r)&=\phi'(h)-\phi'''(h)r^4/35+r^5A(h,r),\\
 \kappa(h,r)&=\phi''(h)+2\phi'''(h)r^2/7+r^3B(h,r),
 \end{split}
\end{equation}
where the remainders are smooth near each positive base height and $r=0$.
If $\phi'''>0$, then $\lambda<\phi'(h)$ and $\kappa_r>0$ throughout
$\Omega_\phi$; if $\phi'''<0$, both inequalities reverse.
Nonstrict third-derivative signs give the corresponding weak inequalities.
\end{lemma}
\begin{proof}
Normalize by $u-h=r^2\xi$, $y-\phi(h)=r\zeta$. To control the remainders,
allow a signed small parameter $\rho$ and use the smooth profile
\[
 \psi(h,\rho,\xi)=\rho\xi\int_0^1\phi'(h+\tau\rho^2\xi)\,d\tau
 =\rho k(h)\xi+\rho^3\phi''(h)\xi^2/2
  +\rho^5\phi'''(h)\xi^3/6+\rho^7\mathcal R.
\]
All heights are restricted to a fixed compact subset of $(0,U)$.
At the zero normalized profile, the arc is $(1-\zeta^2)/2$, and
\eqref{cr:variation} gives the responses
\[
 \begin{array}{c|cc}
 \text{profile variation}&\delta M&\delta M_r\\\hline
 \xi&1/3&2/3\\
 \xi^2/2&1/15&4/15\\
 \xi^3/6&1/105&6/105
 \end{array}
 \qquad
 \begin{pmatrix}1/3&1/15\\2/3&4/15\end{pmatrix}^{-1}
 \binom{1/105}{6/105}=\binom{-1/35}{2/7}.
\]
The matrix determinant is $2/45$. The implicit function theorem applied
to the finite-dimensional family consisting of the quadratic part plus
an independent multiple of the cubic remainder gives smooth fitted
normalized coefficients. Taylor's formula in that independent parameter
shows that their differences from $(\rho k,\rho^3\phi'')$ are
$\rho^5\phi'''(-1/35,2/7)+\rho^6$ times a smooth vector.
Since these coefficients equal $(r\lambda,r^3\kappa)$ at $\rho=r>0$,
\eqref{cr:fit-small} follows. In particular,
$\kappa_r=4\phi'''(h)r/7+O(r^2)$ with controlled derivative remainder.

Assume $\phi'''>0$. Fix a target arc $(h,r)$ and follow it from its peak
$t$ outwards. Put $c(s)=\phi''(h(s))$, which strictly decreases.
Fix its \emph{terminal} fitted curvature $\kappa_*$. The comparison with
this one constant curvature gives
\[
 E'+aE=-s\sech^2\theta\,[c(s)-\kappa_*],\quad a>0,\quad E(0+)=0,
 \qquad z'=E/(\Theta_{\kappa_*})_z.
\]
Its positive-kernel formula is \eqref{cr:error-kernel}, now with an
ordinary forcing $b(s)=-s\sech^2\theta[c(s)-\kappa_*]$.
At the terminal half-width, the equal-midpoint fit identifies
$J_{\kappa_*}(r,\theta(r))=\lambda(h,r)$.
The equal-$M_r$ fit and the angular integrals imply, with $z_*=z(r)$,
\begin{equation}\label{cr:integral-match}
 \int_0^r\frac{\sinh(2\theta(s))}{s}\,ds
 =\int_0^r\frac{\sinh(2\Theta_{\kappa_*}(s,z_*))}{s}\,ds.
\end{equation}
If $\kappa_*\ge c(0)$, then $b>0$ and $E>0$, so $z$ increases and
$\theta(s)<\Theta_{\kappa_*}(s,z_*)$ for $s<r$, contradicting
\eqref{cr:integral-match}. If $\kappa_*\le c(r)$, the reverse strict
comparison gives the same contradiction. Otherwise $b$ is first negative,
then positive, with one change at $s_0$. The integrating-factor multiple
of $E$ is strictly increasing after $s_0$, and $E<0$ before $s_0$.
If $E(r)\le0$, the positive integrating-factor multiple cannot have
crossed zero after $s_0$, since it is strictly increasing there.
Thus $E<0$ on $(0,r)$, again violating
\eqref{cr:integral-match}. Therefore
$k(h)-\lambda(h,r)=E(r)>0$.

For curvature monotonicity put $C(\lambda,\kappa,r)=R/\mathcal J$.
Differentiate the second equation of \eqref{cr:fit-transport} with respect
to $r$ at fixed $h$, and use $\lambda_r=-(Q/P)\kappa_r$:
\begin{equation}\label{cr:shape-transport}
 D(\kappa_r)=T\kappa_r-C_r(\lambda-k),\qquad
 T=\alpha_r+C\frac QP-\left(C_\kappa-\frac QP C_\lambda\right)(\lambda-k).
\end{equation}
Here $C_r$ denotes differentiation with the model parameters fixed.
Along the fixed-peak characteristic,
$w(s)=\kappa_r(h(s),s)$ satisfies
\[
 w'=-\frac T\alpha w+\frac{C_r}{\alpha}(\lambda-k).
\]
The last term is nonnegative by Lemma~\ref{cr:model-inequality} and the
strict slope inequality just proved. Formula \eqref{cr:fit-small} gives
$w(s)>0$ initially because $h(s)\to t>0$ and $\phi'''(t)>0$.
An integrating factor on a compact interval from such an initial
positive half-width to $r$ preserves strict positivity for either sign of $T$;
the coefficient $T/\alpha$ need only be integrable on that interval.
Reflection of the profile proves the strictly negative case.

For $\phi'''\ge0$, apply the strict result to
$\phi_\varepsilon=\phi+\varepsilon u^3$. On any compact subset of
$\Omega_\phi$, all peak heights stay in a compact positive-height
interval below $U$. Smooth parameter dependence, endpoint transversality,
and $r_t>0$ put these arcs and their fixed-half-width inverses in a common
domain for small $\varepsilon$. Their data converge with every fixed
finite derivative order. The global smooth model inverse passes this
convergence to $(\lambda,\kappa)$, including $\kappa_r$.
Letting $\varepsilon\to0^+$ gives the weak inequalities.
Reflection gives the nonpositive case. This convergence uses only local
profile data.
\end{proof}

\begin{lemma}\label{cr:axis-shape}
Let $\phi\in C^1((0,U))\cap C([0,U))$, $\phi(0)=0$.
If $k=\phi'$ is convex, its fitted curvature $\kappa(h,r)$ is
nondecreasing in $r$ at each positive base height, and
$\lambda(h,r)\le k(h)$. If $k$ is concave the inequalities reverse.
The axis fit $(\lambda(0,r),\kappa(0,r))$ is finite and continuous on
$I_\phi$, and is the locally uniform limit of positive-height fits.
Its curvature is likewise nondecreasing for convex $k$ and nonincreasing
for concave $k$.
\end{lemma}
\begin{proof}
Fix a positive base $h$ and two actual half-widths $r_1<r_2$, with peaks
$T_1<T_2<U$. Choose a compact interval $(a,b)$ containing $h$ and these
peaks, with a larger compact neighborhood still inside $(0,U)$.
Convolve $\phi$ on this larger neighborhood by a nonnegative smooth
unit-mass kernel of shrinking support and add $\varepsilon_nu^3$,
where $(\varepsilon_n)$ is a positive decreasing sequence tending to zero. The approximants converge in $C^1$
on $[a,b]$, and their third derivatives are strictly positive: convolution
preserves convexity of $k$, and the added third derivative is $6\varepsilon_n$.
For a $C^1$ scalar differential equation, the solutions and first
variations converge uniformly under $C^1$ convergence of the right-hand
side on a compact tube. Transversality and $r_t>0$ thus preserve these
actual half-widths and give convergence of $(M,M_r)$ at them.
Translate the base of the compact interval to zero to apply
Lemma~\ref{cr:shape} wholly inside that interval; this translation does
not change the relative midpoint or its fit. It gives
$\kappa_n(h,r_1)\le\kappa_n(h,r_2)$ and $\lambda_n(h,r_i)\le\phi_n'(h)$.
Passing to the limit through the continuous inverse of Lemma~\ref{cr:fit}
proves the same inequalities for the continuous fitted curvature.
Apply the argument to $-\phi$ for concavity.

On a compact interval of half-widths at the axis the data $(M,M_r,r)$ lie in a compact
subset of $|M|<r$, $|M_r|<1$. Lemma~\ref{cr:arcs} provides common
small positive heights and uniform convergence of these data.
Uniform continuity of the model inverse gives finite continuous axis
parameters and uniform convergence of the fits. Apply the positive-height
ordering to two fixed half-widths and let $h\to0^+$ to obtain the weak
ordering at the axis. The smoothing argument takes place entirely at
positive heights, before this limit, so continuity of
$\phi$ at zero is sufficient.
\end{proof}

\begin{proof}[Proof of Theorem~\ref{cr:two} and Corollary~\ref{cr:two-derivative}]
Under the displayed direction in \eqref{cr:secants}, the positive inverse
potential branch preserves order, so $k_+$ is convex; the negative branch
reverses order, so $k_-$ is concave. Lemma~\ref{cr:flow} and
Proposition~\ref{cr:matching} apply. On the common actual interval $J$,
Lemma~\ref{cr:axis-shape} gives the continuous nondecreasing function
$Q_*(r)=\kappa_+(0,r)-\kappa_-(0,r)$.
Using the exact fits and \eqref{cr:conditional}, the fundamental theorem
of calculus in midpoint and curvature separately gives
\begin{equation}\label{cr:delta-linear}
 \Delta'=A(r)\Delta+B(r)Q_*(r),\qquad B(r)>0,
\end{equation}
where, writing $M_\pm=M_\pm(0,r)$ and $\kappa_\pm=\kappa_\pm(0,r)$,
\begin{align*}
 A(r)&=\int_0^1\widehat v_M(r,M_-+t(M_+-M_-),\kappa_+)\,dt,\\
 B(r)&=\int_0^1\widehat v_\kappa(r,M_-,\kappa_-+t(\kappa_+-\kappa_-))\,dt.
\end{align*}
These coefficients are continuous on all of $J$: the midpoint segment
lies in $(-r,r)$, and the curvature has no restricted range.
Multiplication by $\mu(r)=\exp[-\int_{r_*}^r A(s)\,ds]>0$ gives
$(\mu\Delta)'=\mu BQ_*$. Lemma~\ref{cr:count} gives at most two isolated
zeros, including all zero-derivative cases and excluding zero plateaus.
The matching proposition completes the proof. Reversing convexity gives
a nonincreasing $Q_*$ and the same count.

For the derivative version, define on either side
$H_0=(f'g-fg')/g^3$. Under the stated local absolute continuity assumptions,
$H_0$ is locally absolutely continuous and
\begin{equation}\label{cr:derivative-identities}
 H_0'=\frac{\mathcal N}{g^4},\qquad
 k_\pm'(u)=H_0(x_\pm(u)),\qquad
 \phi_\pm'''(u)=\left.\frac{\mathcal N}{g^5}\right|_{x=x_\pm(u)}
\end{equation}
almost everywhere where applicable. Thus $\mathcal N\ge0$ makes $H_0$
nondecreasing in $x$, $k_+'$ nondecreasing in $u$, and $k_-'$
nonincreasing in $u$. The opposite sign reverses everything.
These are precisely the required convexity and concavity conditions.
\end{proof}

\begin{remark}
Absolute continuity cannot be omitted from an almost-everywhere derivative
test: integrating $\phi''=-C$, where $C$ is a nonconstant Cantor staircase,
gives a $C^2$ profile with $\phi'''=0$ almost everywhere but with
$\phi''$ decreasing. An almost-everywhere sign alone does not control
its singular variation. The convex--concave statement itself includes
continuous singular variation and has no such omission.
\end{remark}

\section{Proof of the common-intercept criterion}\label{cr:sec-intercept}
The common-intercept criterion follows from the global fit, its
transport identities, and its limit as the half-width tends to zero, without using
Lemma~\ref{cr:model-inequality} or the preceding shape comparison.
All comparisons take place on finite arcs within the potential ranges
of the original profiles.

\begin{lemma}\label{cr:intercept-positive}
Let $\phi\in C^\infty((0,U))$, $k=\phi'$. If
$T(u)=(k(u)-d)/u$ has $T'(u)>0$, then
$\lambda(h,r)-h\kappa(h,r)<d$ throughout the actual domain $\Omega_\phi$.
If $T'<0$, the fitted-intercept inequality reverses.
\end{lemma}
\begin{proof}
Follow an actual half-orbit of peak $t<U$ outward and set
$\mathscr I=\lambda-h\kappa-d$. The limit as $r\to0^+$ in
\eqref{cr:fit-small} gives
\[
 \mathscr I(0+)=k(t)-tk'(t)-d=-t^2T'(t)<0.
\]
At any positive-height zero of $\mathscr I$, $\lambda=d+h\kappa$.
We claim $\kappa>T(h)$. If not, for every $h<u<U$,
$k(u)>d+u\kappa=\lambda+\kappa(u-h)$.
Thus
\[
 \delta(z)=\phi(h+z)-\phi(h)-\lambda z-\kappa z^2/2>0
 \quad(0<z<U-h),\qquad \delta(0)=0.
\]
Translate the base to zero. The actual arc $z_a$ and the fitted quadratic
arc $z_q$ have identical endpoints $A=M-r$, $B=M+r$, and are positive
inside. Their difference satisfies
\[
 w'=\left[\lambda+\tfrac\kappa2(z_a+z_q)\right]w+\delta(z_a),
 \qquad w(A)=w(B)=0.
\]
A positive integrating factor gives $w(B)>0$, a contradiction.
Only actual heights enter $\delta$; the model peak is allowed to exceed
$U-h$. This proves $\kappa>T(h)$ and hence $\lambda-k(h)>0$.
The transport identities yield
\begin{equation}\label{cr:intercept-transport}
 D\mathscr I=\frac{S+hR}{\mathcal J}(\lambda-k(h)).
\end{equation}
At a zero this is strictly positive, so the outward derivative
$-D\mathscr I/\alpha$ is negative. A first zero reached from the initial
negative values would have nonnegative derivative. This contradiction
proves the claim along every actual arc. Reflection of $\phi,d$
proves the reversed case.
\end{proof}

\begin{lemma}\label{cr:intercept-axis}
Assume additionally that $\phi(0)=0$ continuously and $k(u)\to-\infty$
as $u\to0^+$. In the positive case of
Lemma~\ref{cr:intercept-positive}, every transverse axis fit satisfies
$\lambda(0,r)<d$. If $T'<0$ and $k(u)\to+\infty$, then $\lambda(0,r)>d$.
\end{lemma}
\begin{proof}
Fix a peak with axis half-width $r_0$. The data
$(M(h,r(h)),M_r(h,r(h)),r(h))$ converge to an interior point of the
model fitting domain by Lemma~\ref{cr:arcs}.
The joint inverse in Lemma~\ref{cr:fit} gives finite limits of
$\lambda$ and $\kappa$; in particular $h\kappa\to0$.
Since $k(h)\to-\infty$, $\lambda-k(h)>0$ for $0<h\le h_1$ with
$h_1$ sufficiently small. Along this fixed-peak arc,
\eqref{cr:intercept-transport} gives $d\mathscr I/dh>0$ there.
Thus $\mathscr I(h)\le\mathscr I(h_1)=-\eta<0$ for $0<h<h_1$,
where $\eta$ is independent of $h$ for this arc. Let $h\to0^+$
to obtain $\lambda(0,r_0)-d\le-\eta<0$.
Reversing the signs proves the other case.
\end{proof}

\begin{proof}[Proof of Theorem~\ref{cr:intercept}]
The derivative in \eqref{cr:intercept-condition} equals
$\mathcal I_d/(G^2g^2)$, whose denominator is positive on both sides.
For the profiles,
\[
 \frac d{du}\frac{k_\pm(u)-d}{u}
 =\left.\frac1{g(x)}\frac d{dx}
 \frac{f(x)/g(x)-d}{G(x)}\right|_{x=x_\pm(u)}.
\]
It is positive on the right and negative on the left.
Because $f(x)\to F'(0)<0$ and $g(x)\to0$ with the sign of $x$,
$k_+(u)\to-\infty$ and $k_-(u)\to+\infty$.
These limits hold for degenerate restoring forces as well. The preceding
two lemmas prove
\begin{equation}\label{cr:intercept-separation}
 \lambda_+(0,r)<d<\lambda_-(0,r)\qquad(r\in J).
\end{equation}
The original smooth system has unique equilibrium solutions, so matching
applies. At a zero of $\Delta$, the two model midpoints agree.
Since $H_\lambda,H_\kappa>0$, \eqref{cr:intercept-separation} forces
$\kappa_+>\kappa_-$: otherwise increasing first $\lambda_+$ to
$\lambda_-$ and then $\kappa_+$ to $\kappa_-$ would strictly increase
an equal midpoint. The conditional derivative
$\widehat v_\kappa=\mathcal J/P>0$ now gives $v_+>v_-$, hence
$\Delta'>0$ and $0<\Pi'<1$. Lemma~\ref{cr:count} proves uniqueness.
The transformation $(t,y,F,d)\mapsto(-t,-y,-F,-d)$ proves the repelling
case. Since each comparison takes place below a finite peak, the two
potential intervals may have finite, unequal endpoints.
\end{proof}

\begin{remark}
Separate constants $d_+,d_-$ suffice only if they can be ordered as
$d_+\le d_-$. In that case any $d\in[d_+,d_-]$ is already a common
constant, since $((k-d)/u)'$ increases with $d$. Unordered separate
bounds do not imply \eqref{cr:intercept-separation}. Likewise the theorem
is an upper-bound result: $F=-x$, $g=x$ satisfies its strict condition
with $d=0$, while $(x^2+y^2)'/2=x^2$ excludes all nonconstant periodic orbits.
\end{remark}

\section{Proof of uniqueness with slope interfaces}\label{cr:sec-interface}
We now allow jumps in the profile slopes. The main step is to prove
continuous differentiability of the endpoint data despite these jumps;
the resulting variation formulas then extend the comparison argument.

\begin{lemma}\label{cr:bv-arcs}
Let $\phi\in C([0,U))\cap\operatorname{Lip}_{\mathrm{loc}}((0,U))$,
$\phi(0)=0$, and suppose its almost-everywhere derivative has a
left-continuous locally BV representative $k$.
The endpoint, half-width, strict first-variation, and axis-limit assertions of
Lemma~\ref{cr:arcs} still hold. The relative midpoint is locally Lipschitz
in $(h,r)$, and $M_r$ is continuous; its axis value is $C^1$ in $r$.
Along a fixed-peak arc, $h,L$ are $C^1$, $\theta$ is locally absolutely
continuous, and \eqref{cr:angle-ode} holds with its middle equation
interpreted almost everywhere. The initial terms in \eqref{cr:angle-start}
use $k(t-)$, and the angular integral remains valid.
\end{lemma}
\begin{proof}
On a compact positive-height interval, local Lipschitz continuity gives
unique scalar solutions. If $p=u_y=\phi(u)-y$ vanishes at $y_0$, then
$u(y)-u(y_0)=o(|y-y_0|)$, so
$p(y)=-(y-y_0)+o(|y-y_0|)$. Every critical point is therefore a strict
nondegenerate maximum; each level is visited at most twice.
The bounded-profile argument in Lemma~\ref{cr:arcs} gives finite
transverse endpoints.

Fix an initial coordinate $y_0$ and initial height $\eta$.
Nonintersecting solutions give difference quotients governed by the
bounded coefficient
\[
 c_d(y)=\frac{\phi(u(y;\eta+d))-\phi(u(y;\eta))}
 {u(y;\eta+d)-u(y;\eta)}.
\]
A locally BV function has at most countably many discontinuities; at
every continuity point of $k$, $\phi$ is differentiable with derivative
$k$. The reference solution visits the exceptional heights at most
twice each. Thus $c_d\to k(u(y;\eta))$ almost everywhere.
Dominated convergence in the exponential formula for the difference
quotient gives
\begin{equation}\label{cr:bv-variation}
 u_\eta(y;\eta)=\exp\!\left(\int_{y_0}^y k(u(s;\eta))\,ds\right)>0.
\end{equation}
The convergence is uniform on compact $y$ intervals because its integral
error is bounded by an $L^1$ error. If $\eta$ varies, dominated convergence
away from those same exceptional visits proves continuity of this derivative.
Together with continuous $u_y$, this yields a $C^1$ solution family at a
fixed initial coordinate.

The peak coordinate $\phi(t)$ may not be differentiable. Instead,
fix $t_0$ and $y_0=\phi(t_0)$. When the peak becomes $t_0+d$, its initial
coordinate moves by $O(d)$. The integral equation and a local Lipschitz
bound give $|u(y;t_0+d)-(t_0+d)|\le C|y-\phi(t_0+d)|^2$ near the peak;
this follows by applying Gronwall's inequality to
$|u_y|\le L|u-(t_0+d)|+|y-\phi(t_0+d)|$.
Consequently $u(y_0;t_0+d)=t_0+d+O(d^2)$.
Using \eqref{cr:bv-variation} from this fixed coordinate gives the
continuous peak variation
\[
 u_t(y;t)=\exp\!\left(\int_{\phi(t)}^y k(u(s;t))\,ds\right)>0.
\]
To prove continuity in $t$, note that the contribution of the moving
lower integration endpoint is at most a local bound for $|k|$ times
$|\phi(t)-\phi(t_0)|$. On a common
fixed integration interval, the integrands converge away from the
countably many times at which the limiting orbit visits a discontinuity
height of $k$; dominated convergence applies. Thus the peak derivative
is continuous even when the peak itself crosses an interface.
The implicit function theorem at transverse endpoints proves that $A,B$
are jointly $C^1$, including when the peak lies at a slope interface.
Their derivatives and the formula for $v$ are those of
Lemma~\ref{cr:arcs}. Subtracting $\phi(h)$ makes $M$ only locally
Lipschitz in height, but does not alter its continuous derivative with respect to half-width.
The inverse equation $Y_u=(\phi(u)-Y)^{-1}$ near nonzero axis endpoints
has the same bounded smooth dependence on $Y$ as before, so the complete
axis-extension argument of Lemma~\ref{cr:arcs} applies unchanged.

Differentiating the endpoint identities in base height almost everywhere
gives the first two angular equations. For the third, the inverse equation
has a continuous height coefficient even at an interface; hence
$(A_t)_h=A_t/(\phi(h)-A)^2$ and
$(B_t)_h=B_t/(\phi(h)-B)^2$ at every such height. Thus
$L'=2M/(s^2-M^2)$ continuously. Both branches near a peak lie below it;
using the left trace in
$\phi(u)-\phi(t)=k(t-)(u-t)+o(|u-t|)$ gives
$u=t-\eta^2/2-k(t-)\eta^3/6+o(\eta^3)$.
The inversion argument of Lemma~\ref{cr:angles} gives the stated initial
terms and integral, with a pointwise remainder at each peak.
\end{proof}

\begin{proof}[Proof of Theorem~\ref{cr:interface}]
A finite convex or concave function is locally Lipschitz, and its left
derivative is a left-continuous locally BV representative of its
a.e. derivative. Strict convexity of $\phi-Ku^2/2$ makes $k-Ku$
strictly increasing: if its values agreed at two heights, monotonicity
would make the derivative constant almost everywhere between them,
forcing an affine segment. Conversely strict derivative ordering gives
strict convexity by integration. Thus Lemma~\ref{cr:bv-arcs} applies to
both profiles, with $k_+-Ku$ increasing and $k_--Ku$ decreasing.

To compare the profiles across jumps, we use the slope measure. With
the reference field of Lemma~\ref{cr:reference}, set
$E=k(h)-J_K(s,\theta)$, $a=(J_K)_\theta>0$ and $\eta=k(h)-Kh$.
Since $h$ decreases and $k$ is left-continuous, $\eta$ is right-continuous
in $s$. At zero half-width set $\eta(0)=k(t-)-Kt$; there is no atom there,
even if the peak is a slope interface. The angular equations give
\[
 dE+aE\,ds=d\eta,
 \qquad z'=E/(\Theta_K)_z\quad\hbox{almost everywhere}.
\]
In this measure equation, a height jump
$\Delta k(c)=k(c+)-k(c-)$ at $0<c<t$ contributes the atom
$-\Delta k(c)$ when the descending base crosses $c$.
The initial expansions give $E(0+)=0$ as before. The integrating-factor
formula is exactly \eqref{cr:error-kernel}, now for a possibly discontinuous
finite-variation measure. Its kernel is positive, bounded by one, and
has a positive lower bound on each interior compact interval.
Strict increase of $k-Ku$ makes $E<0$ and the absolutely continuous label
strictly decreasing; strict decrease reverses the signs. We recover the
strict terminal-angle comparison \eqref{cr:angle-compare}.

We next establish geometric uniqueness for the original system.
Off the axis the scalar solution families in Lemma~\ref{cr:bv-arcs}
give $C^1$ orbit coordinates, and the original time integral has a
continuous nonzero denominator. At nonzero axis points
\eqref{cr:axis-graph} still applies. Near zero, monotonicity of the slopes
bounds $f=F'$ above almost everywhere, as in Lemma~\ref{cr:flow}.
Integrate this inequality on intervals away from zero and then use
continuity of $F$ to reach zero. This gives $xF(x)\le Lx^2$ and hence
backward uniqueness by the energy proof of that lemma. The reversed
case gives forward uniqueness. The matching proof therefore applies, including the $C^1$ endpoint multiplier and
the correspondence between isolated zeros and limit cycles.

At a matching half-width $R$, terminal reference labels agree, and
$\theta_+(s)>\Theta_K(s,z_R)>\theta_-(s)$ on $(0,R)$.
The angular integrals give $v_+>v_-$; thus $\Delta'(R)>0$ and
$0<\Pi'<1$. Lemma~\ref{cr:count} proves uniqueness and attraction.
All signs reverse in the other case. Since the vector field is continuous,
the endpoint variations give the transverse multiplier across the interfaces.
\end{proof}

\begin{proof}[Proof of Corollary~\ref{cr:interface-piecewise}]
A compact interval away from zero meets finitely many interfaces.
Strict monotonicity on its pieces and the directed jumps make a
one-sided representative of $q-KG$ strictly monotone on the whole
interval. The positive inverse potential branch preserves this order
and the negative branch reverses it. Integrating the corresponding
profile derivatives gives the strict convexity and concavity in
\eqref{cr:interface-condition}; Theorem~\ref{cr:interface} applies.
\end{proof}

\begin{remark}
If a continuous profile has a jump $j$ in its first derivative at $u=c$,
its distributional third derivative contains $j\delta'_c$, not a
nonnegative measure even when $j>0$. The following example shows that
nonnegative third derivatives on the pieces and upward slope jumps on the positive side, with reversed signs
on the negative side, can produce three limit cycles. Set
\[
 \psi(u)=u-u^2+10(u-3/2)_+,\qquad
 P(x)=\operatorname{sgn}(x)\psi(x^2/2),\qquad P(0)=0,
\]
and consider $\dot x=y-\varepsilon P(x)$, $\dot y=-x$.
Its two profiles are $\pm\varepsilon\psi$, with exactly those piecewise
third-derivative signs and directed jumps. Nevertheless it has at least
three hyperbolic limit cycles for every sufficiently small
$\varepsilon>0$. Indeed the first radial return coefficient is
\[
 I(r)=\int_0^{2\pi}\cos\theta\,P(r\cos\theta)\,d\theta
 =\frac{4r^2}{3}A(r^2),\qquad
 A(z)=1-\frac25z+10(1-3/z)_+^{3/2}.
\]
For the hinge term, the substitution $w=\sin\theta$ on a quadrant gives
$\int_0^{2\pi}|\cos\theta|(r^2\cos^2\theta/2-c)_+\,d\theta
=4r^2(1-2c/r^2)_+^{3/2}/3$; the polynomial terms give the other two
coefficients. There is a simple root $z_1=5/2$. For $z>3$,
\[
 \begin{gathered}
 A(3)=-1/5,\qquad A(6)=5/\sqrt2-7/5>0,\qquad A(z)\to-\infty,\\
 A'(z)=-2/5+B(z),\qquad B(z)=45z^{-2}\sqrt{1-3/z}.
 \end{gathered}
\]
Since $B'/B=(-4z+15)/(2z(z-3))$, $B$ first increases from zero to
$B(15/4)=16/(5\sqrt5)>2/5$ and then decreases to zero.
Thus $A$ has one minimum below zero and one maximum above zero on
$(3,\infty)$, and exactly two simple roots there.
All three corresponding circles cross every encountered switching line
$x=0,\pm\sqrt3$ transversely. Piecewise smooth flows and their transverse
hitting times therefore compose to a locally smooth return map near each
circle, with
$\Pi_\varepsilon(r)=r-\varepsilon I(r)+O(\varepsilon^2)$ in $C^1$.
The implicit function theorem at the three simple roots gives three
fixed points and multipliers
$1-\varepsilon I'(r_i)+O(\varepsilon^2)\ne1$.
Thus the piecewise third-derivative signs and directed slope jumps do
not imply a two-cycle bound.
\end{remark}

\section{A joint weak-damping criterion with an explicit remainder}
\label{jt:proof-section}

The proof of Theorem~\ref{jt:weak} uses the first-order approximation
of the fitted-curvature difference in the damping parameter. We bound
the remainder uniformly and confine the half-widths of all cycles to a fixed
interval. All differentiations in a model half-width hold the model parameters
fixed unless a total derivative
is explicitly indicated. For the fitted curvatures of the two actual
half arcs put
$\mathcal Q(r)=\kappa_+(r)-\kappa_-(r)$.
The fitting and matching lemmas proved above imply
\begin{equation}\label{jt:scalar}
 \Delta'=A\Delta+B\mathcal Q,\qquad B>0,
 \qquad \Delta=M_+-M_-.
\end{equation}
To express its coefficients, solve $H(\lambda,\kappa,r)=m$ for $\lambda$
and set
$V(r,m,\kappa)=H_r(\lambda,\kappa,r)$. With
$P=H_\lambda$, $Q_m=H_\kappa$, $R_m=H_{\lambda r}$,
$S=H_{\kappa r}$ and $J=PS-Q_mR_m$, one has
$V_m=R_m/P>0$ and $V_\kappa=J/P>0$.
Integrating first in $m$ and then in $\kappa$ between the two actual data
gives \eqref{jt:scalar}, even when $\Delta$ or $\mathcal Q$ vanishes.
Multiplication by $\exp(-\int A)$ reduces counting to the sign changes of
$\mathcal Q$. If $\mathcal Q$ is strictly monotone, the transformed
displacement is strictly monotone on each of at most two intervals and has
at most two zeros. For a weak one-sign-change assumption the same proof
bounds isolated zeros, allowing zero intervals.

\subsection{A uniform analytic estimate}
\begin{lemma}\label{jt:analytic}
Under the holomorphic bound on $p$ in Theorem~\ref{jt:weak}, the auxiliary
system $\dot z=y-\varepsilon p(z)$, $\dot y=-z$ has transverse half-arc
data on the whole interval $\mathcal B=[1/16,4]$ for sufficiently small
real $\varepsilon$. Its curvature difference satisfies
\begin{equation}\label{jt:remainder}
 |\mathcal Q'(r,\varepsilon)-\varepsilon C_p'(r)|
 <10^{154}|\varepsilon|^2,\qquad
 r\in\mathcal B,\quad |\varepsilon|\le5\cdot10^{-101}.
\end{equation}
The auxiliary arcs may extend beyond the original physical strip; any
physical arc of such a half-width coincides with the corresponding auxiliary arc.
\end{lemma}
\begin{proof}
For the explicit parameter range in Theorem~\ref{jt:weak}, consider
the auxiliary right-half-plane damping
$F_\sigma(z)=\varepsilon p(\sigma z)$, $\sigma=\pm1$. In polar variables
$z=r\rho\sin\theta$, $y=r\rho\cos\theta$, its equation is
\begin{equation}\label{jt:polar}
 \rho_\theta=-\frac{\rho\sin\theta\,Z}{1-\cos\theta\,Z},
 \qquad Z=\frac{\varepsilon p(\sigma r\rho\sin\theta)}{r\rho},
 \qquad 0\le\theta\le\pi.
\end{equation}
The choice $p(\sigma z)$, with no additional factor $\sigma$, represents
the profile $\varepsilon p(\sigma\sqrt{2u})$ on either side.
Let $T(s)$ be the terminal value for $\rho(0)=s$.
The endpoints are $rs$ and $-rT(s)$. Prescribing half-width $r$ means
$(s+T(s))/2=1$, and the normalized midpoint is $m=M/r=s-1$.

Complexify $r,s,\varepsilon$, keeping $\theta$ real. Put
\[
 D_j=\{r:\operatorname{dist}(r,\mathcal B)<1/(2^{7-j})\},
 \quad j=1,2,\qquad L=10^{40},\quad E=10^{-100}.
\]
Thus $D_2$ has radius $1/32$, $D_1$ radius $1/64$,
and $1/32<|r|<5$ on $D_2$.
For $|\rho|<2$ the argument of $p$ is in $|w|<10$ and $|Z|\le L|\varepsilon|$.
If $|\varepsilon|<E$, $|s-1|<1/2$, continuation of
\eqref{jt:polar} gives
\[
 |\rho_\theta|\le4L|\varepsilon|,\qquad
 |T(s)-s|\le16L|\varepsilon|=:d.
\]
Indeed, the second bound leaves a distance greater than $1/2-d$ from
$|\rho|=2$, closing the first-exit argument. Picard iteration gives
holomorphic dependence throughout the finite integration interval.
On $|s-1|\le1/4$, Cauchy's estimate yields
$|\partial_s(T-s)|\le4d$. The equation
$s=1-(T(s)-s)/2$ is a contraction of that disk, with constant $2d<1$.
Its fixed point is holomorphic and
\begin{equation}\label{jt:data}
 |m_\sigma|\le8L|\varepsilon|\quad(D_2),\qquad
 \max(|m_\sigma|,|v_\sigma|)<10^{44}|\varepsilon|\quad(D_1),
\end{equation}
where $v_\sigma=M_\sigma'=m_\sigma+r\partial_rm_\sigma$.
The latter estimate follows from the distance $1/64$ between the two tube
boundaries: $|v_\sigma|\le321\cdot8L|\varepsilon|$.
For real data $\rho>0$ and the angular speed is positive, so these solutions
are the required transverse half arcs.

We now construct and estimate the model inverse on a complex neighborhood.
The quadratic potential profile corresponds to
$F(z)=\lambda z^2/2+\kappa z^4/8$. Put $a=\lambda r$ and $b=\kappa r^3$.
Then $Z=a\rho\sin^2\theta/2+b\rho^3\sin^4\theta/8$.
For $|a|,|b|<10^{-3}$, $|Z|<2\cdot10^{-3}$ while $|\rho|<2$.
The preceding continuation and contraction argument gives a holomorphic
normalized midpoint $m(a,b)$ with $|m|<1/100$; explicitly
\[
 |T(s)-s|\le\frac{4\pi\cdot10^{-3}}{1-2\cdot10^{-3}}<\frac1{50}.
\]
The normalized fitting map is
\[
 \Psi(a,b)=(m, m+am_a+3bm_b).
\]
Cauchy's estimate on $|a|,|b|<10^{-4}$ gives
$|m_a|,|m_b|\le100/9$, and consequently $\|\Psi\|_\infty<1$.
On $|a|,|b|\le5\cdot10^{-5}$ every second partial derivative of each
component is smaller than $10^9$: the pure and mixed estimates are
$8\cdot10^8$ and $4\cdot10^8$, respectively. Linearizing
\eqref{jt:polar} at zero gives
\[
 J_0=D\Psi(0,0)=
 \begin{pmatrix}1/3&1/15\\2/3&4/15\end{pmatrix},\qquad
 J_0^{-1}=\begin{pmatrix}6&-3/2\\-15&15/2\end{pmatrix}.
\]
For $\delta=10^{-14}$,
$\|D\Psi-J_0\|_\infty\le4\cdot10^9\delta=4\cdot10^{-5}$.
The fixed-point equation
\[
 \zeta=J_0^{-1}w-J_0^{-1}(\Psi(\zeta)-J_0\zeta)
\]
therefore has contraction constant at most $9\cdot10^{-4}$ on
$\|\zeta\|_\infty\le\delta$. If $\|w\|_\infty\le10^{-16}$ it maps
this polydisk into itself, since
$(0.225+0.0009)\delta<\delta$. Its unique holomorphic solution satisfies
\[
 \|\zeta\|_\infty\le
 \frac{22.5}{1-0.0009}\|w\|_\infty<23\|w\|_\infty.
\]
By \eqref{jt:data}, actual targets have norm below $10^{-56}$ for
$|\varepsilon|<E$. Thus $a_\sigma,b_\sigma$ are holomorphic on
$D_1\times\{|\varepsilon|<E\}$ and
$|b_\sigma|<10^{46}|\varepsilon|$. For real parameters this local inverse
equals the global inverse of Lemma~\ref{cr:fit}, by uniqueness.

Expand $b_\sigma=\varepsilon b_{\sigma,1}+\mathcal R_\sigma$.
The Cauchy series bound, for $|\varepsilon|\le E/2$, is
\[
 |\mathcal R_\sigma|\le2\cdot10^{146}|\varepsilon|^2.
\]
Since $\mathcal Q=(b_+-b_-)/r^3$ and $|r|>1/32$, its remainder has
modulus below $2^{17}10^{146}|\varepsilon|^2$ on $D_1$.
A second Cauchy estimate, on a disk of radius $1/64$ about any real
$r\in\mathcal B$, bounds the derivative remainder by
$2^{23}10^{146}|\varepsilon|^2<10^{154}|\varepsilon|^2$.
For each open disk, apply Cauchy's estimate at smaller radii and then
take the supremum.

It remains to identify the linear term. Differentiating the normalization
at $\varepsilon=0$ gives
\[
 m_{\sigma,1}(r)=\frac1{2r}\int_0^\pi\sin\theta\,
                 p(\sigma r\sin\theta)\,d\theta.
\]
Hence the difference of the first-order normalized midpoints is
$d_p=D_p/r$, and the difference of the first-order $v$ data is
$d_p+rd_p'$. The second row of $J_0^{-1}$ gives
\[
 \mathcal Q_1=\frac{15}{2r^3}(rd_p'-d_p)
             =\frac{15}{2r^3}(D_p'-2D_p/r)=C_p.
\]
This proves \eqref{jt:remainder}. The even part cancels in the linear
term, and its nonlinear contribution is included in the remainder estimate.
\end{proof}

\subsection{Localization of every possible cycle}
\begin{proof}[Proof of Theorem~\ref{jt:weak}]
Writing $g(x)=xq(x)$ near zero gives
$2G(x)=x^2\ell(x)$ with smooth positive $\ell$ and
$\ell(0)=g'(0)$. Thus $z=x\sqrt{\ell(x)}$ extends smoothly through zero.
It has positive derivative $z'=g/z$ off zero and $z'(0)=\sqrt{g'(0)}$.
The transformation $z=z(x)$, $d\tau=z'(x)\,dt$ is a diffeomorphism and
a positive time change, reducing the system to
\[
 z_\tau=y-\varepsilon p(z),\qquad y_\tau=-z,\qquad -a<z<b.
\]
On any compact periodic orbit the time factor has positive finite bounds;
periodicity, isolation and the return multiplier are preserved in both
directions.

Every such orbit lies in $|z|<3$ and has
$|\varepsilon p(z)|\le3\cdot10^{40}\varepsilon=:\Lambda<1$.
For either half arc its peak height is $t<9/2$. Its axis endpoints
$A<0<B$ obey $B\le\Lambda+\sqrt{2t}$ and
$-A\le\Lambda+\sqrt{2t}$. Indeed the peak ordinate is
$\phi(t)\in[-\Lambda,\Lambda]$. If $B>\Lambda$, the descending branch
therefore contains the entire ordinate interval $[\Lambda,B]$.
Integrating $u_y\le\Lambda-y$ on this interval, with $u(B)=0$, gives
$(B-\Lambda)^2/2\le u(\Lambda)\le t$; the case $B\le\Lambda$ needs
no estimate. Reflection gives the lower-endpoint bound. Consequently
$r\le\Lambda+\sqrt{2t}<4$.
The peak-height bound in Lemma~\ref{cr:arcs} gives $2t_\pm\le r^2$.
If $r\le1/16$, the whole orbit lies where $zp(z)>0$, whereas
\[
 \frac d{d\tau}(z^2+y^2)/2=-\varepsilon zp(z)<0
\]
off the axis. This contradicts periodicity. All possible cycles therefore
have $1/16<r<4$.

Lemma~\ref{jt:analytic} implies
$\mathcal Q'\ge\varepsilon(1-10^{-56})>0$ in the increasing case;
the decreasing case is identical with reversed signs.
The common domain of the actual half-width parameter is an interval by
Lemma~\ref{cr:arcs}; its intersection with $(1/16,4)$ remains an interval.
Equation~\eqref{jt:scalar} and the matching correspondence give the
bound for cycles with half-widths in this interval. Every cycle has such a
half-width, and the bound includes isolated nonhyperbolic cycles.
\end{proof}

\section{Compensation under two-derivative perturbations}
\label{jt:compact-proof}

Theorem~\ref{jt:compact} follows from stability of the endpoint data and
their first two derivatives with respect to half-width. We establish a strict curvature bound
for the reference profiles and then estimate its change under a $C^2$
perturbation.

\begin{lemma}\label{jt:strict-axis}
If $\phi'''>0$ at positive height and $\phi'\in L^1(0,\eta)$, then
$\kappa_r>0$ at every transverse axis half-width. Reversing the signs reverses
the conclusion.
\end{lemma}
\begin{proof}
Follow the fixed-peak characteristic from positive base height $h_1$ down
to zero. In the notation of the shape calculation in
Lemma~\ref{cr:shape}, $w=\kappa_r$ initially satisfies $w(h_1)>0$ and
$e=\lambda-\phi'<0$. Its exact transport equation is
\[
 Dw=Tw-C_r e,\qquad C=R_m/J,\qquad
 T=\alpha_r+C Q_m/P-(C_\kappa-Q_m C_\lambda/P)e.
\]
Here $D=\partial_h-\alpha\partial_r$, and $C_r$ is a model partial
derivative. Transverse-axis convergence keeps the model parameters and
their coefficients in a compact subset of the inverse-model domain.
Consequently $|T|\le C_0(1+|\phi'|)$, which is integrable.
Lemma~\ref{cr:model-inequality} gives $C_r\le0$.
After reversing height orientation, the source $C_r e$ is nonnegative.
An integrating factor yields
\[
 w(0)\ge w(h_1)\exp\left(-\int_0^{h_1}|T(h)|\,dh\right)>0.
\]
Smooth transverse-axis dependence gives the value $w(0)$.
Reflection proves the opposite case.
\end{proof}

\begin{lemma}\label{jt:jets}
Fix a reference profile $\phi_0$, a compact half-width interval
$[r_0,R]\Subset I_{\phi_0}$, and $a>0$. Take a height $T>a$ above all
reference peaks on a neighborhood of this interval. For perturbations
$\psi$ vanishing on $[0,a]$ and defined on a common height neighborhood,
there are $C,\epsilon_0>0$ such that
\begin{equation}\label{jt:jet-bound}
 \sup_{[r_0,R]}\|(M_\psi,M_\psi',M_\psi'')-(M_0,M_0',M_0'')\|_\infty
 \le C\|\psi\|_{C^2([a/2,T])}
\end{equation}
whenever the norm is below $\epsilon_0$. The perturbed transverse arcs
then exist at all these half-widths. The constants can be bounded using the
reference equations and transverse-root estimates given below.
\end{lemma}
\begin{proof}
Label peak height by $t$, write $z=y-\phi(t)$, and solve
\[
 U_z=\phi(U)-\phi(t)-z,\qquad U(0,t)=t.
\]
The first two peak derivatives $V=U_t$, $W=U_{tt}$ satisfy
\begin{align}
 V_z&=\phi'(U)V-\phi'(t),& V(0,t)&=1,\label{jt:jets-eq}\\
 W_z&=\phi'(U)W+\phi''(U)V^2-\phi''(t),&W(0,t)&=0.\nonumber
\end{align}
Choose a cut height $d<a/2$ below all reference peaks, at uniformly
transverse crossings. Cover the compact family between its peaks and $d$
by finitely many tubes for $(U,V,W)$. If $L_0$ bounds the Lipschitz
constant of the reference extended field and $E_0$ bounds the coefficient
of its perturbation on these tubes, Gronwall's inequality gives
\[
 \|Z_\psi(z)-Z_0(z)\|\le
 E_0\|\psi\|_{C^2}|z|e^{L_0|z|}.
\]
Choosing this smaller than the tube clearance validates the estimate by
a first-exit argument. The estimate uses the Lipschitz constant of the
reference field and the $C^2$ bound on the difference of the fields at
the same point.

At $U(z_e(t),t)=d$, let $f_e=\phi(d)-\phi(t)-z_e$. Implicit
differentiation gives
\begin{align*}
 z_e'&=-V/f_e,\\
 z_e''&=-\{W+2(\phi'(d)V-\phi'(t))z_e'
              +(\phi'(d)f_e-1)(z_e')^2\}/f_e.
\end{align*}
The physical endpoint is $Y=\phi(t)+z_e$.
Positive lower bounds for $|f_e|$ control these formulas and shifted-root
evaluation. Below $d$ the perturbation vanishes. The inverse branch obeys
\begin{align*}
 Y_u&=(\phi_0(u)-Y)^{-1},\\
 (Y_t)_u&=(\phi_0(u)-Y)^{-2}Y_t,\\
 (Y_{tt})_u&=(\phi_0(u)-Y)^{-2}Y_{tt}
                  +2(\phi_0(u)-Y)^{-3}Y_t^2.
\end{align*}
Uniform transverse rectangles bound every denominator away from zero
and propagate the three endpoint jets to $u=0$, where continuity of
$\phi_0$ suffices.
Write $m(t)=(B(t)+A(t))/2$, $r(t)=(B(t)-A(t))/2$. Then
\[
 M=m,\qquad M_r=m_t/r_t,\qquad
 M_{rr}=(m_{tt}r_t-m_t r_{tt})/r_t^3.
\]
The positive reference lower bound for $r_t$ persists. Enlarge the peak
interval slightly so its half-width image contains $[r_0,R]$ with positive
clearance. Monotonicity and the intermediate value theorem then give
perturbed inverse peaks at each target half-width. Reference third endpoint
jets bound the change from the reference peak to this inverse peak.
Together these estimates prove \eqref{jt:jet-bound}; they also determine
bounds for $C,\epsilon_0$ in finitely many steps.
\end{proof}

\begin{proof}[Proof of Theorem~\ref{jt:compact}]
Differentiation of the fitting equations gives
\[
 P\lambda_r+Q_m\kappa_r=0,\qquad
 R_m\lambda_r+S\kappa_r=M_{rr}-H_{rr},
\]
so that
$\kappa_r=(P/J)(M_{rr}-H_{rr})=\mathscr T(M,M_r,M_{rr},r)$.
The reference difference is strictly positive by
Lemma~\ref{jt:strict-axis}; compactness permits the positive margin $m$.
On $[r_0,R]$, Lemma~\ref{jt:jets} and the mean value inequality on the
chosen tubes imply
\[
 \mathcal Q'\ge m-L_+C_+\delta_+-L_-C_-\delta_->0.
\]
For $0<r<r_0$, the peak-height bound bounds every peak by
$r^2/2<a$. The candidate and reference profiles coincide there, so their
arcs and fitted curvatures coincide. Conversely, every reference arc in
this range is a candidate arc. Thus the whole interval $(0,R]$ is a
common transverse interval and $\mathcal Q'>0$ throughout. Equation
\eqref{jt:scalar} proves the stated count, including an isolated root at
$R$. If all cycles lie in this half-width interval, the bound is global.
\end{proof}

A rapidly oscillating perturbation illustrates the advantage of the
$C^2$ estimate. For a smooth cutoff $\chi$ supported in a compact
positive-height interval,
take $\psi_n(u)=\epsilon n^{-5/2}\chi(u)\sin(nu)$.
Its $C^2$ norm tends to zero, while on a plateau of $\chi$ its third
derivative is $-\epsilon n^{1/2}\cos(nu)$ and has both signs with
unbounded size. For all sufficiently large $n$, either or both profile
sign tests may fail while Theorem~\ref{jt:compact} applies. To determine
a numerical threshold for $n$, one must bound the reference constants
in the theorem.

\section{Branch curvature and the direction of neutral cycles}\label{sk:neutral-proof}
The sign of the second return derivative in Proposition~\ref{sk:neutral}
is determined by comparing transverse variations along the two sides of
the neutral orbit. We construct the return map from the corresponding
half-return maps.

\begin{proof}[Proof of Proposition~\ref{sk:neutral}]
Fix a neutral periodic orbit. It intersects each half of the section
$x=0$ exactly once, at $q>0$ and $-m<0$. Indeed, its extrema of $x$
satisfy $\ddot x=-g(x)$; thus every maximum is positive and every
minimum is negative. The crossings of the positive section all have
the same direction, whereas crossings of a Jordan curve along a ray
alternate. The intersections are finite because they are transverse
and the orbit is compact. The same argument applies to the negative
section. There is exactly one turning point in each open half-plane,
since every critical point of $x$ in the right half-plane is a strict
maximum and every such point on the left is a strict minimum.

Write the velocities on the upper and lower arcs as $P(x)>0$ and
$-N(x)<0$, respectively, between turning points $x_-<0<x_+$.
Then
\begin{equation}\label{sk:velocity}
 P_x=-f-g/P,\quad N_x=f-g/N,\quad
 (P^2-N^2)_x=-2f(P+N).
\end{equation}
The identity $(\dot x^2/2+G)'=-f\dot x^2$ forces $x_+>x_f$.
The last expression in \eqref{sk:velocity} changes sign once, from
positive to negative, and $P^2-N^2$ vanishes at both turning points.
Consequently $P>N$ in the interior.

Starting at $(0,q)$, let $E=q\exp(-\int f\,dt)$ along the orbit.
Neutrality gives $\oint f\,dt=0$, so $E$ closes periodically. Its upper
and lower restrictions $E_+,E_-$ agree at each turning point and obey
\[
 (E_+)_x=-fE_+/P,\qquad (E_-)_x=fE_-/N.
\]
The same single-sign-change argument applied to
$(\log(E_+/E_-))_x=-f(1/P+1/N)$ proves $E_+>E_->0$.
At the section write $E_+(0)=q$, $E_-(0)=e$.

Set
\[
 V=P+N,\quad W=P-N,\quad t=PN,\quad z=W/V,\quad V_0=q+m.
\]
Here $t$ denotes the velocity product $PN$. From
\eqref{sk:velocity},
\begin{equation}\label{sk:V}
 V_x=-gV/t,\qquad h_V=-t/V,\qquad t=V^2(1-z^2)/4.
\end{equation}
Thus on each side $V$ increases from zero at its turning point to
$V_0$ at the section. Transverse variation with respect to $q$, at
fixed $x$, gives
\[
 \delta P=E_+/P,\quad\delta N=E_-/N,\quad
 U:=\delta V=E_+/P+E_-/N,\quad
 U_x=g(E_+/P^3+E_-/N^3).
\]
Let
\begin{equation}\label{sk:HC}
 U_0=1+e/m,\quad \alpha=U_0/V_0,\quad
 H=(NE_++PE_-)/V=tU/V,\quad C=H-\alpha t,\quad K_0=\alpha qm.
\end{equation}
At fixed normalized velocity $s=V/V_0$, the identity
$V_x\delta x+U=sU_0$ yields $\delta x=C/g$.

To compute the second return derivative, launch nearby points
$(0,\eta)$ forward through the right side and backward through the
left side. Let $-m_i(\eta)$ be their lower endpoints and $J_i(\eta)$
their damping integrals, with the left integral taken in backward
physical time. Thus both integrations start at the positive section,
and the transverse-flow formula gives
\[
 m_i'(\eta)=\frac{\eta}{m_i(\eta)}e^{-J_i(\eta)}\quad(i=L,R).
\]
The forward first return and its damping integral are therefore
\[
 \Pi=m_L^{-1}\circ m_R,\qquad
 I(\eta):=\int_0^{T(\eta)}f\,dt=J_R(\eta)-J_L(\Pi(\eta)).
\]
At the neutral orbit $m_L=m_R=m$, $m_L'=m_R'=e/m$, and $J_L=J_R$.
For each nearby half-arc separately,
\[
 J_i=V_{i,0}\int_0^1L(x_i(s,\eta))\,ds,\qquad V_{i,0}=\eta+m_i(\eta).
\]
At the neutral orbit both scale derivatives are $U_0$ and
$J_R=J_L$; their contribution is
$U_0(J_R-J_L)/V_0=0$. The remaining derivatives use
$L_x\delta x=A C$. Hence
\begin{equation}\label{sk:Iprime}
 I'(q)=\int_0^{V_0}(A_RC_R-A_LC_L)\,dV,\qquad \Pi''(q)=-I'(q).
\end{equation}
The second identity follows by differentiating
$\Pi'(\eta)=\eta e^{-I(\eta)}/\Pi(\eta)$ at the neutral fixed point.
The integrals in \eqref{sk:Iprime} converge. If
$f_0=f(0)<0$ and $t_0=qm$, then at the section
\[
 V_0-V=\frac{g_1V_0}{2t_0}x^2+O(x^3),\quad
 U-U_0=O(x^2),\quad C=O(x^2),\quad A=O(x^{-3}).
\]
Uniformly for nearby initial points,
$x_i(s,\eta)=\pm\sqrt{2\eta m_i(\eta)/g_1}\sqrt{1-s}+O(1-s)$ and
$\partial_\eta x_i=O(\sqrt{1-s})$. Thus both $L$ and
$L_x\partial_\eta x_i$ are dominated by a multiple of
$(1-s)^{-1/2}$. At the nondegenerate turning points, the square-root
coordinate for the original flow makes $x_i$ and its variation continuously differentiable
in $V$. These estimates justify differentiation under the integral sign
in \eqref{sk:Iprime}.

Define the response variables
\[
 D_E=E_+-E_-,\quad\Lambda=\log(E_+/E_-),\quad d=D_E/H,\quad
 S=\Lambda-d,\quad \mathcal R=(D_E/V-LH)/K_0.
\]
Their section values are common to both sides and carry a subscript
zero. At the turning points $\Lambda,d,S$ tend to zero. Differentiating
\eqref{sk:velocity} and \eqref{sk:HC} gives
\begin{equation}\label{sk:response}
\begin{gathered}
 \Lambda_V=L,\quad L_V=-tA/V,\quad (D_E)_V=LH,\quad H_V=K_0z\mathcal R,\\
 \mathcal R_V+(1+VLz)\mathcal R/V=tAH/(VK_0),\quad
 S_V=dzK_0\mathcal R/H.
\end{gathered}
\end{equation}
In particular, with $\mathfrak I(V)=V\exp(\int_0^VLz\,dv)$,
\begin{equation}\label{sk:response-integral}
 \mathcal R(V)=\mathfrak I(V)^{-1}
 \int_0^V\mathfrak I(v)t(v)A(v)H(v)/(vK_0)\,dv.
\end{equation}
On the right $A>0$, so $\mathcal R_R>0$ and $S_R$ strictly increases;
in particular $S_0>0$. On the left $L>0$ and $\Xi>0$ imply
$A_h>A^2/L\ge0$, whence $A_L$ strictly decreases with $V$.
The integral in \eqref{sk:response-integral} therefore changes sign at
most once, from positive to negative. At the section $L_L\to\infty$
and $\mathcal R_L\to-\infty$. If it were never positive, $S_L$ could
not attain $S_0>0$. Thus $\mathcal R_L$ changes sign exactly once,
and $S_L$ strictly increases to a maximum $S_{\max}>S_0$ and then
strictly decreases to $S_0$. It stays positive throughout. This also
implies $A_{L,t}>0$ at the turning point; at the response zero $A_L<0$,
so that maximum is nondegenerate.

The algebraic relation between the response variables is
\begin{equation}\label{sk:response-relation}
 d=\frac{2\tanh(\Lambda/2)}{1-z\tanh(\Lambda/2)},\qquad
 z=\coth(\Lambda/2)-\frac2{\Lambda-S}.
\end{equation}
In particular $0<z<1$ and $0<d<\Lambda$. We now compare the two
branches using the strict curvature inequalities.
Where $L>0$, set
\[
 \chi=A/L^3,\quad \Psi=-g\chi_x/f^2,\quad
 Y=VL,\quad\eta=(1-z^2)Y.
\]
Changing the independent variable to $\Lambda$ gives
\begin{equation}\label{sk:Lambda}
 \Psi-2\chi^2=-\Xi/L^6,\quad
 Y_\Lambda=1-(1-z^2)\chi Y^2/4,\quad
 z_\Lambda=(1-z^2)/2-2z/Y,\quad\chi_\Lambda=\eta\Psi/4.
\end{equation}
Consequently
\begin{equation}\label{sk:trap}
 \mathcal T=z\chi Y^2+2Y-2\Lambda,\qquad
 \mathcal T_\Lambda=zY^2\eta(\Psi-2\chi^2)/4,
 \qquad \mathcal T_L<0<\mathcal T_R.
\end{equation}
Compare the entire left arc with the right outer segment $L_R>0$,
which covers $0<\Lambda\le\Lambda_0$ before the right arc reaches its
maximum of $\Lambda$. Put $\Delta z=z_L-z_R$, $\Delta Y=Y_L-Y_R$.
On the boundary $\Delta Y=0$, $\Delta z\ge0$, \eqref{sk:trap} gives
$z_L\chi_L<z_R\chi_R$. Since $\chi_R>0$, whether $\chi_L$ is positive
or nonpositive this implies
$(1-z_L^2)\chi_L<(1-z_R^2)\chi_R$, and hence $\Delta Y_\Lambda>0$.
On $\Delta z=0$,
\[
 \Delta z_\Lambda=2z\Delta Y/(Y_LY_R).
\]
At the corner $\Delta z=\Delta Y=0$, we also have
$\Delta z_{\Lambda\Lambda}=(2z/Y^2)\Delta Y_\Lambda>0$.
These first-contact inequalities show that the difference path cannot
leave the closed positive quadrant after entering it.

At the terminal value $\Lambda_0$, the right outer point precedes the
section, so its $S$ is smaller than $S_0$. At fixed $\Lambda$, formula
\eqref{sk:response-relation} makes $S$ strictly decrease with $z$.
Therefore $\lim_{\Lambda\to\Lambda_0^-}\Delta z<0$.
At a turning point the equations above yield the finite Taylor expansions
\begin{align}\label{sk:turn-expansion}
 \chi&=\chi_t+\Psi_t\Lambda^2/8+o(\Lambda^2),\notag\\
 Y&=\Lambda-\chi_t\Lambda^3/12
 +(12\chi_t^2+2\chi_t-9\Psi_t)\Lambda^5/1440+o(\Lambda^5),\\
 z&=\Lambda/6-(2\chi_t+1)\Lambda^3/360
 +(6\chi_t^2+8\chi_t-9\Psi_t+2)\Lambda^5/30240+o(\Lambda^5).\notag
\end{align}
They follow by inserting $x-x_t=-V^2/(8g_t)+O(V^4)$ and
$z=L_tV/6+O(V^3)$ into \eqref{sk:Lambda}; alternatively multiply
the singular equation for $z$ by $\Lambda^2$ and integrate from zero.
The assumed regularity suffices for these expansions. If
$\chi_{L,t}<\chi_{R,t}$, both
differences start positive. If they are equal, then
$\Psi_{L,t}<\Psi_{R,t}$ by the strict curvature signs, and the
fifth-order terms again make both differences positive. Either case
contradicts the terminal inequality. Hence $\chi_{L,t}>\chi_{R,t}$
and both differences start negative. A first zero of $\Delta z$ would
have $\Delta Y\ge0$ and enter the forbidden quadrant; at its corner
the strictly positive second derivative also contradicts approach
from negative values. We have proved
\begin{equation}\label{sk:Lambda-order}
 z_L(\Lambda)<z_{R,\mathrm{out}}(\Lambda),\qquad
 S_L(\Lambda)>S_{R,\mathrm{out}}(\Lambda),\quad0<\Lambda<\Lambda_0.
\end{equation}
The first-contact arguments apply on compact subintervals of
$(0,\Lambda_0)$; the ordering up to the section follows by taking
one-sided limits, with $Y_L\to\infty$ there.

Let the outer and inner left branches refer to the increasing and
decreasing parts of $S_L$, in that order along $V$. For $0<S\le S_0$,
the left outer point has $\Lambda_L(S)<\Lambda_0$. By
\eqref{sk:Lambda-order}, the right arc attains the same $S$ later.
If this occurs before its maximum of $\Lambda$, its $\Lambda$ is
larger; if afterwards, its $\Lambda\ge\Lambda_0$ is again larger.
At fixed $S$,
\begin{equation}\label{sk:horizontal}
 z_\Lambda=2/(\Lambda-S)^2-\tfrac12\operatorname{csch}^2(\Lambda/2)>0,
\end{equation}
since $0<\Lambda-S<\Lambda<2\sinh(\Lambda/2)$. Thus
\begin{equation}\label{sk:S-order}
 \Lambda_{L,\mathrm{out}}(S)<\Lambda_R(S),\quad
 z_{L,\mathrm{out}}(S)<z_R(S),\qquad0<S\le S_0.
\end{equation}

To compare the weighted areas, including both parts of the left
branch, define
\[
 k=P^{-2}+N^{-2},\qquad \mathcal Q=kH/(K_0dz)>0.
\]
Equations \eqref{sk:response} give
\[
 AC=K_0(V\mathcal R/t)_V+K_0k\mathcal R+\alpha VL_V.
\]
Integrate first on a truncated interval and combine the boundary terms:
\[
 K_0V\mathcal R/t+\alpha VL=D_E/t-VLC/t.
\]
This tends to zero at the turning point and to $(q-e)/(qm)$ at the
section. The latter is $\alpha d_0$, since $H_0=\alpha qm$.
Since $\int L\,dV=\Lambda_0$, each branch separately gives
$\int AC\,dV=K_0\int k\mathcal R\,dV-\alpha S_0$.
The combined boundary terms thus have finite limits. By
\eqref{sk:Iprime} and \eqref{sk:response},
\begin{equation}\label{sk:area}
 \frac{I'}{K_0}=
 \int_0^{S_0}(\mathcal Q_R-\mathcal Q_{L,\mathrm{out}})\,dS
 -\int_{S_0}^{S_{\max}}(\mathcal Q_{L,\mathrm{out}}
                              -\mathcal Q_{L,\mathrm{in}})\,dS.
\end{equation}

We establish both signs by Green's formula. In coordinates $(\lambda,s)$
put $d=\lambda-s$, $z=\coth(\lambda/2)-2/d$, and
$d_3(z)=2z(3+z^2)/(1-z^4)$. On $s>0$, $d>0$, $0<z<1$ define
\[
 \omega=d\log\frac{1+z^2}{(1-z^2)^2}
       -\frac{1-z^2}{2z}\,d\lambda-d\log d+\frac{ds}{d}
       =a\,d\lambda+b\,ds.
\]
The identities $d\log(V^2z)=(1-z^2)d\Lambda/(2z)$ and
$d\log H=dS/(\Lambda-S)$ show that $d\log\mathcal Q=\omega$
along each orbit path. To compare the values of
$\mathcal Q$, we integrate $\omega$ between the paths. Differentiation gives
\[
 a=d_3z_\lambda-(1-z^2)/(2z)-1/d,\quad
 b=2(d-d_3)/d^2,\quad b_\lambda-a_s=1/(z^2d^2)>0.
\]
Writing $\tau=\tanh(\lambda/2)$, one also obtains
\[
 a=\frac{\tau(9z^4+8z^2-1)-z(1+12z^2+3z^4)}
 {2\tau z(1-z^4)}<0.
\]
For $z\le1/3$, the coefficient $9z^4+8z^2-1$ is nonpositive,
so the numerator is negative. For $z>1/3$, using $0<\tau<1$
gives the strict upper bound
$-(1-z)^2(3z^3-3z^2+3z+1)<0$; the last cubic equals
$1+3z(z^2-z+1)>0$.

For $S_0\le S<S_{\max}$, follow the outer left path to its maximum,
return on the inner left path to level $S$, and close horizontally.
This is a clockwise simple curve; its interior has horizontal slices
$\Lambda_{L,\mathrm{out}}(s)<\lambda<\Lambda_{L,\mathrm{in}}(s)$,
$S<s<S_{\max}$. The complete region lies in the domain of $\omega$
by \eqref{sk:horizontal}. Green's formula gives
\[
 \log\frac{\mathcal Q_{L,\mathrm{in}}(S)}{\mathcal Q_{L,\mathrm{out}}(S)}
 =\int_{\Lambda_{L,\mathrm{out}}(S)}^{\Lambda_{L,\mathrm{in}}(S)}a(\lambda,S)\,d\lambda
 -\iint\frac{d\lambda\,ds}{z^2d^2}<0.
\]
For $0<S\le S_0$, follow the whole remaining left path to the
section, then the right path backwards to its level $S$, and close
horizontally. This is again clockwise and simple, by
\eqref{sk:S-order}; above $S_0$ its slices are the entire left tail,
and below $S_0$ they run from $\Lambda_{L,\mathrm{out}}(s)$ to
$\Lambda_R(s)$. The two densities have the same value at the
section because $P,N,E_+,E_-$ coincide there. Thus
\[
 \log\frac{\mathcal Q_R(S)}{\mathcal Q_{L,\mathrm{out}}(S)}
 =\int_{\Lambda_{L,\mathrm{out}}(S)}^{\Lambda_R(S)}a(\lambda,S)\,d\lambda
 -\iint\frac{d\lambda\,ds}{z^2d^2}<0.
\]
For fixed $S>0$ these regions stay away from $d=0$, $z=0,1$.
The boundaries are piecewise $C^1$ in the original $x$ coordinate:
\[
 \Lambda_x=-fV/t,\qquad S_x=-dz(gd-fV)/t.
\]
Their section limits agree and are finite; at the left maximum of
$S$, $\Lambda_x\ne0$. Thus Green's formula applies, including at the
maximum where an inverse branch may have a square-root derivative.
Finally, the turning-point expansions are
\[
 \mathcal R=\frac{A_tH_t}{12K_0}V^2+o(V^2),\quad
 S=\frac{A_tL_t^2}{360}V^5+o(V^5),\quad
 \mathcal Q=\frac{48H_t}{K_0L_t^2V^4}(1+o(1)).
\]
Here $A_t,L_t,H_t>0$ on both sides. Thus
$\mathcal Q=O(S^{-4/5})$ is integrable. At the section,
$k\mathcal R=O((V_0-V)^{-1/2})$ is integrable; at the left maximum
$\mathcal Q$ is finite and continuous. We may therefore pass to the
endpoints in \eqref{sk:area}. Both contributions are strictly
negative, so $I'<0$ and $\Pi''>0$ as asserted.
\end{proof}

\section{Connected parameter families and the two-cycle bound}\label{sk:family-proof}
To pass from the direction of neutral cycles to a global count, we
parametrize periodic orbits by their section half-width. We first establish
that the half-arcs in Theorem~\ref{sk:family} exist at every positive
half-width. Write
$\phi_i(h,\theta)=F(x_i(h),\theta)$. Under (S1), both inverse branches
exist for all $h>0$, and
\begin{equation}\label{sk:profile-origin}
 x_{R,L}(h)=\pm\sqrt{2h/g_1}+O(h),\quad
 \phi_{R,L}(h)=\pm f(0,\theta)\sqrt{2h/g_1}+O(h),\quad
 \phi_{i,h}=O(h^{-1/2}).
\end{equation}

\begin{lemma}\label{sk:half-arcs}
Under (S1), the half-arcs and midpoint functions in
Theorem~\ref{sk:family} exist at every positive half-width, are $C^2$ in
half-width and parameter, and satisfy $|M_i|<r$, $|M_{i,r}|<1$.
Every periodic orbit corresponds to exactly one zero of $\Delta$.
\end{lemma}
\begin{proof}
Fix a side and parameter. Solve
$h_y=\phi(h)-y$ from $h(\phi(t))=t>0$. Every positive-height critical
point has $h_{yy}=-1$. Thus this is the unique maximum; both branches
decrease away from it and remain in $0<h\le t$. Since $\phi$ is
bounded on $[0,t]$, neither branch can continue to infinite $|y|$
while remaining positive. Their finite endpoints satisfy $A\le0\le B$.
An endpoint $(h,y)=(0,0)$ would give an original orbit confined to
one half-plane and asymptotic to the equilibrium. This is impossible:
with $f_0=f(0,\theta)$ and $\omega=\sqrt{g_1-f_0^2/4}>0$, the local
coordinates $X=x$, $Z=(y-f_0x/2)/\omega$ have angular speed
$-\omega+O(\sqrt{X^2+Z^2})$. An orbit approaching the origin must
cross $X=0$ repeatedly. Hence $A<0<B$.

The endpoints are reached in finite original time, since
$dt_{\rm time}/dy=-1/g(x_i(h))=O(h^{-1/2})$ and
$h\asymp y-A$ or $B-y$ there. They are transverse for the $C^7$
original flow, so are $C^2$ functions of $t,\theta$. At fixed $y$ the
height variation $H=h_t$ satisfies $H_y=\phi_h(h)H$ and
$H(\phi(t))=1$. The coefficient is integrable at both endpoints,
so $H$ has finite positive limits. Differentiating the endpoint
conditions gives $H(A)=AA_t$, $H(B)=BB_t$, hence $A_t<0<B_t$.
Local rotation and a linear radial bound give $A,B\to0$ as $t\to0$,
uniformly on the compact parameter interval.

On a height slice write $y_\pm=\phi(h)+\widehat M\pm\widehat r$.
Since $dy/dh=1/(\phi-y)$,
\[
 (\widehat r^2)_h=-2\widehat r^2/(\widehat r^2-\widehat M^2)\le-2.
\]
Integrating to the maximum proves $r(t)^2\ge2t$. Thus $r(t)$
increases bijectively from zero to infinity. Its $C^2$ inverse gives
$M(r,\theta)$, with $|M_r|<1$ because both $B$ and $-A$ strictly
increase with $r$.

To identify all returns, let $m_i=-A_i$, $q_i=B_i$.
Since $q_i(t)>\phi_i(t)\to\infty$, the geometric half-return
$n_i:m_i\mapsto q_i$ is an increasing $C^2$ bijection
$(0,m_i^\infty)\to(0,\infty)$, where $m_i^\infty$ may be finite.
The actual first return is
\[
 \Pi=n_L\circ n_R^{-1},\qquad
 \operatorname{dom}\Pi=\{q>0:n_R^{-1}(q)<m_L^\infty\}.
\]
Its domain is an initial interval and $\Pi$ is strictly increasing.
Every periodic orbit crosses the section transversely. Its positive
section points form a finite periodic orbit of $\Pi$; an increasing
interval map has no periodic points of period greater than one.
Thus a periodic orbit consists of exactly one arc on each side.
The two arcs match precisely when their half-widths and midpoints agree,
which is exactly $\Delta=0$.
\end{proof}

Two integral formulas give the signs of the midpoint difference and
its parameter derivative. At fixed half-width,
center each arc by $y=M_i+z$, $-r\le z\le r$. Put
$w=h_R-h_L$, $\delta(h)=\phi_R(h)-\phi_L(h)$, and
\[
 a(z)=\int_0^1\phi_{R,h}((1-s)h_L+s h_R)\,ds,\qquad
 K(z)=\exp\left(-\int_{-r}^za(v)\,dv\right).
\]
Then $w_z=aw+\delta(h_L)-\Delta$. The estimate
\[
 \int_0^1\frac{ds}{\sqrt{(1-s)v+su}}=\frac2{\sqrt u+\sqrt v}
\]
and \eqref{sk:profile-origin} show that $a=O((r-|z|)^{-1/2})$.
Therefore $K$ has finite positive endpoint values, and integration gives
\begin{equation}\label{sk:kernel}
 \Delta(r,\theta)=\frac{\int_{-r}^rK(z)\delta(h_L(z),\theta)\,dz}
 {\int_{-r}^rK(z)\,dz}.
\end{equation}
For a single arc with endpoints $A=M-r$, $B=M+r$, differentiation
at fixed $r$ and fixed original $y$ gives
$v_y=\phi_hv+\phi_\theta$, $v(A)=AM_\theta$, $v(B)=BM_\theta$.
Using $K_i(y)=\exp(-\int_A^y\phi_{i,h}(h_i(s))\,ds)$ yields
\begin{equation}\label{sk:parameter-kernel}
 M_{i,\theta}=
 \frac{\int_A^BK_i\phi_{i,\theta}(h_i)\,dy}{BK_i(B)-AK_i(A)}.
\end{equation}
All integrals converge by the same endpoint estimate, and the
denominator is positive. Assumption (S3) implies
$M_{R,\theta}>0>M_{L,\theta}$ and $\Delta_\theta>0$.
At the lower endpoint, (S3) and $f(0,\theta_-)<0$ make the integrand
in \eqref{sk:kernel} nonpositive and strictly negative near both arc
endpoints. Consequently $\Delta(r,\theta_-)<0$ for all $r>0$.

\begin{proof}[Proof of Theorem~\ref{sk:family}]
For each half-width $r$, strict parameter monotonicity and the negative
lower boundary imply that an interior zero exists precisely when
$\Delta(r,\theta_+)>0$, and then it is unique. The implicit function
theorem patches these zeros into one $C^2$ graph $\theta=\beta(r)$
over the stated interval. Lemma~\ref{sk:half-arcs} shows that it
exhausts all periodic orbits.

At fixed parameter, the half-return maps for nearby initial points
define $r_R,r_L$ by
\[
 q=r_R+M_R(r_R),\quad r_L-M_L(r_L)=r_R-M_R(r_R),\quad
 \Pi(q)=r_L+M_L(r_L).
\]
Thus
\[
 \Pi'(q)=\frac{(1+M_L'(r_L))(1-M_R'(r_R))}
 {(1-M_L'(r_L))(1+M_R'(r_R))}.
\]
At a fixed point $r_R=r_L=r$,
\[
 \Pi'-1=-\frac{2\Delta_r}{(1-M_L')(1+M_R')}.
\]
Therefore $\beta'=0$ exactly at a neutral orbit. If
$\nu=M_L'=M_R'$ there, differentiation of the neighboring-orbit
formula, with $(r_R)_q=(r_L)_q=(1+\nu)^{-1}$, gives
\begin{equation}\label{sk:graph-curvature}
 \Pi''=-\frac{2\Delta_{rr}}{(1-\nu)(1+\nu)^2},\qquad
 \beta''=\frac{(1-\nu)(1+\nu)^2}{2\Delta_\theta}\Pi''>0.
\end{equation}
Every stationary point of the graph is consequently a nondegenerate
minimum. More explicitly, $\beta'$ changes from negative to positive
at each of its zeros; two such zeros would force an intervening zero
with the opposite crossing direction. Thus there is at most one
stationary point. Rolle's theorem gives at most two intersections
with any horizontal line. The inner and outer intersections, if both
exist, have $\beta'<0$ and $\beta'>0$, respectively, giving
$0<\Pi'_{\rm inner}<1<\Pi'_{\rm outer}$. A neutral cycle has
$\Pi''\ne0$, hence multiplicity exactly two; this also excludes
a period annulus. Orbits of bounded half-width remain in a compact subset
of the strip:
$t_i\le r^2/2$, $|y|<2r$, and $x$ remains in the compact interval
$[x_L(r^2/2),x_R(r^2/2)]\Subset J$.
\end{proof}

\section{A same-primitive curvature criterion at a Hopf boundary}\label{sk:hopf-proof}
At the Hopf boundary, comparison of branches with the same value of
$F$ gives the return derivative directly. This proves
Theorem~\ref{sk:hopf} without the neutral-cycle comparison.

\begin{proof}[Proof of Theorem~\ref{sk:hopf}]
The half-arc constructions above apply because $f(0)=0$.
Let $L_0=f'(0)/g_1>0$, $\ell_1=L_x(0)>0$. The same-$F$ branches
are defined for every $\omega>0$, and
\[
 D(\omega):=L_R(\omega)-L_L(\omega)
 =2\ell_1\sqrt{2\omega/f'(0)}+O(\omega)>0
\]
near zero. Since $D''<0$ and $D(0)=0$, $D(\omega)/\omega$
strictly decreases. Thus either $D>0$ everywhere or it crosses zero
exactly once, from positive to negative. Put $\lambda_i=1/L_i$.

Suppose a periodic orbit has section endpoints $q,-m$. Its upper
and lower branches in $\omega=F(x)$ are
$\zeta_i(\omega)>\omega>\eta_i(\omega)$, up to turning heights
$\omega_i^*>0$, and satisfy
\begin{equation}\label{sk:hopf-arcs}
 y_\omega=\lambda_i(\omega)/(\omega-y),\qquad
 \zeta_i(0)=q,\quad\eta_i(0)=-m.
\end{equation}
Integration of $g\,dx+(y-F)\,dy=0$ around the orbit, followed by
integration by parts, yields
\begin{equation}\label{sk:hopf-area}
 \int_0^{\omega_R^*}(\zeta_R-\eta_R)\,d\omega
 =\int_0^{\omega_L^*}(\zeta_L-\eta_L)\,d\omega.
\end{equation}
The boundary terms vanish because $\omega=0$ at the section and the
branch separation vanishes at a turning point.

For either pair of graphs $d=y_R-y_L$ satisfies
\[
 d_\omega=a d+b,\quad
 a=\frac{\lambda_R}{(\omega-y_R)(\omega-y_L)},\quad
 b=\frac{\lambda_R-\lambda_L}{\omega-y_L}.
\]
On every common compact interval before a turning point the
integrating-factor formula is valid. Initially $\lambda_R<\lambda_L$,
so $\zeta_R>\zeta_L$ and $\eta_R<\eta_L$. If there is a crossing
height $\omega_c$ of the $L$ branches, these strict inequalities still
hold at $\omega_c$. Above it, each upper graph difference can cross
zero only from positive to negative and each lower difference only
from negative to positive. If no $\omega_c$ existed, the initial
ordering, and any extra interval up to the larger turning height,
would contradict \eqref{sk:hopf-area}; a smaller right turning
height would already contradict the upper-graph ordering.

These same one-crossing directions exclude
$\omega_R^*>\omega_L^*$ by the area identity: at the left turning
point the right graphs enclose the left point, so no earlier
crossing could have occurred. If the turning heights agree above
$\omega_c$, the expansions from the regular inverse equation
$\omega_y=(\omega-y)/\lambda_i$ are
\[
 \zeta_i=\omega^*+\sqrt{2\lambda_i(\omega^*)(\omega^*-\omega)}
 +O(\omega^*-\omega),\quad
 \eta_i=\omega^*-\sqrt{2\lambda_i(\omega^*)(\omega^*-\omega)}
 +O(\omega^*-\omega).
\]
Since $\lambda_R>\lambda_L$ there, the initial ordering holds again
near the common endpoint, precluding prior crossings and again
contradicting the area identity. Equality at or below $\omega_c$
is excluded by the initial ordering and the same identity. Hence
\begin{equation}\label{sk:hopf-turns}
 \omega_c<\omega_R^*<\omega_L^*.
\end{equation}

Let $\mu=\omega_L^*/\omega_R^*>1$ and rescale the left graphs by
$\widetilde y_L(\omega)=y_L(\mu\omega)/\mu$. Their coefficient is
$\widetilde\lambda_L(\omega)=\lambda_L(\mu\omega)/\mu$.
Since $(\omega L_L)_\omega>0$,
$(\lambda_L/\omega)_\omega<0$, and therefore
$\widetilde\lambda_L<\lambda_L$. The scaled upper graph lies below
the original upper graph: their difference is negative at both
endpoints, and any interior zero would have strictly positive
derivative. The scaled lower graph similarly lies above the original
one, because its endpoint differences are positive and its zero
derivatives are negative. Up to $\omega_c$, the right graphs enclose
the original left graphs and thus the scaled ones. Above $\omega_c$,
$\lambda_R>\lambda_L>\widetilde\lambda_L$, so the one-crossing
directions, together with the square-root expansions at their common
turning point, rule out every intersection. Consequently
\[
 \zeta_R>\widetilde\zeta_L>\omega>
 \widetilde\eta_L>\eta_R\qquad(0<\omega<\omega_R^*).
\]
Using $f\,dt=d\omega/(y-\omega)$ with the time orientation of each arc,
the scaling cancels from numerator and denominator and gives
\[
 \oint f\,dt=\int_0^{\omega_R^*}\left(
 \frac1{\zeta_R-\omega}-\frac1{\widetilde\zeta_L-\omega}
 +\frac1{\widetilde\eta_L-\omega}-\frac1{\eta_R-\omega}\right)d\omega<0.
\]
Every term is integrable, with at most an inverse-square-root
singularity at the turning point. Thus every periodic orbit is
hyperbolic and repelling.

To deduce uniqueness, we determine the midpoint difference near zero.
From $F_h=L$,
\[
 \phi_{R,L}(h)=L_0h\pm\tfrac23\ell_1\sqrt{2/g_1}\,h^{3/2}+O(h^2).
\]
For fixed small half-width $r$, $t_i\le r^2/2$ and the centered arc
equation gives $M_i=O(r^2)$,
$h_i(z)=(r^2-z^2)/2+O(r^3)$, and $K=1+O(r)$ uniformly.
Substitute into \eqref{sk:kernel}, using
$\int_{-1}^1(1-s^2)^{3/2}ds=3\pi/8$, to obtain
\begin{equation}\label{sk:small-width}
 \Delta(r)=\frac{\pi\ell_1}{8\sqrt{g_1}}r^3+O(r^4)>0.
\end{equation}
At any positive zero, the return derivative formula in the preceding
proof gives $\Delta_r<0$. A differentiable function initially
positive whose every zero has negative derivative has at most one
zero: after its first downward crossing a return to nonnegative
values would force a zero without negative derivative. This proves
uniqueness and excludes a period annulus.
\end{proof}

\section{Proof of the branch-curvature two-cycle criterion}\label{sk:coefficient-proof}
\begin{proof}[Proof of Theorem~\ref{sk:two}]
Conditions (S1) and (S3) give the geometric and lower-boundary
hypotheses of Theorem~\ref{sk:family}. At each interior parameter,
(S2) gives the hypotheses of Proposition~\ref{sk:neutral}.
At the upper endpoint $L_H>0$ on both sides and
\[
 (L_H)_{FF}=\Xi_H/L_H^3.
\]
Thus the left branch is strictly convex and the right branch strictly
concave as functions of $F_H$; their difference has negative second
derivative. Since $f_H<0$ on the left, (S4) is exactly
$(\omega L_L)_\omega>0$ there. All hypotheses of
Theorem~\ref{sk:hopf} hold. Its proof gives
$\Delta_H>0$ near zero and $\Delta_{H,r}<0$ at every zero, so
$\{r>0:\Delta_H(r)>0\}$ is one initial interval. Now apply
Theorem~\ref{sk:family}. Its graph and multiplier conclusions give
all interior assertions. Formula \eqref{sk:kernel} excludes cycles
at the lower endpoint, while Theorem~\ref{sk:hopf} gives the upper
endpoint assertions.
\end{proof}

\section{Applications}\label{sec:applications}
\subsection{Classical Li\'enard systems of degree four}\label{ap:quartic-section}
For classical Li\'enard systems, we take the degree to be that of $F$; the damping $f=F'$ has degree one less. Theorem~\ref{cr:unique} gives the following proof of the quartic result of Li and Llibre \cite{LiLlibre2012}.
\begin{corollary}\label{ap:quartic}
For arbitrary real coefficients, the system
\[
 \dot x=y-F(x),\qquad \dot y=-x,\qquad
 F(x)=a_0+a_1x+a_2x^2+a_3x^3+a_4x^4
\]
has at most one limit cycle in the plane. A limit cycle can exist only if $a_1a_3<0$, and is then hyperbolic, attracting if $a_3>0$ and repelling if $a_3<0$. If $a_1a_3\ge0$ and $(a_1,a_3)\ne(0,0)$, there are no nonconstant periodic orbits. If $a_1=a_3=0$, every periodic orbit is nonisolated.
\end{corollary}
\begin{proof}
The global translation $y\mapsto y-a_0$ removes the constant term. The coefficients are smooth, $G=x^2/2$, and $xg=x^2>0$ off the origin. Taking $K=8a_4$ on both sides gives
\[
 \mathcal U_{8a_4}(x)=xF''(x)-F'(x)-8a_4x^3=-a_1+3a_3x^2.
\]
If $a_1a_3<0$, this has the strict sign of $a_3$ throughout both half-lines. Theorem~\ref{cr:unique} gives the asserted bound, hyperbolicity and stability.

For the remaining coefficients write the potential profiles as
\[
 \phi_\pm(u)=E(u)\pm p(u),\quad E(u)=2a_2u+4a_4u^2,\quad
 p(u)=a_1\sqrt{2u}+a_3(2u)^{3/2}.
\]
If $a_1a_3\ge0$ and the two odd coefficients are not both zero, $p$ has one strict sign for $u>0$. A periodic orbit would give positive half-arcs $u_\pm(y)$ on a common interval $(A,B)$, both vanishing at its endpoints and satisfying $u_\pm'=\phi_\pm(u_\pm)-y$. Their difference $w$ would satisfy
\[
 w'=[2a_2+4a_4(u_++u_-)]w+p(u_+)+p(u_-),\qquad w(A)=w(B)=0.
\]
The coefficient is continuous on $[A,B]$ and the last two terms have the same strict sign inside. The integrating-factor formula from $A$ to $B$ therefore gives $w(B)\ne0$, a contradiction.

If $a_1=a_3=0$, the two profiles coincide. Their endpoint and midpoint functions consequently coincide on their common transverse domain. Every matching half-arc therefore belongs to an open family of matching half-arcs; any periodic orbit has neighboring periodic orbits and is not isolated. This also covers polynomials of lower degree. The condition $a_1a_3<0$ remains a necessary condition for existence.
\end{proof}

\subsection{The Selkov system}\label{sk:application}
Consider the polynomial Selkov system
\begin{equation}\label{sk:original}
 \dot X=-X+(\alpha+X^2)Y,\qquad
 \dot Y=b-(\alpha+X^2)Y,
\end{equation}
where $b>0$. We use the normalized parameter $a=\alpha/b^2$.
We first allow negative $a$, a mathematical extension beyond the
physical range of the reaction parameters. This range is considered in
Conjecture 3.1 of Cao, Chen, Llibre and Tang
\cite[p.~13]{CaoChenLlibreTang2025}. We then turn to positive parameters.

\subsubsection{Negative parameters: the global upper bound}
Put
\[
 p=-a,\quad p_*=(3+2\sqrt3)^{-1},\quad c=\sqrt p,\quad
 B=b^2,\quad B_H=(1+p)/(1-p)^2.
\]
\begin{corollary}\label{sk:negative}
If $1/9<p<p_*$ and $0<B<B_H$, system~\eqref{sk:original}, with
$\alpha=-pb^2$, has at most two limit cycles in the whole plane.
Each nonhyperbolic cycle has multiplicity two, and for fixed $p$
there is at most one parameter $b$ with a nonhyperbolic cycle. If
two cycles coexist, the inner one is hyperbolic and attracting and
the outer one is hyperbolic and repelling. On $B=B_H$ there is at
most one cycle, necessarily hyperbolic and repelling. The bound two
is attained for parameters with $p<p_*$ sufficiently close to $p_*$.
\end{corollary}
In the normalization of \cite{CaoChenLlibreTang2025}, its conjecture
asks for at most two cycles when
$-(3+2\sqrt3)^{-1}<a<-1/9$, $b>0$ and
$b^2(1+a)^2+a-1<0$. These are exactly the open parameter conditions
in Corollary~\ref{sk:negative}.

Every periodic orbit of \eqref{sk:original} lies in $X>bc$.
Indeed, integration of $\dot X+\dot Y=b-X$ over a period gives
$\langle X\rangle=b>bc$, whereas on $X=bc$ the derivative
$\dot X=-bc<0$. A periodic function cannot cross a level exclusively
in one direction, and tangency is excluded there. The mean value
therefore places the orbit in the right component. On that component
the changes
\begin{equation}\label{sk:transformation}
 u=X/b,\quad y=X+Y-b-\frac1{b(1-p)},\qquad
 d\tau=b(u^2-p)\,dt
\end{equation}
are smooth, with a positive time factor and spatial determinant
$1/b$. They yield
\begin{equation}\label{sk:negative-lienard}
 u'=y-F_b(u),\quad y'=-g(u),\qquad
 F_b(u)=b(u-1)+\frac1b\left(\frac{u}{u^2-p}-\frac1{1-p}\right),
 \quad g(u)=\frac{u-1}{u^2-p}.
\end{equation}
Here $u>c$, primes denote $\tau$ derivatives, and
$f_b(u)=b-(u^2+p)/(b(u^2-p)^2)$. Periodic orbits and their
multipliers are preserved. To apply Theorem~\ref{sk:two}, translate
$x=u-1$ and take $\theta=b$.

The potential is
\begin{equation}\label{sk:potential}
 G(u)=\frac12\log\frac{u^2-p}{1-p}
 -\frac1{2c}\log\left(\frac{u-c}{u+c}\frac{1+c}{1-c}\right).
\end{equation}
It decreases on $(c,1)$, increases on $(1,\infty)$, and tends to
$+\infty$ at both endpoints. Its inverse branches are denoted by
$u_L(h)$ and $u_R(h)$.

\begin{lemma}\label{sk:threshold}
For $1/9<p<p_*$ set $k=1/c$ and
\[
 T_k(t)=\frac{k^2\sinh^2t\,[\cosh(kt)\cosh t-1]}
 {2[\cosh(kt)-\cosh t]^2},\quad t>0.
\]
This function has exactly one critical point, a nondegenerate
minimum, whose value $m(c)$ satisfies $9/8<m(c)<B_H<2$.
There is an increasing bijection $t=t(h)$ of $(0,\infty)$ such that
\begin{equation}\label{sk:threshold-difference}
 F_b(u_R(h))-F_b(u_L(h))
 =\frac{u_R(h)-u_L(h)}b[B-T_{1/c}(t(h))].
\end{equation}
For $0<B\le m(c)$ there is no periodic orbit.
\end{lemma}
\begin{proof}
Set $u=c\coth\vartheta$, $\vartheta_0=\operatorname{artanh}c$,
and $W(\vartheta)=\vartheta/c-\log\sinh\vartheta$.
Then $G(c\coth\vartheta)=W(\vartheta)-W(\vartheta_0)$.
Equal-potential angles are $z-t/2,z+t/2$, where
\[
 z=\operatorname{artanh}\frac{\tanh(t/2)}{\tanh(kt/2)}>t/2.
\]
They straddle $\vartheta_0$. As the potential increases, the smaller
angle decreases to zero and the larger increases without bound,
so their difference $t$ increases from zero to infinity.
For distinct $u,v>c$,
\[
 F_b(u)-F_b(v)=\frac{u-v}{b}
 \left[B-\frac{uv+p}{(u^2-p)(v^2-p)}\right].
\]
Substitution of the displayed angles gives
\eqref{sk:threshold-difference}.

To determine the critical points of $T_k$, differentiate to obtain
\[
 T_k'(t)=\frac{k^2\sinh t}{2[\cosh(kt)-\cosh t]^3}N_k(t),
\]
\begin{align*}
8N_k(t)={}&(3-k)\cosh((2k+2)t)-(k+1)\cosh((k+3)t)+2\cosh(2kt)\\
 &+(7k-15)\cosh((k+1)t)+(k+3)\cosh((2k-2)t)+6\cosh(2t)\\
 &-(7k+15)\cosh((k-1)t)+(k-1)\cosh((3-k)t)+18.
\end{align*}
For $\sqrt{3+2\sqrt3}<k<3$, the positive frequencies, in decreasing
order, are
$2k+2,k+3,2k,k+1,2k-2,2,k-1,3-k$.
The corresponding signs, including the constant last, are
$+,-,+,+,+,+,-,+,+$. Expansion into exponentials gives eight sign
variations. The exponential Descartes rule bounds the number of real
zeros, with multiplicity, by this number. It follows by induction:
divide by the least exponential and differentiate; Rolle's theorem
gives the usual loss of one zero, while when the first two
coefficients have the same sign the matching signs sufficiently
far to the left give one additional derivative zero before the first
zero of the function. This proves the induction also with repeated
zeros. Now
\[
 N_k(0)=N_k''(0)=N_k^{(4)}(0)=0,\quad
 N_k^{(6)}(0)=-30(k^2-1)(k^4-6k^2-3)<0.
\]
The positive leading exponential makes $N_k$ positive for large
$t$. Its sixfold zero at zero and one zero on each side exhaust the
eight-zero bound, so the positive zero is unique and simple.
Moreover $T_k(0+)=B_H$ and
$T_k(t)\sim k^2e^{(3-k)t}/8$ at infinity.

For the strict lower bound, put $v=\cosh t$, $w=\cosh(kt)$.
Differentiation gives
\[
 \partial_k\log T_k=\frac2k-t\sinh(kt)
 \frac{wv+v^2-2}{(wv-1)(w-v)}<0,
\]
because
\[
 \frac{wv+v^2-2}{(wv-1)(w-v)}-\frac1{w-1}
 =\frac{(v-1)[(2w-1)v+w-2]}{(wv-1)(w-v)(w-1)}>0
\]
and $kt\coth(kt/2)>2$. The triple-angle identity gives
$T_3(t)=9/8+9/(32\cosh^2t)>9/8$. Since $k<3$ and the minimum
occurs at finite $t$, $m(c)>9/8$. Finally $p<p_*<4/25$ gives
$B_H<725/441<2$.

If $B\le m(c)$, the two profiles have nonpositive difference,
strictly negative near zero. Suppose a periodic orbit existed.
Its two actual half-arcs share endpoints $-m,q$ and satisfy
$h_i' =\phi_i(h_i)-y$, $h_i(-m)=h_i(q)=0$. Their difference obeys
$w'=aw+\delta(h_L)$, where the same secant coefficient used in
\eqref{sk:kernel} is integrable at both endpoints. Its positive
integrating factor gives
$0=\int_{-m}^qK\delta(h_L)\,dy<0$, a contradiction.
The argument applies for every $0<B\le m(c)$ because the endpoint
estimate is along the hypothetical compact orbit.
\end{proof}

We verify (S1)--(S4) on $b\in[\sqrt{m(c)},\sqrt{B_H}]$.
The functions are analytic there, $(u-1)g(u)>0$ and
$g'(1)=1/(1-p)>1$. Formula~\eqref{sk:potential} gives properness of
$G$, while
$F_b(u)\sim1/[2b(u-c)]$ at $c$ and $F_b(u)\sim bu$ at infinity.
Since $b>1$ and $B_H<2$,
$-1<f_b(1)=b-B_H/b\le0$, so $f_b(1)^2<4g'(1)$.
Also
\[
 (f_b)_u=\frac{2u(u^2+3p)}{b(u^2-p)^3}>0.
\]
Thus there is exactly one damping zero $u_f$, with $u_f>1$ below
the Hopf boundary and $u_f=1$ on it.

For the three curvature inequalities, set $\eta=1-6p-3p^2>0$ and
$x=u-1>0$. For $L=f_b/g$,
\[
 L_u=\frac{N(x)}{b x^2[(x+1)^2-p]^2},
\]
\begin{align*}
N(x)={}&Bx^6+4Bx^5+[(5-p)B+1]x^4+4x^3+c_2x^2+c_1x+c_0,\\
c_2={}&6+4p-(5-p)(1-p)B,\\
c_1={}&4[1+p-(1-p)^2B],\qquad
c_0=(1-p)[1+p-(1-p)^2B].
\end{align*}
For $B\le B_H$, $c_0,c_1\ge0$ and $c_2\ge\eta/(1-p)>0$.
Hence $A=L_u/g>0$ on the right, including the Hopf endpoint.
With $d_u=u^2-p$, direct differentiation gives
\begin{equation}\label{sk:Xi-certificate}
 b^2\Xi=\frac{2P(u,p,B)}{(u-1)^6},
\end{equation}
where
\begin{align*}
P={}&(1-p)(u^2-p)^3B^2\\
 &-[2p^3-3p^2u+p^2+12pu^4-26pu^3+12pu^2
       +2u^6-3u^5+3u^4]B\\
 &+p^2+3pu^2-3pu+2p+u^3.
\end{align*}
This is strictly convex in $B$. Put $B_z=(u^2+p)/(u^2-p)^2$.
Then
\begin{align*}
 P(B_z)&=-2u^2(u^2+3p)^2(u-1)^2/(u^2-p)^2<0,\\
 P(B_H)&=(u-1)^3R_p(u)/(1-p)^3,\\
 R_p(u)&=(3p-1)(p+1)u^3+6p(p+1)u^2
          +3p(3p-1)(p+1)u+2p(3p-1).
\end{align*}
The discriminant of this cubic is
\[
 -108p^2(p-1)^2(p+1)^2(3p-1)(3p+1)^2(3p^2+6p-1)<0.
\]
Its unique real zero is negative, since its leading and constant
coefficients are negative; hence $R_p(u)<0$ for $u>0$.
On the left $B_z>B_H$, $P(B_H)>0>P(B_z)$, so the smaller root of
the convex quadratic $P$ lies between $B_H$ and $B_z$. This proves
$P(B)>0$ for every $B\le B_H$. On the right where $f_b>0$,
$B_z<B\le B_H$ and both endpoint values of $P$ are negative;
convexity proves $P(B)<0$. These give the required signs of
$\Xi$.

Differentiation with respect to the parameter gives
\[
 \partial_bF_b(u)=(u-1)
 \left[1+\frac{u+p}{b^2(1-p)(u^2-p)}\right].
\]
Together with Lemma~\ref{sk:threshold} this proves (S3).
At $b=\sqrt{B_H}$ cancellation of the removable factor gives
\begin{align*}
 L_H(u)&=(u+1)\left[b+\frac{2p}{b(1-p)(u^2-p)}\right],\\
 L_0=L_H(1)&=\frac{2(1+3p)}{b(1-p)^2},\qquad
 (L_H)_u(1)=\frac{\eta}{b(1-p)^3}>0.
\end{align*}
On the left $L_H$ is strictly convex in $\omega=F_H$ because
$(L_H)_{\omega\omega}=\Xi_H/L_H^3>0$. The supporting-line
inequality at a positive $\omega$, compared with the limit at zero,
gives $\omega(L_H)_\omega\ge L_H-L_0$. Moreover
\[
 L_H(u)>b(1+c),\qquad b(1+c)-L_0/2=\frac{c}{b(1+c)^2}>0.
\]
It follows that $(\omega L_H)_\omega\ge2L_H-L_0>0$. Since
$f_H<0$ on the left, $(F_HL_H)_u<0$. This verifies (S4).

\begin{proof}[Proof of the upper-bound assertions in Corollary~\ref{sk:negative}]
Theorem~\ref{sk:two} therefore applies on the closed band
$m(c)\le B\le B_H$. Lemma~\ref{sk:threshold} excludes all smaller
positive $B$, independently of the focus condition. The barrier and
transformation \eqref{sk:transformation} account for every periodic
orbit of the original polynomial system and preserve multiplicity
and stability.
\end{proof}

At $p=1/9$ the finite minimum defining $m(c)$ is unavailable, but
$T_3(t)>9/8$ permits the lower endpoint $B_-=9/8$ instead; all
strict coefficient inequalities above still hold. The same
two-cycle bound therefore holds at $p=1/9$.
At $p=p_*$, the condition $(L_H)_u(1)>0$ fails. Here $N_k$ has a
zero of order at least eight at zero, exhausting the exponential
zero bound. Its leading coefficient is positive, so $T_k'(t)>0$
for $t>0$, and $T_k(t)>B_H$. The profile comparison above
therefore excludes periodic orbits for $0<B\le B_H$.

\subsubsection{Attainment of the bound}
We construct parameters near $p=p_*$ and $B=B_H$ for which two small
cycles coexist. Introduce near $u=1$ the coordinate
$z=\operatorname{sgn}(u-1)\sqrt{2G(u)}$ and the positive local time
factor $d\sigma/d\tau=g(u)/z$. In velocity coordinates $w=u'$, the
system becomes
\[
 \frac{dz}{d\sigma}=w,\qquad
 \frac{dw}{d\sigma}=-z-h(z)w,\qquad h(z)=zf_b(u(z))/g(u(z)).
\]
On $B=B_H$ write $h(z)=\sum_{j\ge1}h_jz^j$. With $A_0=1-p$,
finite Lagrange inversion gives
\begin{equation}\label{sk:focus-coefficients}
 h_2=\frac{1-6p-3p^2}{bA_0^{5/2}},\qquad
 h_4=-\frac{(3p+1)(15p^2+102p-13)}{36bA_0^{7/2}}.
\end{equation}
These coefficients follow from the coefficient formula
\[
 h_{n-1}=A_0^{n/2}[x^{n-1}]
 \left\{f_b(1+x)\left(\frac{2A_0G(1+x)}{x^2}\right)^{-n/2}\right\},
 \qquad n=3,5.
\]
At $p=p_*$, $h_2=0$ and
$h_4=-16p_*^2/[3b_*(1-p_*)^{7/2}]<0$.
For the positive-section return map, the first nonzero focus
coefficients are $v_3=-\pi h_2/4$ and, when $h_2=0$,
\begin{equation}\label{sk:v5}
 v_{5*}=-\pi h_4/8=\frac{2\pi p_*^2}{3b_*(1-p_*)^{7/2}}>0.
\end{equation}
To justify these formulas, retaining only the odd part of $h$
gives a reversible local center under $(z,w,\sigma)\mapsto
(-z,w,-\sigma)$. The first nonzero even term changes the radial
equation by $-h_{2j}R^{2j+1}\sin^{2j}\theta\cos^2\theta$
to lowest order. The center solution and its propagation factor
are $R=q+O(q^2)$ and $1+O(q)$, so integration gives
$-\pi h_2q^3/4$ and $-\pi h_4q^5/8$, respectively. Lower odd
terms cannot alter the first of these nonzero orders.

Let $\varepsilon=p_*-p>0$, $\delta=B-B_H(p)$ and
\[
 \alpha_*=\frac{\pi\sqrt{1-p_*}}{2b_*},\qquad
 s_* =\frac{3\pi(1+p_*)}{4b_*(1-p_*)^{5/2}}.
\]
The analytic return map has linear coefficient
$e^{-\pi f_b(1)/\omega}-1=-2\alpha_*\delta+O(\delta^2+|\varepsilon\delta|)$,
where $\omega^2=(1-p)^{-1}-f_b(1)^2/4$.
Write $\Pi(q)-q=\sum_{j\ge1}c_j(\varepsilon,\delta)q^j$.
On the Hopf curve $c_2=0$ and
$c_3=-2s_*\varepsilon+O(\varepsilon^2)$; at the degenerate
point $c_4=0$ by the center comparison just given. Analytic
dependence consequently gives
\[
\begin{aligned}
 c_2&=O(|\delta|),&
 c_3&=-2s_*\varepsilon+O(\varepsilon^2+|\delta|),\\
 c_4&=O(|\varepsilon|+|\delta|),&
 c_5&=v_{5*}+O(|\varepsilon|+|\delta|).
\end{aligned}
\]
The remaining Taylor terms are $O(q^6)$ and their $q$ derivative
is $O(q^5)$, uniformly for parameters in a fixed small neighborhood.
Thus, with
$q=\sqrt\varepsilon R$ and $\delta=-\kappa\varepsilon^2$,
\[
 \frac{\Pi(q)-q}{q\varepsilon^2}
 =2\alpha_*\kappa-2s_*R^2+v_{5*}R^4+O(\sqrt\varepsilon),
\]
uniformly with its first derivative for $R$ in compact positive
intervals. Choose $0<\kappa<s_*^2/(2\alpha_*v_{5*})$.
The limiting polynomial has two distinct positive roots in $R$.
The implicit function theorem preserves both as simple positive
zeros for sufficiently small $\varepsilon>0$. These give two
hyperbolic periodic orbits with $B<B_H(p)$ and $1/9<p<p_*$.
The global bound makes them the only two cycles, completing
the attainment assertion in Corollary~\ref{sk:negative}.

\subsubsection{Positive parameters: uniqueness and existence}
For positive parameters we apply the branch criterion of Chen and Chen
\cite[Proposition~2.2, p.~3542]{ChenChen2015}. This argument is
independent of the two-cycle criterion.

\begin{theorem}\label{sk:positive}
If $a>0$, $b>0$ and $b^2(1+a)^2+a-1<0$, then
\eqref{sk:original}, with $\alpha=ab^2$, has exactly one limit cycle
in the whole plane. This cycle is hyperbolic and attracting.
\end{theorem}

We recall the required form of their criterion. Suppose
$xg(x)>0$ for $x\ne0$ on a finite
strip containing zero, $F(0)=0$, and let $F_1,F_2$ be the right and
left potential branches on $(0,z_1)$ and $(0,z_2)$. Suppose they
are not identical near zero. If some $\alpha_0$ smaller than both
endpoints satisfies
\begin{enumerate}
\item $F_1\le0\le F_2$ for $0<z<\alpha_0$, and $F_1>0$ for
 $\alpha_0<z<z_1$;
\item $F_2'(z)\le0$ whenever $z>\alpha_0$ and $F_2(z)<0$;
\item $F_1F_1'$ is nondecreasing on $(\alpha_0,z_1)$;
\item for every pair $\alpha_0<z<z_1$, $\alpha_0<v<z_2$ with
 $z\le v$ and $F_1(z)=F_2(v)$, one has $F_1'(z)\ge F_2'(v)$,
\end{enumerate}
then there is at most one periodic orbit in the strip, and it is
hyperbolic and attracting if it exists.

Every periodic solution of \eqref{sk:original} with $\alpha>0$
lies in the first quadrant: at a minimum of $Y$ one has
$Y=b/(\alpha+X^2)>0$, and at a minimum of $X$ one then has
$X=(\alpha+X^2)Y>0$. On $u=X/b>0$, the analogue of
\eqref{sk:transformation}, with $1-p$ replaced by $1+a$ and
positive time factor $b(a+u^2)$, gives
\[
 u'=y-\mathcal F(u),\quad y'=-\gamma(u),\quad
 \mathcal F(u)=b(u-1)+\frac1b\left(\frac{u}{a+u^2}-\frac1{1+a}\right),
 \quad\gamma(u)=\frac{u-1}{a+u^2}.
\]
Put $\mathcal G(u)=\int_1^u\gamma(v)\,dv$,
$\mathcal K=\mathcal F'/\gamma$, $B=b^2$, and
$B_H^+=(1-a)/(1+a)^2$. The assumptions imply $0<a<1$ and
$0<B<B_H^+$. The factorization
\[
 b\mathcal F(u)=\frac{(u-1)[B(1+a)(a+u^2)+a-u]}{(1+a)(a+u^2)}
\]
gives three distinct zeros $0<\ell<1<\rho$, with signs
$-,+,-,+$ on the four successive subintervals of $(0,\infty)$.
Indeed the quadratic factor is positive at zero, negative at one,
and positive at infinity. The numerator of
$b\mathcal F'=B+(a-u^2)/(a+u^2)^2$ is quadratic in $u^2$ and
has exactly two positive zeros. Rolle's theorem places these in
$(\ell,1)$ and $(1,\rho)$, so $\mathcal F'>0$ on
$(0,\ell)\cup(\rho,\infty)$.

The product and slope conditions require two algebraic inequalities.
The first is
\begin{equation}\label{sk:left-product}
 \frac1\gamma(\mathcal F\mathcal K)'
 =\mathcal K^2+\frac{\mathcal F}{\gamma}\mathcal K'>0
 \qquad(0<u<\ell).
\end{equation}
To prove this inequality, write
$d=a+u^2$, $\overline F=b\mathcal F$, $Q=b\mathcal K$ and
$\kappa=(u-a)/[(1+a)d]$. Then
\begin{align*}
 \overline F&=(u-1)(B-\kappa),\qquad
 Q=\frac{Bd^2+a-u^2}{(u-1)d},\\
 Q'&=\frac{Bd^2(u^2-2u-a)+R}{(u-1)^2d^2},\qquad
 R=u^4-4au^2+4au-a^2.
\end{align*}
Put $Z=Q^2+\overline FQ'/\gamma$. The relevant domain is
$0<u<1$, $B_L:=\max(0,\kappa)<B<B_H^+$, and
$Z_{BB}=4ud/(u-1)<0$. At $B_L=\kappa$, $Z=Q^2\ge0$.
If $B_L=0$, put $c_0=a-u\ge0$; the numerator
$N=(u-1)^2d^2Z$ satisfies
\[
 (1+a)N(a,u,0)=(1-u)\{(1+6u)c_0^2
 +2u(1+3u-u^2)c_0+u^2(1-u^2)\}>0.
\]
The upper endpoint is also positive. The domain implies
$w=(u+2)a-u>0$, since
$B_H^+-\kappa=(1-u)w/[(1+a)^2d]$.
Direct expansion at $B_H^+$ yields
\[
 Z(a,u,B_H^+)=\frac{\sum_{j=0}^4A_j(u)w^j}
 {(u+2)^4(1+a)^4d^2},
\]
where
\begin{align*}
 A_4={}&3(4u^2+7u+1),\qquad
 A_3=2(3u^5+10u^4+35u^3+60u^2+32u+4),\\
 A_2={}&(u+1)(2u^7+9u^6+20u^5+4u^4+46u^3+135u^2+68u+4),\\
 A_1={}&2u(1-u)(u+1)^2(3u^4+12u^3+33u^2+32u+4),\\
 A_0={}&4u^2(1-u)^2(u+1)^4.
\end{align*}
All are positive for $0<u<1$. For
$\theta=(B-B_L)/(B_H^+-B_L)\in(0,1)$, the exact quadratic
interpolation is
\[
 Z(B)=(1-\theta)Z(B_L)+\theta Z(B_H^+)
 +\frac{2ud}{1-u}\theta(1-\theta)(B_H^+-B_L)^2>0.
\]
Dividing by $b^2$ proves \eqref{sk:left-product}.
Evaluation at $B=0$ is valid because $Z$ is a polynomial in $B$;
the system itself is considered only for $b>0$.

On the right, $\mathcal K'(u)>0$ for all $u>1$. To see this, put
$h=u-1>0$; its positive denominator is
$bh^2[a+(h+1)^2]^2$ and its numerator is
\begin{align*}
 &Bh^6+4Bh^5+[(a+5)B+1]h^4+4h^3+c_2h^2+c_1h+c_0,\\
 &c_2=6-4a-(a+5)(a+1)B,\quad
 c_1=4[1-a-(a+1)^2B],\quad c_0=(a+1)[1-a-(a+1)^2B].
\end{align*}
Here $c_0,c_1>0$ and
$c_2>(1+6a-3a^2)/(1+a)>0$.

Condition (iv) requires a comparison of slopes at possibly different potential
heights. Suppose
$0<p_0<1<q_0$ and $\mathcal F(p_0)=\mathcal F(q_0)$. Let
$s=p_0+q_0$, $t=p_0q_0$, $D_p=a+p_0^2$, $D_q=a+q_0^2$.
Subtracting the two values gives
\[
 B=(t-a)/(D_pD_q),\quad 0<a<t,\quad
 b\mathcal F(p_0)=-\frac{(p_0-1)(q_0-1)[t-a(s+1)]}{(1+a)D_pD_q}.
\]
Thus a negative common value means $t/(s+1)<a<t$, and a positive
common value means $0<a<t/(s+1)$. Substitution into the slope
difference gives the affine expression
\begin{equation}\label{sk:cross-certificate}
 b[\mathcal K(q_0)-\mathcal K(p_0)]
 =\frac{(q_0-p_0)H(a)}{D_pD_q(p_0-1)(q_0-1)},\quad
 H(a)=2t^2-st+a(-s^2-2t+3s).
\end{equation}
The denominator is negative. If the common value is negative and
$\mathcal G(p_0)\le\mathcal G(q_0)$, then $s>2$: the strict
comparison
$\mathcal G(1-h)>\mathcal G(1+h)$ for $0<h<1$ follows by
integrating $v/[a+(1\mp v)^2]$, and otherwise would contradict the
assumed energy order. At $a_*=t/(s+1)$ and at $t$,
\[
 H(a_*)=2ts(p_0-1)(q_0-1)/(s+1)<0,\qquad
 H(t)=-ts(s-2)<0.
\]
Affine interpolation proves $H(a)<0$ in the entire relevant range.

If the common value is positive and
$\mathcal G(q_0)\le\mathcal G(p_0)$, then $2t<s$.
Indeed, let $E_a=\mathcal G(q_0)-\mathcal G(p_0)$ and
$E_0=\log(q_0/p_0)+1/q_0-1/p_0$. Then
\[
 E_a-\frac{E_0}{1+a}
 =\frac{a}{1+a}\int_{p_0}^{q_0}
 \frac{(u-1)^2(u+1)}{u^2(a+u^2)}\,du>0.
\]
The inequality $\log w>2(w-1)/(w+1)$ for $w>1$ follows by
differentiation; it gives $E_0>(q_0-p_0)(2t-s)/(st)$.
Thus $2t\ge s$ would contradict $E_a\le0$.
Now $H(0)=t(2t-s)<0$ and $H(a_*)<0$, so again $H(a)<0$.
Both cases prove the strict slope comparisons
\begin{equation}\label{sk:cross-results}
\begin{aligned}
 \mathcal F(p_0)=\mathcal F(q_0)<0,\quad
 \mathcal G(p_0)\le\mathcal G(q_0)&\ \Longrightarrow\ \mathcal K(q_0)>\mathcal K(p_0),\\
 \mathcal F(p_0)=\mathcal F(q_0)>0,\quad
 \mathcal G(q_0)\le\mathcal G(p_0)&\ \Longrightarrow\ \mathcal K(q_0)>\mathcal K(p_0).
\end{aligned}
\end{equation}

\begin{proof}[Proof of Theorem~\ref{sk:positive}]
Choose $0<\epsilon<\ell<1<\rho<U$ so that any proposed two
periodic orbits are contained in $\epsilon<u<U$. The potential
inverse branches $p(z)<1<q(z)$ are smooth on their open domains,
with $p'=1/\gamma(p)<0$ and $q'=1/\gamma(q)>0$.
If $\mathcal G(\ell)<\mathcal G(\rho)$, reflect space by
$(\widehat x,\widehat y)=(1-u,-y)$ without reversing time.
The right and left profiles become $\Psi_1=-\mathcal F(p(z))$,
$\Psi_2=-\mathcal F(q(z))$, with
$\alpha_0=\mathcal G(\ell)$. They satisfy condition (i) by the
four-interval sign table. Wherever $\Psi_2<0$ beyond $\alpha_0$,
$q>\rho$ and $\Psi_2'=-\mathcal K(q)<0$, proving (ii).
Formula~\eqref{sk:left-product} gives
$(\Psi_1\Psi_1')'>0$, proving (iii). For every pair $z\le v$
with $\Psi_1(z)=\Psi_2(v)>0$, we have a negative common value of
$\mathcal F$ and $\mathcal G(p(z))\le\mathcal G(q(v))$.
The first implication in \eqref{sk:cross-results} gives
$\Psi_1'(z)>\Psi_2'(v)$, proving (iv).

If $\mathcal G(\ell)\ge\mathcal G(\rho)$, retain the original
coordinates, take $\alpha_0=\mathcal G(\rho)$, and set
$P_1=\mathcal F(q(z))$, $P_2=\mathcal F(p(z))$.
Condition (i) again follows from the sign table; if $P_2<0$,
$p<\ell$ and $P_2'=\mathcal K(p)<0$, proving (ii).
For $z>\alpha_0$, $q>\rho$, and the already proved
$\mathcal K'(q)>0$ gives
$(P_1P_1')'=\mathcal K(q)^2+\mathcal F(q)\mathcal K'(q)/\gamma(q)>0$.
For the pairs in (iv), the common value is positive and
$\mathcal G(q(z))\le\mathcal G(p(v))$; the second implication in
\eqref{sk:cross-results} applies. In the equality case
$\mathcal G(\ell)=\mathcal G(\rho)$, such pairs do not exist:
$v>\alpha_0$ forces $P_2(v)<0<P_1(z)$. In both configurations
the branch endpoints exceed $\alpha_0$, and the profiles have
opposite strict signs near zero. Thus all hypotheses of the cited
criterion hold on the finite strip. This excludes two cycles and
shows that any cycle is hyperbolic and attracting. Since every
periodic orbit lies in $u>0$, uniqueness holds in the whole plane.

For existence, the first quadrant of \eqref{sk:original} is
positively invariant. Define
\[
 V=X+Y-\tfrac12\log(1+X),\qquad
 m_0=\inf_{X\ge0}\frac{\alpha+X^2}{1+X}
     =2(\sqrt{1+\alpha}-1)>0.
\]
Here $X/2+Y\le V\le X+Y$, and
\[
 \dot V=b-X+\frac{X}{2(1+X)}
 -\frac{\alpha+X^2}{2(1+X)}Y
 \le b+\tfrac12-X-\tfrac12m_0Y
 \le b+\tfrac12-c_0V,
\]
with $c_0=\min(1,m_0/2)>0$. All forward orbits are bounded.
Choose $T,M>0$ such that $X(t)\le M$ for $t\ge T$, and put
$c_1=\alpha+M^2$. Since $\dot Y\ge b-c_1Y$, scalar comparison gives
\[
 Y(t)\ge Y(T)e^{-c_1(t-T)}+\frac b{c_1}\bigl(1-e^{-c_1(t-T)}\bigr)
 \qquad(t\ge T).
\]
Thus $Y(t)\ge b/(2c_1)$ for $t\ge T_1:=T+(\log2)/c_1$.
The inequality $\dot X\ge-X+\alpha Y$ then yields
\[
 X(t)\ge X(T_1)e^{-(t-T_1)}+
 \frac{\alpha b}{2c_1}\bigl(1-e^{-(t-T_1)}\bigr)
 \qquad(t\ge T_1).
\]
Consequently $X(t)\ge\alpha b/(4c_1)>0$ for $t\ge T_1+\log2$.
Every positive orbit eventually lies in a compact subset of the
open quadrant. The unique equilibrium has positive determinant
$1/(1+a)$ and positive trace $-\mathcal F'(1)$ in the transformed
system, and hence is a hyperbolic source. A sufficiently small
Lyapunov ellipse about it has strictly outward vector field.
An orbit starting on that ellipse cannot reenter it; its compact
omega-limit set contains no equilibrium. The Poincare--Bendixson
theorem produces a periodic orbit. Together with
uniqueness and hyperbolicity, this yields exactly one attracting cycle.
\end{proof}

The strip inequalities (3.17)--(3.18) of
\cite{CaoChenLlibreTang2025} fail for some parameters in this region:
for $a=1/2$, $b=1/100$, the height difference on their left side is
negative although their right side is positive. Inequalities
\eqref{sk:left-product} and \eqref{sk:cross-results} provide the
product and slope comparisons used here throughout the stated
region. The cases $a=0$, $b=0$, and
$b^2(1+a)^2+a-1\ge0$ are outside the scope of this result.

\subsection{A uniqueness conjecture in the Journal of Mathematical Analysis and Applications}\label{ap:jmaa-section}
Chen, Li, Xiao and Yang \cite{ChenLiXiaoYang2025}, in Conjecture~2.1, ask whether their system~(1.8) has at most one limit cycle when
\[
 m,n\in\mathbb N_{>0},\quad n<m,\quad a>0,\quad
 \widetilde\varphi(a,m,n)<b<a^{(2n-1)/(2m)}.
\]
Here $\widetilde\varphi$ denotes the threshold on the homoclinic branch. By their Theorem~2.4,
$\widetilde\varphi(a,m,n)>a^{(2n-1)/(2m)}/(2n)>0$.
We prove the conjecture and extend the upper bound to a larger parameter set, allowing an arbitrary positive damping factor $\delta$. Their global two-cycle Conjecture~7.1 is discussed in Section~\ref{ap:other-section}.

\begin{lemma}\label{ap:mechanical-graphs}
A nonconstant periodic orbit of a $C^1$ system $\dot x=v$, $\dot v=A(x,v)$ consists of two graphs $v=u(x)>0$ and $v=-\ell(x)<0$ between its extreme abscissae $x_L<x_R$. Both functions are continuous on $[x_L,x_R]$, $C^1$ in its interior, and vanish precisely at its endpoints.
\end{lemma}
\begin{proof}
A least-period orbit is a regular Jordan curve. Fix $c$ and consider the ray $\{(c,v):v>0\}$. Every intersection is transverse, and its time-oriented tangent has positive horizontal component $\dot x=v$. The interior of an oriented Jordan curve lies locally on the same side of its tangent at every point. Thus traversal up this ray would cross every intersection in the same direction, either always into or always out of the interior. Successive transverse crossings must instead alternate between these two directions. There is therefore at most one intersection on this ray. The ray $\{(c,v):v<0\}$ has the analogous property, since all its horizontal tangent components are negative. At any regular value $c$ strictly between the minimum and maximum abscissae, periodicity supplies both an increasing and a decreasing crossing of $x=c$, so there are exactly two intersections, one on each ray.

At a point where $v=0$, $\dot v\ne0$, since the orbit cannot contain an equilibrium. Every turning point of $x$ is therefore nondegenerate. Compactness gives finitely many such points, and two distinct turning points cannot have the same abscissa, since they would be the same phase point. An interior turning point would change the number of intersections of neighboring regular vertical lines by two, contradicting the count just established. There are only the two extreme turning points. The implicit-function theorem gives the asserted graphs.
\end{proof}

\begin{lemma}\label{ap:left-localization}
For $m,n\in\mathbb N_{>0}$ and $\xi,\beta>0$, every nonconstant periodic orbit of
\begin{equation}\label{ap:normalized-jmaa}
 \dot x=v,\qquad \dot v=-g(x)-f(x)v,\qquad
 g(x)=x^{2m+1}-x,\quad f(x)=\xi(\beta+x^{2n-1})
\end{equation}
is contained in $x<0$.
\end{lemma}
\begin{proof}
At the extreme points of a periodic orbit, $g(x_L)<0$ and $g(x_R)>0$. If it crosses $x=0$, the sign of $g$ forces $x_L=-L<-1$ and $x_R=R>1$. The graphs in Lemma~\ref{ap:mechanical-graphs} satisfy
\[
 u'=-g/u-\xi(\beta+x^{2n-1}),\qquad
 \ell'=-g/\ell+\xi(\beta+x^{2n-1}).
\]
For $0\le x<T:=\min(L,R)$ put $D_u=u(x)-u(-x)$ and $D_\ell=\ell(x)-\ell(-x)$. Oddness gives
\[
 D_u'=\frac{g(x)}{u(x)u(-x)}D_u-2\xi\beta,\qquad
 D_\ell'=\frac{g(x)}{\ell(x)\ell(-x)}D_\ell+2\xi\beta.
\]
Both initial values are zero. On every compact subinterval of $[0,T)$ the coefficients are continuous, so positive integrating factors yield $D_u<0<D_\ell$ for $x>0$. If $R>L$, the limit of $D_u$ at $L$ is $u(L)>0$, a contradiction. If $R<L$, the limit of $D_\ell$ at $R$ is $-\ell(-R)<0$, also a contradiction. If $R=L$, then $g(x)>0$ on $(1,R)$ and $D_u'<-2\xi\beta$. Integrating from any $c\in(1,R)$ prevents $D_u$ from tending to zero at $R$, whereas both endpoint speeds vanish there.

A periodic orbit contained in $x>0$ is excluded by the strictly negative divergence $-\xi(\beta+x^{2n-1})$ in this simply connected half-plane. An orbit touching $x=0$ without crossing it would have $v=0$ there and pass through the equilibrium $(0,0)$, which is impossible. Every periodic orbit is therefore contained in $x<0$.
\end{proof}

\begin{theorem}\label{ap:jmaa}
If $m,n\in\mathbb N_{>0}$, $n\le m+1$, $\xi>0$ and $0<\beta<1$, system~\eqref{ap:normalized-jmaa} has at most one nonconstant periodic orbit in the whole plane. If it exists, it is a hyperbolic attracting limit cycle.
\end{theorem}
\begin{proof}
The global diffeomorphism $y=v+F(x)$, where $F(x)=\xi(\beta x+x^{2n}/(2n))$, gives the Li\'enard form. Lemma~\ref{ap:left-localization} reduces the count to $x<0$. Set
\[
 z=x+1,\quad w=y-F(-1),\quad I=(-\infty,1),\quad
 \mathcal F(z)=F(z-1)-F(-1),\quad \mathfrak g(z)=g(z-1).
\]
The coefficients are smooth, vanish at zero, and $z\mathfrak g(z)>0$ on $I\setminus\{0\}$. In particular, $\mathfrak g'(0)=2m>0$ and $\mathcal F'(0)=\xi(\beta-1)<0$. Its potential is
\[
 \mathcal G(z)=-\frac{(z-1)^2}{2}+\frac{(z-1)^{2m+2}}{2m+2}
                 +\frac{m}{2(m+1)}.
\]
The left endpoint energy is infinite and the right endpoint energy is $m/[2(m+1)]$. These unequal endpoint energies are allowed in Theorem~\ref{cr:unique}, which requires no invariance of the strip.

Choose $K=0$ in Theorem~\ref{cr:unique}. Put $t=-x>0$, so that the shifted coordinate is $z=1-t$, and let $r=2m+1$ and $s=2n-1$. Then $1\le s\le r$, and differentiation gives
\begin{equation}\label{ap:jmaa-certificate}
 \frac{\mathcal U_0(1-t)}\xi
 =(r-s)t^{r+s-1}+(s-1)t^s-r\beta t^{r-1}+\beta.
\end{equation}
Write $B(t)=(r-s)t^{r+s-1}+(s-1)t^s>0$. The nonnegative weights $r-s,s-1,1$ sum to $r$, and
$(r-s)(r+s-1)+(s-1)s=r(r-1)$. Weighted arithmetic--geometric mean, with zero-weight terms omitted, yields
\[
 B(t)+1\ge r t^{r-1},\qquad
 \mathcal U_0(1-t)/\xi
 =\beta[B(t)+1-rt^{r-1}]+(1-\beta)B(t)>0.
\]
This includes $s=1$ and $s=r$. The sets $t>1$ and $0<t<1$ are precisely the two sides of the shifted equilibrium. The strict inequality therefore holds on both sides of the entire left well. Theorem~\ref{cr:unique} and Lemma~\ref{ap:left-localization} give the bound in the whole plane.
\end{proof}

\begin{corollary}\label{ap:conjecture21}
The system
\begin{equation}\label{ap:original-jmaa}
 \frac{dX}{dt}=V,\qquad
 \frac{dV}{dt}=aX-X^{2m+1}-\delta(b+X^{2n-1})V
\end{equation}
has at most one nonconstant periodic orbit in the plane if
\[
 m,n\in\mathbb N_{>0},\quad n\le m+1,\quad a,\delta>0,\quad
 0<b<a^{(2n-1)/(2m)}.
\]
Any such orbit is hyperbolic and attracting. In particular, Conjecture~2.1 of \cite{ChenLiXiaoYang2025} holds in its entire stated parameter range.
\end{corollary}
\begin{proof}
The positive constant time scaling and invertible linear change
\[
 X=a^{1/(2m)}x,\qquad V=a^{(m+1)/(2m)}v,\qquad \tau=\sqrt a\,t
\]
produce~\eqref{ap:normalized-jmaa} with
$\beta=b a^{(1-2n)/(2m)}\in(0,1)$ and
$\xi=\delta a^{(2n-1-m)/(2m)}>0$.
They preserve periodic orbits, stability and multipliers. The scaling remains valid when $m=2n-1$. Theorem~\ref{ap:jmaa} applies. Since the original conjecture assumes $n<m$ and has a strictly positive homoclinic threshold, its full parameter range is included.
\end{proof}
The factor $\xi=a^{(2n-1-m)/(2m)}$ in \cite[(4.2)]{ChenLiXiaoYang2025} corresponds to $\delta=1$ in the scaling above; retaining $\delta>0$ also covers the damping factor in that paper's introductory equation. Lemmas~4.9--4.10 of that paper establish uniqueness for $n=m$ and $n=m+1$ when $1/(2n)<\beta<1$. The present argument treats these cases together with $n<m$ and extends the upper bound to $0<\beta<1$. For $K=0$, the sign condition reduces to monotonicity of $f/g$ on the two sides of the left well. Existence throughout the larger interval $0<\beta<1$ is a separate question.

If a negative overall damping factor is allowed in~\eqref{ap:original-jmaa}, the change $(t,V)\mapsto(-t,-V)$ replaces $\delta$ by $-\delta$. It preserves the number of periodic orbits and reverses their stability. The corresponding upper bound therefore also holds for $\delta<0$, with any cycle repelling. At $\delta=0$ the system is Hamiltonian and its potential wells contain period annuli, so the bound on nonconstant periodic orbits no longer holds.

\subsection{Further examples and applications}\label{ap:other-section}
The first examples compare the uniqueness criteria and show that the two-cycle bound can be attained, including for coefficients of low regularity. Further applications concern predator--prey models and polynomial systems.

\subsubsection{Potential perturbations and uniqueness}
Adding $aG+bG^2$ to $F$ adds $a+2bu$ to each profile slope. Consequently
\begin{equation}\label{ap:gauge}
 \widetilde{\mathcal U}_{K+2b}=\mathcal U_K,\qquad
 \widetilde{\mathcal N}=\mathcal N.
\end{equation}
Indeed $\widetilde F'/g=F'/g+a+2bG$; one differentiation proves the first identity, and a third profile derivative proves the second. The perturbations preserve these two expressions, although they generally change the phase curves.

The following three families satisfy Theorem~\ref{cr:unique} on the entire real axis:
\begin{align}
 g&=x,& F&=\beta x^3-\alpha\arctan x+Ax^2+Bx^4,
 &\mathcal U_{8B}&=3\beta x^2+\frac{\alpha(1+3x^2)}{(1+x^2)^2}>0;
 \label{ap:arctan}\\
 g&=x+x^3,& F&=-\alpha x+\beta x^5+AG+BG^2,
 &\mathcal U_{2B}&=\alpha(1+3x^2)+\beta(15x^4+5x^6)>0;
 \label{ap:cubic-g}\\
 g&=x^{2m+1},&F&=-x+x^{2m+3}+aG+bG^2,
 &\mathcal U_{2b}&=(2m+1)x^{2m}+(2m+3)x^{4m+2}>0.
 \label{ap:degenerate}
\end{align}
In the first two rows $\alpha,\beta>0$ and $A,B\in\mathbb R$; in the last $m$ is any nonnegative integer, $a,b\in\mathbb R$, and $x\ne0$. In every row $G=\int_0^xg$, all coefficients are smooth, and $xg(x)>0$. Thus there is at most one periodic orbit, necessarily hyperbolic and attracting. The identities follow by differentiating the displayed $F$ and applying~\eqref{ap:gauge}. For example, in~\eqref{ap:degenerate}, $F_0'=-1+(2m+3)x^{2m+2}$, which gives the displayed numerator. The origin has degenerate restoring force when $m\ge1$.

These certificates allow arbitrary even potential perturbations. For instance, in~\eqref{ap:degenerate},
\[
 (F'/g)'=\frac{2m+1}{x^{2m+2}}+(2m+3)+2bx^{2m+1}.
\]
If $b\ne0$, this derivative has opposite signs near zero and on one infinite tail. Thus the global constant-sign condition on $(f/g)'$ fails, while the common-curvature condition remains valid.

For~\eqref{ap:degenerate} the upper bound is attained when $|a|+|b|$ is sufficiently small (with $m$ fixed). To construct a trapping annulus, first take $a=b=0$. Then $E=G+y^2/2$ satisfies
$\dot E=x^{2m+2}(1-x^{2m+2})$. A small regular level $E=e<G(1)$ is an outward inner boundary, including its two axis crossings, since $\dot x=y\ne0$ there and $E$ increases on both neighboring pieces. For an outer boundary put
\[
 c(x)=\arctan x,\quad D(x)=F_0(x)-c(x),\quad
 V=G+\frac12(y-c(x))^2.
\]
Direct completion of the square gives
\begin{equation}\label{ap:trapping}
 \dot V=-c'(x)\left(y-\frac{F_0(x)+c(x)}2\right)^2
             -g(x)D(x)+\frac{c'(x)}4D(x)^2.
\end{equation}
For $|x|\ge2$, $gD>0$ and $0<c'D/g<1$, because for $x\ge2$,
$0<D<x^{2m+3}$ and $c'=1/(1+x^2)$; oddness handles $x\le-2$.
Hence $\dot V<-3gD/4<0$ there. On $|x|\le2$, $c'\ge1/5$ and the remaining coefficients are bounded, so $\dot V<0$ for sufficiently large $|y|$. The function $V$ is proper and its gradient vanishes only at the origin. A sufficiently large level $V=L$ is a regular inward boundary enclosing the small level of $E$. The intervening compact positively invariant annulus has no equilibrium. The Poincare--Bendixson theorem supplies a periodic orbit; the already proved upper bound and hyperbolicity make it the unique attracting limit cycle. Its return-map fixed point persists for small $(a,b)$ by the implicit-function theorem. The global upper bound excludes additional cycles.

\subsubsection{The return map for small damping}
The small-damping existence arguments below use the following return expansion. For $g=x$ and $F=\varepsilon P(x)$, put $x=r\cos\theta$, $y=-r\sin\theta$. Then
\begin{equation}\label{ap:polar}
 \dot r=-\varepsilon\cos\theta P(r\cos\theta),\qquad
 \dot\theta=1+\frac{\varepsilon}{r}\sin\theta P(r\cos\theta).
\end{equation}
On a fixed compact annulus, the second expression is positive for small $|\varepsilon|$. If $P$ is $C^1$, the radial equation $r_\theta=\varepsilon a_\varepsilon(r,\theta)$ has $a_\varepsilon\to a_0=-\cos\theta P(r\cos\theta)$ in $C^1$. Its solutions satisfy $r(\theta)=r(0)+O(\varepsilon)$ and its initial-value derivative satisfies
$J_\theta=\varepsilon(\partial_ra_\varepsilon)J$, $J(0)=1$, whence $J=1+O(\varepsilon)$. Integrating these equations proves
\begin{equation}\label{ap:return}
 \frac{\Pi_\varepsilon(r)-r}{\varepsilon}\longrightarrow-I(r)
 \quad\hbox{in }C^1,\qquad
 I(r)=\int_0^{2\pi}\cos\theta P(r\cos\theta)\,d\theta.
\end{equation}
All estimates are uniform on compact parameter sets. A simple zero $\rho$ of $I$ gives a unique nearby fixed point, either by the implicit-function theorem in the smooth case or by endpoint signs and strict derivative bounds in the $C^1$ case. Its multiplier is $1-\varepsilon I'(\rho)+o(\varepsilon)$, hence, for sufficiently small nonzero $\varepsilon$, it is attracting if $\varepsilon I'(\rho)>0$ and repelling if $\varepsilon I'(\rho)<0$. These arguments take place in a slightly larger compact annulus, ensuring existence until one full turn.

\subsubsection{Finite potential endpoints}
For fixed real $A,B$ consider
\[
 \begin{gathered}
 z(x)=\frac{x}{\sqrt{1+x^2}},\quad g(x)=\frac{x}{(1+x^2)^2},\quad
 G(x)=\frac{x^2}{2(1+x^2)},\\
 F(x)=\varepsilon\left[-z+\frac{16}{3}z^3+Az^2+Bz^4\right],\qquad\varepsilon>0.
 \end{gathered}
\]
Both potential endpoints equal $1/2$. Since $g=zz'$ and $z'=(1+x^2)^{-3/2}$,
\[
 \mathcal U_{8\varepsilon B}(x)
 =\frac{\varepsilon(1+16z(x)^2)}{(1+x^2)^{9/2}}>0.
\]
Thus Theorem~\ref{cr:unique} gives at most one periodic orbit for every $\varepsilon>0$. The diffeomorphism $x\mapsto z$ maps $\mathbb R$ onto $(-1,1)$, and the positive time change $d\tau/dt=z'$ gives $z_\tau=y-\varepsilon P(z)$, $y_\tau=-z$ on that strip. Formula~\eqref{ap:return} gives $I(r)=\pi r(-1+4r^2)$, whose zero $r=1/2$ is simple with $I'=2\pi$. For sufficiently small $\varepsilon>0$ the resulting hyperbolic attracting orbit lies in $1/4<r<3/4$, hence stays strictly inside the strip. The inverse map and reciprocal time factor are bounded on that compact orbit, so it yields a finite-period orbit in the original plane. There is therefore exactly one limit cycle for all sufficiently small positive $\varepsilon$.

Finite potential endpoints have appeared in earlier uniqueness theory, for example \cite{CarlettiVillari2005}. Condition (E) of its Theorem~1.1 requires $G+F$ to be unbounded above on every right-hand tail and $G-F$ to be unbounded above on every left-hand tail. Both $F$ and $G$ here are bounded; the preceding return-map proof supplies existence independently of that condition.

\subsubsection{An asymmetric family attaining the two-cycle bound}
\begin{proposition}\label{ap:explicit-two}
Let $g=x$ and
\begin{equation}\label{ap:two-family}
 F=\varepsilon\left(4x-\frac{20}{3}x^3+\frac85x^5+Ax^2+Bx^4+
                  \eta x^6+\lambda x^7+\delta\sinh x\right).
\end{equation}
For $A,B\in\mathbb R$, $\varepsilon>0$, $\lambda,\delta\ge0$ and
$\eta^2\le35\lambda/8$, there are at most two limit cycles. For fixed $A,B$, there is $\varepsilon_0(A,B)>0$ such that there are exactly two hyperbolic limit cycles whenever additionally
$0\le\lambda\le1/8$, $0\le\delta\le1$ and $0<\varepsilon<\varepsilon_0(A,B)$. The inner cycle is repelling and the outer cycle attracting.
\end{proposition}
\begin{proof}
Differentiation gives
\[
 \frac{\mathcal N}{\varepsilon}
 =12+20x^2+x^4(24+48\eta x+105\lambda x^2)
    +\delta[(x^2+3)\cosh x-3x\sinh x]>0.
\]
The quadratic factor is nonnegative since its discriminant is nonpositive, with the case $\lambda=0$ forcing $\eta=0$. The last bracket is bounded below by
$\cosh x[(|x|-3/2)^2+3/4]>0$. Theorem~\ref{cr:two} gives the global upper bound without a small-parameter restriction.

For existence, the even powers integrate to zero in~\eqref{ap:return}, whereas
\[
 I(r)=\pi rQ(r^2),\quad Q(z)=(z-1)(z-4)+\frac{35\lambda}{64}z^3+\delta S(z),\quad
 S(z)=\sum_{j=0}^\infty\frac{z^j}{4^j j!(j+1)!}.
\]
The series follows by integrating the power series of $\sinh$; its derivatives converge uniformly on compact sets. We have $Q''>0$ for $z\ge0$, $Q(0)>0$, $Q(1)\ge0$, $Q(4)\ge0$, and
\[
 S(2)<\frac43,\qquad Q(2)<-2+\frac{35}{64}+\frac43=-\frac{23}{192}<0.
\]
For the bound on $S(2)$ use $(j+1)!\ge2^j$ and sum the resulting exponential series, which is strictly smaller than $\sum_{j\ge0}(1/4)^j$. Strict convexity implies exactly two positive simple zeros $1\le z_-<2<z_+\le4$, with $Q'(z_-)<0<Q'(z_+)$. These zeros and their derivative bounds vary continuously over the compact allowed $(\lambda,\delta)$ set and remain uniformly separated. Applying~\eqref{ap:return} uniformly on that set, including its compact allowed $\eta$ range, produces two simple return fixed points for one common $\varepsilon_0(A,B)$. The multiplier signs give the asserted stability. The global upper bound excludes any additional limit cycles.
\end{proof}
The family is polynomial of degree seven when $\delta=0$, $\lambda>0$, and nonpolynomial when $\delta>0$. The two-cycle bound is restricted to the coefficient conditions in the proposition.

\subsubsection{Low regularity and slope interfaces}
First, Theorem~\ref{cr:two} applies beyond classical third derivatives. With
$P_0(x)=4x-20x^3/3+8x^5/5$, set
\begin{equation}\label{ap:hinge}
 g=x,\quad F=\varepsilon\left[P_0(x)+Ax^2+Bx^4+
       \delta\operatorname{sgn}(x)(x^2/2-1)_+^2\right],\qquad\delta\ge0.
\end{equation}
The slopes of the potential profiles are
\[
 k_\pm(u)=\varepsilon\left[2A+8Bu\pm(k_0(u)+2\delta(u-1)_+)\right],
 \quad k_0=4(2u)^{-1/2}-20(2u)^{1/2}+8(2u)^{3/2}.
\]
Here $k_0''=(12+40u+96u^2)/(2u)^{5/2}>0$. Thus $k_+$ is convex and $k_-$ concave, giving at most two limit cycles for every $\varepsilon>0$.

For any fixed $A,B$ and any bounded interval $0\le\delta\le D$, this bound is attained for all sufficiently small $\varepsilon>0$. In~\eqref{ap:return},
\[
 I_\delta(r)=\pi r[(r^2-1)(r^2-4)+\delta S(r^2)],\quad
 S(z)=\frac1{\pi\sqrt z}\int_0^{2\pi}|\cos\theta|
                (z\cos^2\theta/2-1)_+^2\,d\theta.
\]
For $0\le a\le1/2$, $h_a(z)=z^{-1/2}(az-1)_+^2$ is $C^1$ and convex: on $az>1$,
\[
 h_a''=\frac{3a^2}{4\sqrt z}+\frac{a}{2z^{3/2}}+\frac3{4z^{5/2}}>0,
\]
and its derivative is zero at the junction. Integration proves convexity of $S$. Moreover $S=0$ on $(0,2]$ and $S>0$ above $2$. Hence the bracket is strictly convex, with a simple zero at $1$, value $-2$ at $2$, and nonnegative value at $4$; it has one further simple zero $z_+(\delta)\in(2,4]$. Compactness in $\delta$ and the $C^1$ return calculation give two hyperbolic cycles, inner repelling and outer attracting. For $\delta>0$, $F$ is $C^1$ but not $C^2$ at $x=\pm\sqrt2$, and the outer orbit crosses both lines transversely since its limiting radius is $\sqrt{z_+(\delta)}>\sqrt2$. The curvature jump satisfies the convex-slope condition at the junctions as well as on the smooth pieces.

A different low-regularity example is obtained for $0<p<1$ by
\begin{equation}\label{ap:pullback}
 \begin{gathered}
 g(x)=\operatorname{sgn}(x)|x|^p,\qquad
 Z(x)=\sqrt{\frac2{p+1}}\operatorname{sgn}(x)|x|^{(p+1)/2},\\
 F(x)=\varepsilon[P_0(Z(x))+AZ(x)^2+BZ(x)^4].
 \end{gathered}
\end{equation}
There are at most two limit cycles for all $\varepsilon>0$, and exactly two hyperbolic cycles for $0<\varepsilon<\varepsilon_0(A,B)$, with the same threshold for every $p\in(0,1)$. To justify the possibly singular change, note that $Z$ is an increasing homeomorphism, $Z'=g/Z>0$ off zero, and the oriented phase curves there transform to the normalized system of Proposition~\ref{ap:explicit-two}. At an axis crossing $(0,y_0)$ with $y_0\ne0$, the graph equation
$dy/dx=g(x)/(F(x)-y)$ has a continuous right side locally Lipschitz in $y$, with denominator bounded away from zero; it has a unique crossing and the physical time integral $\int dx/(y-F(x))$ is finite. The normalized equation has the same unique continuation, since $dy/dZ=Z/(\widehat F(Z)-y)$ extends continuously across the axis. Nearby return maps, using $y$ on the axis, are identical. Finally $gF>0$ near zero, so $E=G+y^2/2$ is nonincreasing there and a forward solution cannot leave the equilibrium; a nonconstant periodic orbit cannot pass through it. This proves that the correspondence covers every periodic orbit, and preserves multipliers. Both $g$ and $F$ fail to be locally Lipschitz at zero; nonzero $(A,B)$ also destroys oddness. Thus the normalized example transfers its cycle count and multipliers to these non-Lipschitz coefficients.

For the separate interface theorem, let $W_1(x)=\operatorname{sgn}(x)(|x|-1)_+$ and
\begin{equation}\label{ap:interface-example}
 g=x,\qquad F=\varepsilon[-x+x^3+Ax^2+W_1(x)].
\end{equation}
For every $\varepsilon>0$ there is at most one periodic orbit, necessarily attracting and hyperbolic. Indeed, with $K=0$, $F'/x$ has derivative $\varepsilon(x^{-2}+3)$ inside $|x|<1$ and $3\varepsilon$ outside; its jumps, in increasing $x$, are both $+\varepsilon$. Theorem~\ref{cr:interface} applies. Existence for small positive $\varepsilon$ follows from
\[
 \begin{gathered}
 I(r)=\pi rQ(r),\qquad Q(r)=-1+3r^2/4+S(r),\\
 S(r)=\frac2\pi\left[\arccos(1/r)-r^{-1}\sqrt{1-r^{-2}}\right]\quad(r>1),
 \end{gathered}
\]
with $S=0$ for $r\le1$. Since $S'=4\sqrt{1-r^{-2}}/(\pi r^2)>0$, there is a unique simple root $1<\rho<2/\sqrt3$. Near these circles every crossing of $x=\pm1$ is transverse. The smooth flows within the finitely many pieces and the implicit-function theorem for crossing times yield a $C^1$ return expansion~\eqref{ap:return}. The resulting hyperbolic attracting cycle crosses both corners.

For continuous piecewise-linear $F$, the global common-curvature condition is more restrictive. If $F$ is globally continuous with finitely many linear pieces and $F'/x-Kx^2/2$ is nondecreasing on both half-lines, its two linear tails force $K=0$. Then every segment slope is nonpositive. Since $F(0)=0$, $-xF(x)\ge0$, so energy is nondecreasing. Any periodic orbit must lie where $F=0$ and belongs to a family of harmonic circles; there are no limit cycles. The nonincreasing version follows by reversing signs. The continuous and possibly non-Lipschitz coefficients above also differ from the Lipschitz odd-$F$ assumptions in \cite[Theorem~1.1]{ChenHanXia}; conversely that theorem permits some discontinuous restoring forces outside our assumptions.

\subsubsection{Comparison of the two uniqueness criteria}
For $x\ne0$ write
\[
 \mathcal C_d(x)=\frac{d}{dx}\left[\frac{F'(x)/g(x)-d}{G(x)}\right].
\]
Its sign equals that of the intercept numerator, since the multiplying denominator is $g^2G^2>0$. For $g=x$, $F=\sum_{j=1}^5a_jx^j$,
\[
 \mathcal C_{2a_2}=\frac2{x^4}(-3a_1-3a_3x^2+5a_5x^4).
\]
Thus $a_1,a_3<0$, $a_5\ge0$ imply at most one periodic orbit for arbitrary $a_2,a_4$; any such orbit is hyperbolic and attracting. The condition $a_5\ge0$ is needed for this global sign test.

Neither uniqueness criterion implies the other. Let $P(x)=-x-10x^3+x^5$. Then
$\mathcal C_0=2(3+30x^2+5x^4)/x^4>0$, but
$\mathcal U_K(\pm1)=-14\mp K$ while $\mathcal U_K(0)=1$.
No constant $K$ gives either global sign. Also
$\mathcal N=-3+30x^2+15x^4$ changes sign, excluding the direct two-cycle test. Conversely, $F=-x+x^3$, $g=x$ satisfies $\mathcal U_0=1+3x^2>0$, but
$\mathcal C_d=(6-6x^2+4dx)/x^4$ is negative on both distant tails for every $d$.

A nonpolynomial example is
\[
 \Psi(x)=P(x)+\sinh x-x-x^3/6.
\]
Termwise differentiation of its locally uniformly convergent series gives
\[
 \mathcal C_0=\frac6{x^4}+\frac{60}{x^2}+10+
       2\sum_{j=2}^\infty\frac{(2j-3)x^{2j-4}}{(2j)!}>0.
\]
Here $\mathcal U_K(\pm1)=-27/2-e^{-1}\mp K$ and its near-zero limit is $1$, while $\mathcal N$ is negative near zero and positive at $1$. Thus Theorem~\ref{cr:intercept} applies, although neither of the other two sign conditions holds.

Both $L=P$ and $L=\Psi$ have exactly one hyperbolic attracting limit cycle for $g=x$. For existence use the annulus proof of~\eqref{ap:trapping}, now taking $c(x)=\int_0^xe^{-s^2}ds$. The small boundary is outward since $L=-x+O(x^3)$. At large $|x|$, $x(L-c)>0$ and $0<e^{-x^2}(L-c)/x<1$, for both polynomial and exponential growth. Thus~\eqref{ap:trapping} gives an inward large regular level, and Poincare--Bendixson gives a periodic orbit. Theorem~\ref{cr:intercept} supplies uniqueness and hyperbolicity. In contrast $L=-x$ passes the same criterion but has no periodic orbit because $\dot E=x^2$.

The criterion also allows nonlinear restoring forces and unequal finite endpoints. For an increasing smooth diffeomorphism $z=z(x)$ with $z(0)=0$ put
\[
 g=zz',\quad G=z^2/2,\quad F=L(z)+aG+bG^2,\quad L=P\ \hbox{or}\ \Psi.
\]
Then $F'(0)=-z'(0)<0$ and
$\mathcal C_a(x)=z'(x)\,\mathcal C_0^{L}(z(x))>0$.
Theorem~\ref{cr:intercept} gives the upper bound for all $a,b$. Choices $z=x\sqrt{1+x^2/2}$ and $z=2\sinh(x/2)$ give $g=x+x^3$ and $g=\sinh x$; at $a=b=0$ the positive time change $d\tau/dt=z'$ transfers the existence proof above on the whole normalized plane. By contrast,
$z=\tanh x/(1+\tanh x/2)$ has range $(-2,2/3)$ and endpoint energies $2$ and $2/9$; the same upper bound holds, but the preceding existence argument does not apply on this strip. A smooth degenerate example is
\[
 g=x^3,\quad F=-x-10x^7+x^{11},\qquad
 \mathcal C_0=28x^{-8}+280x^{-2}+132x^2>0.
\]
It has $F'(0)<0$ and $g'(0)=0$; the degenerate version of Theorem~\ref{cr:intercept} applies.

Finally a common intercept cannot be replaced by unordered constants for the two sides. On $(-4,4)$ set
\[
 F_\varepsilon=\varepsilon p(x),\quad
 p=-36x+\frac{196}{3}x^3-\frac{112}{5}x^5+\frac{64}{35}x^7,\quad g=x.
\]
Then
\[
 \mathcal C_d=\varepsilon(216x^{-4}-392x^{-2}-224+384x^2/5)+4d/x^3.
\]
For $d_+=4000\varepsilon$ and $0<x<4$,
$x^3\mathcal C_{d_+}/(4\varepsilon)
=54/x-98x-56x^3+96x^5/5+4000>24$;
reflection gives positivity with $d_-=-4000\varepsilon$ on the other side. But~\eqref{ap:return} gives
\[
 I(r)=\pi r(r^2-1)(r^2-4)(r^2-9),
\]
with simple roots $1,2,3$, so there are at least three hyperbolic cycles for small positive $\varepsilon$, all within the strip. Thus uniqueness fails if the two constants are allowed to have this order.

\subsubsection{An Ivlev--Leslie--Gower model: a global two-cycle parameter region}
Liao, Zegeling and Huang \cite[Section~4.1]{LZH2024} consider
\begin{equation}\label{ap:ivlev}
 \dot x=rx(1-x/k)-xe^{-\beta x}Y,\qquad
 \dot Y=Y\left(D-\frac{\mu Y}{x+\lambda}\right),\qquad x,Y>0,
\end{equation}
and conjecture on p.~33627 that it has at most two nested limit cycles for all allowed parameters. Their Proposition~4.1 already produces two small cycles at a degenerate Hopf point. We establish a global upper bound near that point, and extend the resulting exact count to an open set with all six parameters strictly positive.

\begin{proposition}\label{ap:ivlev-two}
There is a nonempty open subset of the parameter space $r,k,\mu,D,\beta,\lambda>0$ on which~\eqref{ap:ivlev} has exactly two limit cycles in the whole positive quadrant. Both are hyperbolic; the inner one is repelling and the outer one attracting.
\end{proposition}
\begin{proof}
We first take $\lambda=0$, $\beta=1/4$, $D=1$ and
\[
 \mu=\frac{ke^{-1/4}}{r(k-1)},\quad 1<k<16,\quad
 k_* =\frac{35+3\sqrt{137}}{10},\quad r_* =\frac{4k_*}{k_*-5}.
\]
The equilibrium abscissa is then $1$. The point $(r_*,k_*)$ is the degenerate point already used in \cite{LZH2024}.

To reduce the system to ordinary Li\'enard form, take $d\tau/dt=xe^{-x/4}>0$ and put
\[
 H=re^{x/4}(1-x/k),\quad c=e^{x/4}/x,\quad d=\mu e^{x/4}/x^2,
 \quad A=-H'-c+2dH,\quad B=H(c-dH).
\]
The equations become $x'=H-Y$, $Y'=cY-dY^2$. Differentiating and eliminating $Y$ gives $x''=d(x')^2-Ax'-B$.
Set $m(x)=\exp(-\int_1^xd(t)dt)$ and $u=\psi(x)=\int_1^xm(t)dt$. Then
\[
 u''+\widetilde f(u)u'+\widetilde g(u)=0,\qquad
 \widetilde f(\psi(x))=A(x),\quad\widetilde g(\psi(x))=m(x)B(x).
\]
The phase map $(x,Y)\mapsto(\psi(x),m(x)[H(x)-Y])$ has determinant $-m^2$ and is one-to-one. Every positive periodic orbit lies in $0<x<k$: at its maximum abscissa, $r(1-x/k)=e^{-x/4}Y>0$. Moreover
\[
 (\log(H/x))'=\frac14-\frac1{k-x}-\frac1x\le\frac14-\frac4k<0.
\]
Since $\mu H(1)=1$, we obtain $(x-1)B(x)>0$. Thus the complete strip $\psi((0,k))\times\mathbb R$ is a single well and contains the images of all positive periodic orbits. On each such compact orbit both time factors are bounded above and below by positive constants. The transformation preserves periodicity, isolation and multipliers; no extension to $x=0$ is needed.

Using $d/du=m^{-1}d/dx$ and $m'=-dm$, the exact third-order certificate is
\begin{equation}\label{ap:ivlev-N}
 \widetilde{\mathcal N}
 =B^2A''+(4dB^2-3BB')A'
 +[3(B')^2-5dBB'+(3d^2+d')B^2-BB'']A.
\end{equation}
An explicit finite certificate, given after this proof, establishes at the base point
\begin{equation}\label{ap:ivlev-base}
 \widetilde{\mathcal N}(x)>0\quad(0<x\le176/25,\ x\ne1),\qquad
 \widetilde{\mathcal N}(x)=c_4(x-1)^4+O((x-1)^5),\quad c_4>0.
\end{equation}
The enlarged endpoint $176/25>k_*$ leaves room for small changes in $k$.

To extend the sign condition to nearby parameters, normalize the restoring force locally by
$z=\operatorname{sgn}(u)\sqrt{2\int_0^u\widetilde g}$.
Since $v:=\widetilde g'(0)=B'(1)>0$, this is an analytic diffeomorphism through zero: $z=u\sqrt{R(u)}$, where $R(0)=v$ and $R(u)=2u^{-2}\int_0^u\widetilde g$ is positive analytic near zero. With the positive analytic time factor $\alpha=dz/du=\widetilde g/z$, the mechanical equation becomes
$\ddot z+\widehat f(z)\dot z+z=0$, where $\widehat f=\widetilde f/\alpha=\sum_{j\ge0}A_jz^j$.
At $x=1$ define
\[
 v=B',\quad w=B''-dB',\quad t=B'''-(2d'+d^2)B',\qquad
 a_0=A,\quad a_1=A',\quad a_2=A''+dA'.
\]
Taylor inversion of $z=u\sqrt{R(u)}$ yields
\begin{equation}\label{ap:ivlev-shape}
 A_0=\frac{a_0}{\sqrt v},\qquad
 A_2=\frac{12v^2a_2-12vwa_1+(5w^2-3vt)a_0}{24v^{7/2}}.
\end{equation}
The certificate transforms as
$\widetilde{\mathcal N}=\alpha^5\widehat{\mathcal N}$, and
\[
 \widehat{\mathcal N}=3A_0-A_2z^2+3A_4z^4+O(z^5).
\]
At the base point~\eqref{ap:ivlev-base} gives $A_0=A_2=0$ and $c_*:=A_4>0$. On a common neighborhood, choose $A_4\ge2c_*/3$ and bound the analytic remainder by $c_*z^4$. Thus $A_0\ge0$, $A_2\le0$ imply
$\widehat{\mathcal N}\ge3A_0-A_2z^2+c_*z^4\ge0$ near zero.
Away from $x=1$, use compactness after treating the singular endpoint: $x^{10}\widetilde{\mathcal N}$ extends analytically through $x=0$ and at the base point has value
\[
 \lim_{x\to0^+}x^{10}\widetilde{\mathcal N}
 =\frac{6k_*^5e^{-5/4}}{(k_*-1)^5}>0.
\]
Together with~\eqref{ap:ivlev-base}, this gives a positive minimum on the remaining closed set $[0,176/25]\setminus(1-\eta,1+\eta)$, which persists for nearby parameters. Theorem~\ref{cr:two} therefore gives a global bound of two in the region $A_0\ge0,A_2\le0$.

The two shape parameters are independent. Define
\[
 T=1-\frac{r(k-5)}{4k},\qquad
 C=e^{-3/4}[B'A''-B''A'+2dB'A']_{x=1}.
\]
Then $A(1)=e^{1/4}T$, and on $T=0$,
\[
 C=-\frac{25k^2-175k-2}{4(k-5)^2},\qquad
 \det\frac{\partial(T,C)}{\partial(r,k)}\bigg|_*
 =\frac{50k_*-175}{16k_*(k_*-5)}>0.
\]
On this same curve $A_2=e^{3/4}C/(2v^{5/2})$, so $(A_0,A_2)$ are analytic local coordinates. Choose a small $\epsilon>0$ and parameters satisfying
\[
 A_0=c_*\epsilon^4/2,\qquad A_2=-5c_*\epsilon^2/2.
\]
Then $(r,k)-(r_*,k_*)=O(\epsilon^2)$ and $A_4=c_*+O(\epsilon^2)$.
Scale $z=\epsilon Z$, $\dot z=\epsilon V$. The odd part of the damping is $O(\epsilon)$ and must be retained: the system with just this odd damping is reversible under $(Z,V,t)\mapsto(-Z,V,-t)$. On the section $Z=0$, $V=\rho\in[1/2,5/2]$, its first-return map is the identity. Indeed it is an $O(\epsilon)$ perturbation of the harmonic oscillator on a fixed compact annulus, so every such orbit has a transverse opposite-axis crossing; reflection of that half-orbit completes a periodic orbit. Its trajectory and first variational equation differ from the harmonic ones by $O(\epsilon)$ uniformly in $C^1$ with respect to $\rho$.

The even damping is
$\epsilon^4c_*(1/2-5Z^2/2+Z^4)+O_{C^2}(\epsilon^6)$.
Multiply it by an auxiliary parameter $\theta\in[0,1]$. In vector notation the resulting equation is
$\dot W=X^{\rm odd}_\epsilon(W)+\theta\epsilon^4Y_\epsilon(W)=X_{\epsilon,\theta}(W)$,
where the coefficients and their required derivatives are uniformly bounded on a common compact flow tube. Its first and second $\theta$-variations, with fixed initial data, satisfy
\[
 \begin{aligned}
 \dot Q_1&=DX_{\epsilon,\theta}Q_1+\epsilon^4Y_\epsilon,
       &Q_1(0)&=0,\\
 \dot Q_2&=DX_{\epsilon,\theta}Q_2+
 D^2X_{\epsilon,\theta}[Q_1,Q_1]+2\epsilon^4DY_\epsilon Q_1,
       &Q_2(0)&=0.
 \end{aligned}
\]
Gronwall's inequality, followed by differentiation with respect to the initial ordinate $\rho$, gives $Q_1=O_{C^1}(\epsilon^4)$ and $Q_2=O_{C^1}(\epsilon^8)$. The return section is uniformly transverse, so implicit differentiation of its crossing equation gives the same bounds for the return time. At $\theta=0$ the return map is the identity. Replacing that reversible base orbit and its variational equation by the harmonic ones changes the first variation by only $O_{C^1}(\epsilon^5)$, since both differ by $O(\epsilon)$. Taylor's formula in $\theta$ therefore gives
\[
 P_\epsilon(\rho)-\rho
 =-\frac{\pi c_*}{8}\epsilon^4\rho(\rho^2-1)(\rho^2-4)
       +O_{C^1}(\epsilon^5).
\]
The coefficient follows by integrating $-\rho c_*(1/2-5\rho^2\cos^2\theta/2+\rho^4\cos^4\theta)\sin^2\theta$ and using the integrals $\pi,\pi/4,\pi/8$. Its two simple zeros give two small hyperbolic cycles in the positive quadrant, the inner one repelling and the outer one attracting. The global upper bound excludes additional limit cycles.

It remains to move into strictly positive $\lambda$. For general nearby parameters put
\[
 \begin{gathered}
 H=re^{\beta x}(1-x/k),\quad c=De^{\beta x}/x,\quad d=\mu e^{\beta x}/[x(x+\lambda)],\\
 A=-H'-c+2dH,\qquad B=H(c-dH).
 \end{gathered}
\]
The same transformations and~\eqref{ap:ivlev-N} apply. Near the corner, we use an estimate uniform as $(x,\lambda)\to(0,0)$ with $x>0$ and $\lambda\ge0$. Write $q=x+\lambda$,
$L_0=2H-Dq/\mu$, $M_0=H(Dq/\mu-H)$, so $A=dL_0-H'$, $B=dM_0$. Direct differentiation gives
\[
 \frac1d=O(q^2),\quad
 \frac{d'}{d^2}=\frac{\beta xq-q-x}{\mu e^{\beta x}}=O(q),\quad
 \frac{d''}{d^3}=\frac{(\beta xq-q-x)^2+q^2+x^2}{\mu^2e^{2\beta x}}=O(q^2).
\]
All bounds are uniform in a fixed compact parameter box with $r,\mu>0$. Substitution in~\eqref{ap:ivlev-N}, using bounded derivatives of $L_0,M_0,H$, gives
\begin{equation}\label{ap:ivlev-corner}
 \widetilde{\mathcal N}/d^5=3M_0^2L_0+O(q)=6r^5+O(q)>0
 \quad(0<x+\lambda<\delta)
\end{equation}
for one fixed sufficiently small $\delta>0$. The estimate follows from
$A/d=L_0+O(q^2)$, $A'/d^2=O(q)$, $A''/d^3=O(q^2)$,
$B/d=M_0$, $B'/d^2=O(q)$, $B''/d^3=O(q^2)$; the only surviving order-one term is $3d^2B^2A/d^5$.

Fix $\beta=1/4,D=1$ and take $\mu=(1+\lambda)/H(1)$. The nonzero Jacobian above continues the base point as $r_*(\lambda),k_*(\lambda)$ with $A_0=A_2=0$, $A_4>0$. For small $\lambda\ge0$, the complete interval $0<x<k$ is covered by three uniform regions: the corner~\eqref{ap:ivlev-corner}, a fixed neighborhood of the equilibrium controlled by~\eqref{ap:ivlev-shape}, and the remaining compact set controlled by strict positivity at $\lambda=0$. Moreover
\[
 (\log(H/(x+\lambda)))'=\frac14-\frac1{k-x}-\frac1{x+\lambda}
 \le\frac14-\frac4{k+\lambda}<0
\]
preserves the single well. Choose a fixed small $\lambda_0>0$ and repeat the small-$\epsilon$ construction with the continued $A_4>0$. This yields two hyperbolic cycles with $A_0>0,A_2<0$ and strictly positive $\widetilde{\mathcal N}$ throughout $0<x\le176/25$.
Finally permit all six parameters to vary. The corner estimate remains uniform, positivity persists on the remaining compact interval, and the strict conditions $\beta<4/(k+\lambda)$ and $\mu r>D\lambda$ preserve the unique simple positive restoring zero (which may move). Hyperbolicity preserves the two cycles, and the global criterion excludes additional ones. These conditions define an open set in the full six-dimensional parameter space.
\end{proof}

The finite certificate for~\eqref{ap:ivlev-base} uses exact arithmetic in $\mathbb Q(\sqrt{137})$. Take $k=k_*$, $r=4k/(k-5)$, and let $z$ be an algebraic placeholder for $e^{(x-1)/4}$. Set
\[
 H_0=rz(1-x/k),\quad c_0=z/x,\quad d_0=\frac{(k-5)z}{4(k-1)x^2},\quad
 D_x=\partial_x+\frac z4\partial_z,
\]
\[
 A_0^{\rm alg}=-D_xH_0-c_0+2d_0H_0,\quad
 B_0^{\rm alg}=H_0(c_0-d_0H_0).
\]
Let $N_0^{\rm alg}$ be the right side of~\eqref{ap:ivlev-N} with these three functions and $D_x$ in place of primes. Define the explicit polynomials $P_j$ by the finite identity
\[
 x^{10}(k-5)^3(k-1)^5z^{-5}N_0^{\rm alg}(x,z)=\sum_{j=0}^5P_j(x)z^j.
\]
Two differentiations of the displayed rational functions determine all coefficients. Each $P_j$ has degree at most $11$. Restoring the common exponentials gives
\[
 \widetilde{\mathcal N}(x)=\frac{e^{5x/4}W(x)}{x^{10}(k-5)^3(k-1)^5},\qquad
 W(x)=\sum_{j=0}^5P_j(x)e^{j(x-1)/4}.
\]
Write $P_j(1+h)=\sum_i p_{ji}h^i$ and
\[
 t_n=\sum_{j=0}^5\sum_{i\le n}p_{ji}\frac{(j/4)^{n-i}}{(n-i)!},\qquad
 V_{66}(h)=\sum_{n=4}^{70}t_nh^{n-4}.
\]
Reduction by $25k^2-175k-2=0$ gives $t_0=t_1=t_2=t_3=0$.
Write $p_{ji}=a_{ji}+b_{ji}\sqrt{137}$, $a_{ji},b_{ji}\in\mathbb Q$, and set $L=151/25$. The exponential Taylor remainder gives the rational bound
\[
 \left|W(1+h)/h^4-V_{66}(h)\right|
 \le\sum_{j=1}^5\sum_i (|a_{ji}|+12|b_{ji}|)L^{i-4}
             \frac{6561\,8^{71-i}}{(71-i)!}<10^{-12}
 \quad(|h|\le L).
\]
At $h=0$ the quotient has its removable value. For $i<4$, combine powers as $|h|^{i-4}|h|^{71-i}=|h|^{67}$ before estimating; no negative power is bounded separately. We used $jL/4<8$ and $e^8<3^8=6561$.

For each interval $[a,b]$ below, expand $V_{66}(a+(b-a)s)=\sum_{j=0}^{66}c_js^j$. Its Bernstein coefficients are
$b_i=\sum_{j\le i}c_j\binom ij/\binom{66}j$.
Using the rational enclosure
$q/10^{50}<\sqrt{137}<(q+1)/10^{50}$, $q=\lfloor\sqrt{137\,10^{100}}\rfloor$,
all coefficients have the following integer lower bounds:
\[
\begin{array}{c|r@{\qquad}c|r}
{[-1,-1/2]}&237034&[-1/2,0]&3049464\\
{[0,1/2]}&13805155&[1/2,1]&31356053\\
{[1,2]}&52791984&[2,3]&112578428\\
{[3,4]}&271887926&[4,5]&770976948\\
{[5,151/25]}&2196080761&&
\end{array}
\]
The Bernstein basis is nonnegative and sums to one, so $W(1+h)/h^4>0$ on $[-1,L]$. This proves~\eqref{ap:ivlev-base}, including its positive quartic coefficient. The local production of two cycles was established in \cite{LHF2018,LZH2024}. The global sign condition supplies the corresponding upper bound on the stated open parameter set.

\subsubsection{Further predator--prey applications}
For a system $\dot x=A_1(x)-B_1(x)Y$, $\dot Y=Y(C_1(x)-D_1(x)Y)$ with $B_1>0$, the same Li\'enard reduction applies: set $H=A_1/B_1$, $c=C_1/B_1$, $d=D_1/B_1$, change time by $d\tau/dt=B_1$, and use
\[
 L=-H'-c+2dH,\quad J=H(c-dH),\quad M=\exp(-\int d),\quad u=\int M.
\]
Then $\widetilde f=L$, $\widetilde g=MJ$ and the exact uniqueness numerator is
\begin{equation}\label{ap:ecology-U}
 \mathcal U_\kappa=L'J-L(J'-dJ)-\kappa M^3J^3.
\end{equation}
The third-order numerator is~\eqref{ap:ivlev-N} with $A=L,B=J$. The phase Jacobian is $-M^2\ne0$, as before. The cycle bounds then apply on a single-well strip; for the parameter choices below, the models admit such a strip containing every positive periodic orbit.

For the toxin-determined plant--herbivore model of \cite{Zhao2015},
\[
 \dot x=rx(1-x/K)-C(x)Y,\quad \dot Y=(\beta C(x)-d_0)Y,
 \quad C(x)=\frac{ex}{1+hex}\left(1-\frac{ex}{4G_0(1+hex)}\right),
\]
take
\begin{equation}\label{ap:toxin-parameters}
 (r,K,h,e,G_0,\beta,d_0)=(5/12,4,1,1,1/3,16/5,1).
\end{equation}
Here $C=x(x+4)/[4(x+1)^2]$,
$H=5(4-x)(x+1)^2/[12(x+4)]$,
$c=4(x-1)(5-x)/[5x(x+4)]$, and $d=0$ in~\eqref{ap:ecology-U}.
Every positive periodic orbit has $x<4$, by evaluating the first component at its maximum abscissa. On $(0,4)$, $H>0$ and $(x-1)c>0$, so there is a single well. The certificate is
$\mathcal U_0=P_8(x)/[x^2(x+4)^4]$.
In the degree-eight Bernstein basis on $[0,4]$, the coefficients of $P_8$ are exactly
\[
 \frac{680}3,\ \frac{1192}9,\ \frac{4384}{105},\ \frac{86504}{315},\
 \frac{290344}{525},\ \frac{9832}{105},\ \frac{394832}{315},\
 \frac{89528}9,\ \frac{23560}3.
\]
These positive coefficients specify the polynomial and prove its positivity throughout the interval. Theorem~\ref{cr:unique} gives at most one positive periodic orbit, hyperbolic and attracting if present. Existence follows from \cite[Theorem~2.2(a)]{Zhao2015}: its unimodal-response conditions are
$1/(4h)<G_0<1/(2h)$, a unique interior equilibrium, $H'(x_*)>0$ and $C'(x_*)\ne0$. Here $x_*=1$, $H(1)=1$, $1<K=4<5$, $H'(1)=7/15>0$, and $C'(1)=1/16>0$. Thus there is exactly one limit cycle.

This conclusion persists on an open set in the seven-dimensional space of the original positive parameters. The two simple roots of $\beta C-d_0$ keep their strict order $0<x_*<K<x^*$, and the cited strict existence conditions persist. For the global sign put
\[
 q=1+hex,\quad \ell=4G_0+(4G_0h-1)ex,\quad
 A=\frac{4G_0r}{Ke}(K-x)q^2,\quad B=\beta x\ell-\frac{4G_0d_0}{e}q^2.
\]
Then $H=A/\ell$, $c=B/(x\ell)$. Defining
$S=A'\ell-A\ell'$, $T=S'\ell-2S\ell'$ and
$V=B'x\ell-B(\ell+x\ell')$, direct differentiation in~\eqref{ap:ecology-U} gives
\[
 x^2\ell^5\mathcal U_0=xB(S^2-AT)+ASV+\ell SB^2.
\]
This polynomial has degree at most nine. At~\eqref{ap:toxin-parameters} it equals $(x+4)P_8(x)/243$, with all degree-nine Bernstein coefficients strictly positive. On the moving interval $[0,K]$ these coefficients depend continuously on the original parameters. Positivity therefore persists, as does $\ell>0$. The global uniqueness bound and the cited existence theorem therefore give exactly one cycle throughout an open neighborhood. This neighborhood covers only part of the parameter range in the branch-monotonicity conjecture of \cite[Section~2.5]{Zhao2015}.

For the Allee model with predator competition \cite{Gonzalez2019}, take the parameter slice
\begin{equation}\label{ap:allee}
 \dot x=x[(1-x)(x-\vartheta)-Y],\qquad
 \dot Y=Y(px-c_0-eY),\qquad -1<\vartheta<1.
\end{equation}
Here $H=(1-x)(x-\vartheta)$, $c=p-c_0/x$, $d=e/x$, $M=x^{-e}$,
\[
 L=2(1-e)x+(2e-1)(1+\vartheta)-p+(c_0-2e\vartheta)/x,
 \quad J=H[ex^2+(p-e(1+\vartheta))x+e\vartheta-c_0]/x.
\]
All positive cycles have $x<1$. If $\vartheta>0$, they also have $x>\vartheta$, since the first component is strictly negative at a potential minimum with $x\le\vartheta$. In the strong case below $c_0>p\vartheta$ and $p>c_0$, giving exactly one restoring root on $(\vartheta,1)$. For $\vartheta\le0$, $H>0$ on $(0,1)$ and the quadratic in $J$ has negative constant term and positive value at $1$, hence precisely one root in $(0,1)$. Thus each of the following three parameter choices gives a single well containing all positive periodic orbits:
\[
\begin{array}{c|c|c|c|c|c}
\vartheta&e&p&c_0&x_*&\kappa\\\hline
1/10&1/10&1&191/800&1/4&0\\
0&1/10&1&37/160&1/4&0\\
-1/10&1/3&1&3/25&1/5&10
\end{array}
\]
Formula~\eqref{ap:ecology-U} has respectively positive denominators $x^3,x^2,x^4$. The last uses $M^3=x^{-1}$. For row $i$ of the parameter table, put
$P_i(x)=x^{m_i}\mathcal U_{\kappa_i}(x)$, where
$(m_1,m_2,m_3)=(3,2,4)$. These polynomials have degrees
$(n_1,n_2,n_3)=(6,5,12)$, respectively. On an interval $[a,b]$, write
$P_i(a+(b-a)s)=\sum_{j=0}^{n_i}c_js^j$.
The power coefficients and Bernstein coefficients are given by
\[
 c_j=\frac{(b-a)^j}{j!}P_i^{(j)}(a),\qquad
 \beta_\ell=\sum_{j=0}^{\ell}c_j
 \frac{\binom{\ell}{j}}{\binom{n_i}{j}},\qquad 0\le\ell\le n_i.
\]
Using the rational parameters in the table, these formulas give the
following exact minimum Bernstein coefficients:
\[
\begin{array}{c|c|c}
\vartheta&\text{interval}&\text{minimum coefficient}\\\hline
1/10&[1/10,1]&38259/1280000000\\
0&[0,1/4],\ [1/4,1/2],\ [1/2,1]&
121203/10240000,\ 969/51200,\ 5327/102400\\
-1/10&[0,1/10],\ [1/10,1/5]&649957/20250000000,\ 273937/3125000000\\
-1/10&[1/5,1/2],\ [1/2,1]&1/6250,\ 52949/25000000
\end{array}
\]
The Bernstein basis is nonnegative and sums to one, so the listed minima prove positivity on each interval. Each of the three systems consequently has at most one positive periodic orbit, necessarily attracting and hyperbolic. The three cases correspond respectively to strong, zero and weak Allee effects; their positivity neighborhoods cover only part of the parameter space.

\subsubsection{Additional polynomial and Hamiltonian subfamilies}
Consider
\[
 \dot x=ax+by,\qquad \dot y=cx^3+dx^2y+exy^2,\qquad b\ne0,
\]
and set $A=bc-ad+ea^2/b$, $B=d-2ea/b$, $C=e/b$. The global mechanical change $u=x$, $v=ax+by$ yields
$u''=Au^3+(a+Bu^2)u'+Cu(u')^2$.
Writing $h(u)=e^{-Cu^2/2}$ and $z=\int_0^uh(s)ds$ gives a smooth ordinary Li\'enard equation on the full image interval with
$f(z(u))=-a-Bu^2$, $g(z(u))=-Ah(u)u^3$.
Every original periodic orbit is covered, since the coordinate change is invertible on the entire finite plane onto that strip. If
\[
 A<0,\quad a,C>0,\quad B=-D<0,\qquad
 (5-2\sqrt6)aC<D<(5+2\sqrt6)aC,
\]
then
\[
 \mathcal U_0(z(u))=-Au^2[DCu^4-(D+aC)u^2+3a]>0\quad(u\ne0).
\]
The quadratic in $u^2$ has positive leading coefficient and negative discriminant. The restoring force points to zero, so Theorem~\ref{cr:unique} gives at most one periodic orbit, which is hyperbolic and attracting if it exists. The explicit system $\dot x=x+y$, $\dot y=-2x^3+xy^2$ satisfies these conditions, with certificate $u^2(2u^4-3u^2+3)$. The argument gives an upper bound for this subfamily; existence and a two-cycle bound for the full family remain separate questions.

Gavrilov and Iliev \cite{GavrilovIliev2015} study quadratic perturbations of Hamiltonian two-saddle loops. One admissible perturbation direction is
\[
 \dot x=2xy,\qquad \dot y=-y^2-3x^2+3-\varepsilon(x-1)y,\qquad \varepsilon>0.
\]
The half-plane $x>0$ is invariant. The change $z=2\sqrt x$ is a diffeomorphism there, and its differentiated equation has
$f(z(x))=\varepsilon(x-1)$, $g(z(x))=6\sqrt x(x^2-1)$.
After translating $z=2$ to zero, the restoring force has the correct strict sign on the entire strip. Differentiating with respect to $z$, using $dx/dz=\sqrt x$, gives
\[
 \mathcal U_0(z(x))=-3\varepsilon(x-1)^2(3x+1)<0\quad(x>0,\ x\ne1).
\]
There is at most one periodic orbit in this half-plane, hyperbolic and repelling if present. The transformation degenerates on the boundary containing the two saddles, so the bound is confined to $x>0$. A bound for the general perturbation problem would also have to control the saddle boundary.

The same criterion gives another proof for a previously studied subfamily. The Li\'enard form in \cite[Section~3]{LuZhang2026} has
\[
 g=xe^{-x^2},\qquad f=(a+\alpha x^2)e^{-x^2/2}.
\]
For $\alpha>0$ and $-\alpha<a<0$,
$\mathcal U_0=e^{-3x^2/2}[-a+(a+\alpha)x^2+\alpha x^4]>0$.
Taking $F=\int_0^xf$ therefore gives at most one periodic orbit, which is hyperbolic and attracting if it exists. That paper already proves uniqueness for this subfamily; the calculation gives another proof of the upper bound. Existence requires a separate argument.

\subsubsection{Two existence questions for a quintic vector field}
The two existence questions of Llibre and Salhi \cite{LlibreSalhi2019} can be answered by an inverse integrating factor. This argument applies directly to their quintic vector field, which can be written as
\begin{equation}\label{ap:llibre-salhi}
 \dot x=xR-y\Omega,\quad \dot y=yR+x\Omega,\quad
 S=x^2+y^2,\quad E=x^2/a+y^2/b-1,
\end{equation}
where $ab\ne0$, $R=(1-S)[-\lambda E+2(a-b)xy]$ and
$\Omega=-2ab(SE-S+1)$.
For every $\lambda\ne0$, it has exactly one hyperbolic limit cycle if $a,b<0$, and exactly two if $a,b>1$. These cover the full respective ranges in Open Questions~1 and~2 of the author manuscript of \cite{LlibreSalhi2019}.

To prove this, differentiate the two defining polynomials:
\[
 \dot S=2S(S-1)[\lambda E-2(a-b)xy],\qquad
 \dot E=E[-2\lambda(E+1)(1-S)+4(a-b)xy(1-2S)].
\]
The polynomial $V=S(S-1)E$ satisfies $X\cdot\nabla V=(\operatorname{div}X)V$. A limit cycle must lie in $V=0$. Indeed, if $V\ne0$ along a periodic orbit, the vector field $X/V$ is divergence free in a tubular neighborhood. Its return map preserves a nonvanishing transverse density. Integrating that density gives a local coordinate in which the return map is a translation; a fixed point forces this translation to vanish identically. The orbit is then nonisolated. If an orbit meets $V=0$, the displayed scalar equation for $V$ forces it to remain there.

On the circle $S=1$, write $x=\cos\theta$, $y=\sin\theta$ and
$E_C=\cos^2\theta/a+\sin^2\theta/b-1$. In both parameter ranges $E_C<0$ and $ab>0$, so $\dot\theta=-2abE_C>0$: the entire circle is periodic. The normal variational coefficient is $K_C=2\lambda E_C-4(a-b)\cos\theta\sin\theta$, and
\[
 \frac{K_C}{\dot\theta}=-\frac\lambda{ab}+\frac{d}{d\theta}\log|E_C|.
\]
Its multiplier is $M_C=\exp(-2\pi\lambda/(ab))\ne1$.
When $a,b>1$, the ellipse $E=0$ encloses the circle and has $S>1$, $\dot\theta=2ab(S-1)>0$. Its normal coefficient is
$K_E=2\lambda(S-1)+4(a-b)xy(1-2S)$ and on it
$S_\theta=-2(a-b)Sxy/(ab)$. Therefore
\[
 \frac{K_E}{\dot\theta}=\frac\lambda{ab}
              +\frac{d}{d\theta}\log[S(S-1)],\qquad
 M_E=\exp(2\pi\lambda/(ab))\ne1.
\]
For $a,b<0$ the ellipse has no real points, so the circle is the only nontrivial compact component of $V=0$. For $a,b>1$, the circle and ellipse are its only such components. This proves both exact counts and hyperbolicity. For $\lambda>0$ the circle is attracting and the ellipse repelling; the directions reverse for $\lambda<0$.

\subsubsection{A square-root example with exactly one cycle}
For $a>0$ let $g(x)=\operatorname{sgn}(x)\sqrt{|x|}$ and
$F(x)=x^3-a\operatorname{sgn}(x)\sqrt{|x|}$, with both values zero at the origin. The coefficients are continuous but not locally Lipschitz there. Off zero,
\[
 \frac{F'}g=3x\sqrt{|x|}-\frac a{2x},\qquad
 \left(\frac{F'}g\right)'=\frac92\sqrt{|x|}+\frac a{2x^2}>0.
\]
Theorem~\ref{cr:unique} gives at most one periodic orbit and hyperbolic stability. Existence can be proved without appealing to a nonsmooth Poincare--Bendixson theorem. Set
$G=2|x|^{3/2}/3$, $E=G+y^2/2$. A small level of $E$ is outward because
$\dot E=|x|(a-|x|^{5/2})>0$ off its two transverse axis crossings. For an outer level use $c(x)=x/\sqrt{1+x^2}$ and $V=G+(y-c)^2/2$. Identity~\eqref{ap:trapping} remains valid using only continuity of $F,g$. Since $c'>0$, $g(F-c)>0$ at infinity and $c'(F-c)/g\to0$, a sufficiently large regular level of $V$ is inward, by exactly the two-region estimate following~\eqref{ap:trapping}. The resulting compact annulus avoids the equilibrium.

Uniqueness of trajectories in this annulus holds off the $y$-axis by smoothness, and on that axis by the graph equation $dy/dx=g/(F-y)$, whose denominator is bounded away from zero and whose dependence on $y$ is locally Lipschitz. A trajectory starting on the positive $y$-axis must return to the $y$-axis in finite time. Otherwise it stays on the right, where $y$ decreases to a finite limit and $\int_0^\infty g(x(t))dt<\infty$. Bounded $\dot x$ makes $g(x(t))$ uniformly continuous; a nonnegative integrable uniformly continuous function tends to zero, so $x(t)\to0$. The equation for $\dot x$ then forces $y(t)\to0$, contradicting exclusion of the equilibrium from the annulus. The next crossing is on the negative $y$-axis, and the reflected argument forces a subsequent positive-$y$-axis crossing. These transverse crossings define a continuous self-map of the compact positive-$y$-axis segment in the annulus. Endpoint inequalities and the intermediate-value theorem give a fixed point. Thus there is exactly one hyperbolic attracting cycle. Finally $xF(x)=x^4-a|x|^{3/2}\le0$ near zero, so Lemma~\ref{cr:flow} gives backward uniqueness at the equilibrium and excludes a nonconstant periodic solution through it. The trapping argument uses no differentiability of the flow at the equilibrium.

The simpler $C^1$ example $g=x$, $F(x)=\operatorname{sgn}(x)P(x^2/2)$ with
$P(u)=u-2u^2+3(u-1)_+^2$ has at most two limit cycles: its profile slopes are $\pm[1-4u+6(u-1)_+]$, respectively convex and concave. Theorem~\ref{cr:two} thus applies across a second-derivative jump. The conclusion here is an upper bound.

\subsubsection{A finite-growth Wolkowicz--Rothe--Shafer subfamily}
Consider
\begin{equation}\label{ap:wrs}
 \dot x=x(1-x/k)-\frac{xY}{x^2+ax+1},\qquad
 \dot Y=Y\left[-1+\frac{(a+5/2)x}{x^2+ax+1}\right],
\end{equation}
where $a>-2$, $1/2<k<19/10$. Put $s=1/2$, $d_* =\sqrt{769}$ and
\[
 \begin{gathered}
 a_*=(d_*-33)/32,\qquad k_*=(d_*+33)/32,\\
 H=(x^2+ax+1)(1-x/k),\qquad Q=5/2-x-1/x,\\
 f=-H'-Q,\qquad g=HQ.
 \end{gathered}
\]
There is a neighborhood of $(a_*,k_*)$ in which the conditions $A_0\ge0$, $A_2\le0$ imply at most two positive limit cycles, where
\[
 \begin{gathered}
 v=g'(s)>0,\quad w=g''(s),\quad t=g'''(s),\quad A_0=f(s)/\sqrt v,\\
 A_2=\frac{12v^2f''(s)-12vwf'(s)+(5w^2-3vt)f(s)}{24v^{7/2}}.
 \end{gathered}
\]
Every neighborhood of that point contains an open parameter set with exactly two hyperbolic cycles, inner repelling and outer attracting.

The parameter argument of Proposition~\ref{ap:ivlev-two} applies once the global sign condition and the local nondegeneracy are established. The positive time change $d\tau/dt=x/(x^2+ax+1)$ gives $x'=H-Y$, $Y'=YQ$, hence $x''+fx'+g=0$. The inverse is $Y=H-x'$. Every positive periodic orbit has $x<k$, and
$Q=-(x-1/2)(x-2)/x$ shows that $g$ has the correct single-well sign on $(0,k)$. At the base point $A_0=A_2=0$ and
\begin{equation}\label{ap:wrs-factor}
 \mathcal N(x)=\frac{42(x-1/2)^4}{k_*^3x^4}R(x),
\end{equation}
where
\begin{align*}
 R(x)={}&x^8+\frac{d_*-279}{56}x^7+
 \frac{50483-1027d_*}{7168}x^6+\frac{143793+13997d_*}{28672}x^5\\
 &-\frac{1548299+53663d_*}{57344}x^4+
 \frac{117163+3531d_*}{3584}x^3-\frac{426011+12507d_*}{28672}x^2\\
 &+\frac{103-9d_*}{896}x+\frac{929+33d_*}{512}.
\end{align*}
Using $2773/100<d_*<1387/50$, all degree-eight Bernstein coefficients of $R$ on the four intervals have the following rigorous common lower bounds:
\[
\begin{array}{c|c}
{[0,1/2]}&4501233/2867200\\
{[1/2,1]}&6514533999/5138022400\\
{[1,3/2]}&9630293/5734400\\
{[3/2,19/10]}&175522820427/1792000000000
\end{array}
\]
Thus $R>0$ on the entire enlarged interval. Local potential normalization at $s$ gives
$\widehat{\mathcal N}=3A_0-A_2z^2+3A_4z^4+O(z^5)$ with $c_*:=A_4(a_*,k_*)>0$, by~\eqref{ap:wrs-factor}. Near the center the same uniform Taylor bound as in Proposition~\ref{ap:ivlev-two} gives
$\widehat{\mathcal N}\ge3A_0-A_2z^2+c_*z^4$.
Away from the center, use the polynomial $k^3x^4\mathcal N$; its strict positive minimum on the remaining closed subset of $[0,19/10]$ persists. This controls both the singular left endpoint and the moving upper endpoint. Theorem~\ref{cr:two} supplies the global bound.

For independence of the parameters set $T=H'(s)$, $C=H'''(s)Q'(s)-H''(s)Q''(s)$. Direct differentiation yields
\[
 \det\frac{\partial(T,C)}{\partial(a,k)}
 =\frac{2(16ak+33k-5)}{k^3},\qquad a_*k_*=-5/16.
\]
This determinant is positive at the base point. Since $f(s)=-T$ and on $T=0$,
$vf''-wf'=-H(s)C$, the Jacobian of $(A_0,A_2)$ equals
$H(s)/(2v^3)$ times this determinant and is nonzero. Therefore we may prescribe
$A_0=c_*\epsilon^4/2$, $A_2=-5c_*\epsilon^2/2$.
The return-map calculation about the reversible system in Proposition~\ref{ap:ivlev-two} applies: the coefficients are analytic near a simple restoring zero, $(A_0,A_2)$ are independent, and $A_4=c_*+O(\epsilon^2)$. It gives the same two simple zeros at scaled radii $1,2$ and the stated stability. The cycles remain near $(s,H(s))$ and are positive. The global upper bound excludes additional cycles. Strict inequalities and hyperbolicity give open parameter sets, also after small variations of the original growth, mortality and denominator parameters, since the cleared numerator remains strictly positive and the unique simple restoring root persists.

Broer and Gaiko \cite[Theorem~4.3]{BG2010} give a global two-cycle bound for a model containing this subfamily; the later paper \cite{LZH2024} states a general WRS two-cycle conjecture. For this correspondence, set $M=a+5/2$, $X=Mx$, $Z=Y$ and use their parameters
$\alpha=M^{-2}$, $\beta=a/M$, $\delta=1$, $\lambda=(Mk)^{-1}$ and zero predator competition. Their denominator positivity condition is exactly $a>-2$. Thus the upper bound for this subfamily is already covered by their theorem. Earlier local and global parameter-region results also appear in \cite{ZCW2002,XZ2006}; the calculation above supplies an algebraic proof near the specified finite-growth parameter center.

\paragraph{Nonrational perturbations of the Selkov family}
Compactly supported perturbations of \eqref{sk:negative-lienard}
give nonrational systems to which Theorem~\ref{sk:two} applies.
\begin{corollary}\label{sk:perturbation}
Fix $1/9<p<p_*$ and choose $9/8<B_-<m(\sqrt p)$.
Translate \eqref{sk:negative-lienard} by $x=u-1$, keeping the notation
$F_b,g$, and restrict $b$ to $[\sqrt{B_-},\sqrt{B_H}]$.
For every fixed $\psi\in C_c^\infty((\sqrt p-1,\infty)\setminus\{0\})$
there exists $\epsilon_0>0$ such that
\[
 \dot x=y-F_b(x)-\epsilon\psi(x),\qquad \dot y=-g(x)
\]
has at most two periodic orbits in its natural strip for every
$|\epsilon|<\epsilon_0$ and every parameter in that interval.
All conclusions about multiplicity and stability in
Theorem~\ref{sk:two} remain valid.
\end{corollary}
\begin{proof}
The Selkov coefficients satisfy
the three strict curvature inequalities for every $0<B\le B_H$.
At $B_-$ the profile difference is strictly negative at all positive
heights. Since $b>1$ and $B_H<2$, the focus inequality also holds
on this slightly enlarged closed parameter band.
The perturbation vanishes near zero and near both strip endpoints,
so all local coefficients, endpoint behavior and the Hopf parameter
are preserved. It is independent of $b$, so $x\partial_bF_b>0$ is unchanged.
On the compact support times the closed parameter band, right
$A>0$, left $f<0$, left $\Xi>0$, and the upper-endpoint inequality
$(F_HL_H)_x<0$ have uniform strict margins.
To control the region $f>0$ under perturbation, note that at an
original zero of $f$ on the support, $\Xi=-A^2<0$. Hence the curvature
expression remains uniformly negative on the original closed positive
domain and on a neighborhood of its boundary. Every new positive-$f$ point lies in
that neighborhood when the perturbation is small. Preserving right
$A>0$ and the signs near zero and infinity also preserves exactly
one damping zero. Finally the affected potential heights form a
compact subinterval of $(0,\infty)$, where the negative
lower-boundary profile difference is bounded away from zero. A sufficiently
small $\epsilon_0$ preserves all of (S1)--(S4), and the theorem applies.
\end{proof}
The allowed size depends on the fixed choices of $p,B_-$ and $\psi$;
the estimates need not remain uniform as these approach their
boundary values. A nonzero compactly supported smooth perturbation
is generally nonrational.
If the unperturbed member has the two hyperbolic cycles constructed
above, sufficiently small perturbations preserve them, so this class
also contains systems attaining the upper bound.

\subsubsection{A joint criterion with sign-changing third derivatives}
\label{jt:example}
Let
\[
 p(z)=4z-\frac{20}{3}z^3+\frac85z^5+\cos(8z)-1,\qquad
 h(x)=\frac{33(e^x-1)}{11e^x+12}.
\]
For every $0<\varepsilon\le10^{-210}$ the full-plane system
\begin{equation}\label{jt:example-system}
 \dot x=y-\varepsilon p(h(x)),\qquad \dot y=-h(x)h'(x)
\end{equation}
has exactly two limit cycles. Both are hyperbolic; the inner one repels and
the outer one attracts. The third derivatives of both profiles change sign, but their combined contribution to the displacement is estimated by Theorem~\ref{jt:weak}.

Indeed, $h'(x)=759e^x/(11e^x+12)^2>0$, and $h$ maps $\mathbb R$
diffeomorphically onto $(-11/4,3)$. Its restoring force has
$G=h^2/2$, with $U_-=121/32$ and $U_+=9/2$. The positive change
$z=h(x)$, $d\tau=h'(x)\,dt$ gives the normalized strip system of Theorem~\ref{jt:weak}, with a one-to-one correspondence between compact cycles.
For $|w|<10$,
\[
 |p(w)/w|\le4+\frac{2000}{3}+16000+8e^{80}<10^{40}.
\]
Here $(\cos(8w)-1)/w=-8\int_0^1\sin(8tw)\,dt$ supplies the bound
at zero as well. For $0<z\le1/16$,
\[
 p(z)/z\ge4-\frac{20}{3\cdot256}-32z>1,
\]
using $\cos(8z)-1\ge-32z^2$. For negative $z$ the even term divided
by $z$ is nonnegative and the same positivity follows.
The standard integrals of $\sin^2\theta$, $\sin^4\theta$ and
$\sin^6\theta$ give
\[
 D_p(r)=\frac{\pi r}{2}(4-5r^2+r^4),\qquad
 C_p'(r)=\frac{15\pi}{4}\left(3+\frac5{r^2}+\frac{12}{r^4}\right)>30.
\]
The theorem therefore gives at most two cycles in the plane for~\eqref{jt:example-system}.

To prove existence, fix a real $1/16\le r\le4$.
The midpoint bound \eqref{jt:data} gives
$|\Delta(r,\varepsilon)/\varepsilon|\le64\cdot10^{40}$ on
$|\varepsilon|<E=10^{-100}$. The quotient is holomorphic at zero.
Cauchy's coefficient estimates, summed as a geometric series for
$|\varepsilon|\le E/2$, therefore give
\[
 |\Delta(r,\varepsilon)-\varepsilon D_p(r)|
 \le128\cdot10^{140}\varepsilon^2<10^{143}\varepsilon^2.
\]
At $r=1/2,3/2,5/2$ the respective values of $D_p$ are
$45\pi/64,-105\pi/64,945\pi/64$. The error is smaller than each
absolute linear term, giving the signs $+,-,+$.
For all $1/2\le r\le5/2$, the polar estimates give
$|\rho-1|\le24\cdot10^{40}\varepsilon$ and hence
$|z|<2.501<11/4$. Both half arcs stay in the physical strip throughout
this interval. The intermediate value theorem yields two matching roots.
By \eqref{jt:scalar} and strict $\mathcal Q'>0$, the transformed
displacement is first decreasing and then increasing; the two roots lie
on its strict branches. Thus $\Delta'<0$ at the inner root and
$\Delta'>0$ at the outer root. The return multiplier in
Proposition~\ref{cr:matching} gives hyperbolicity and the stated stability.

Finally, in the normalized variable the third-order expression, divided
by $\varepsilon$, is
\begin{align*}
 N_p(z)&=z^2p'''(z)-3zp''(z)+3p'(z)\\
 &=12+20z^2+24z^4+(512z^2-24)\sin(8z)+192z\cos(8z).
\end{align*}
It is positive at zero, whereas
\[
 N_p(\pi/8)=12+\frac{5\pi^2}{16}+\frac{3\pi^4}{512}-24\pi<0,
\]
\[
 N_p(-\pi/4)=12+\frac{5\pi^2}{4}+\frac{3\pi^4}{32}-48\pi<0.
\]
The bounds $3<\pi<4$ suffice for both negative signs. The pullback
identity is $\mathcal N(x)=\varepsilon h'(x)^5N_p(h(x))$.
Since $h'>0$ and $\phi_\pm'''=\mathcal N(x_\pm)/g(x_\pm)^5$,
each separate profile third derivative changes sign. Thus neither profile satisfies the separate third-order sign condition; the upper bound follows from the joint criterion and its estimate of the full nonlinear remainder.

\subsubsection{Families outside the sign conditions}
The sign conditions do not apply to every member of the polynomial families considered above. For instance, the JMAA family with $m=1,n=5,\xi=1,\beta=2$ has $\mathcal N(-1)=12$ and $\mathcal N(-2)=-23994$. The whole left well fails the displayed two-cycle sign test; its Conjecture~7.1 is not proved here. The Chen--Han--Xia family $F=ax+bx^3+x^5$, $g=cx+x^3$, $a,c>0$ has $\mathcal N(0)=3ac^2>0$ but $\mathcal N(x)=-15x^8+O(x^6)$ at infinity. The two-cycle conclusions of \cite{ChenHanXia} for that family therefore require a different criterion.

The sextic family $F=\alpha x+x^3+x^6$, $g=x$, has
$\mathcal N=3\alpha-3x^2+48x^5$, whose two tails have opposite signs. Hence the common-sign condition fails on the whole plane. Uniform bounds for general sextic and septic Li\'enard polynomials remain beyond these criteria.

\bibliographystyle{elsarticle-num-alpha}
\bibliography{references}
\end{document}